\documentclass[10pt, leqno]{amsart}
\usepackage[margin=1in]{geometry}
\usepackage{amsmath,amssymb,amsthm,mathtools}
\usepackage{comment}
\usepackage[colorlinks=true,linkcolor=blue,citecolor=blue,urlcolor=blue]{hyperref}

\numberwithin{equation}{section}
\newtheorem{theorem}{Theorem}[section]
\newtheorem{lemma}[theorem]{Lemma}
\newtheorem{proposition}[theorem]{Proposition}
\newtheorem{corollary}[theorem]{Corollary}

\theoremstyle{definition} \newtheorem{definition}[theorem]{Definition}

\newcommand{\ph}{\varphi}
\newcommand{\ten}{\otimes}
\newcommand{\id}{\mathrm{id}}
\newcommand{\norm}[1]{\left\|#1\right\|}
\newcommand{\ip}[2]{\left\langle#1,#2\right\rangle}

\DeclareMathOperator{\ran}{ran}
\DeclareMathOperator{\dom}{dom}
\title{Selflessness of $q$-Araki-Woods $C^\ast$-algebras}

\author[Kumar]{Manish Kumar}  
\author[Lu]{ Xiao-Qi Lu}
\author[Wang]{Simeng Wang} 

\address[Manish Kumar]{KU Leuven, Department of Mathematics, Celestijnenlaan 200B, 3001 Leuven, Belgium}
\email{manish.kumar@kuleuven.be}
\address[Xiao-Qi Lu]{School of Mathematics and Statistics, University of Glasgow, University Place, Glasgow G12 8QQ, UK}
\email{xiaoqi.lu@glasgow.ac.uk}
\address[Simeng Wang]{Institute for Advanced Study in Mathematics, Harbin Institute of Technology, Harbin, 150001, China}
\email{simeng.wang@hit.edu.cn}



\begin{document}
	
	\begin{abstract}
		We prove that the $q$-Araki-Woods $C^*$-algebra associated with $q\in(-1,1)$ and a strongly continuous one-parameter group $U$ of orthogonal transformations on a  real Hilbert space $\mathcal H_{\mathbb R}$ is completely selfless whenever $\dim_{\mathbb R} \mathcal H_{\mathbb R}\ge 2$. In particular, the $C^\ast$-algebra is simple, and it is purely infinite when $U$ is nontrivial. This settles the simplicity question in full generality. The proof begins in finite dimension, where finitely many modular eigenoperators generate the algebra. For $C^\ast$-algebras with such generators, we use free Fisher information to obtain a commutator spectral gap in the centralizer and combine it with a rapid decay estimate  to establish a general simplicity criterion. 
	\end{abstract}
	
	\maketitle

	\section{Introduction}
	The $q$-commutation relation interpolates between the bosonic and
	fermionic relations. It appeared in Frisch and Bourret's 1970 work
	on stochastic processes and later in Greenberg's study of quon
	statistics \cite{FB70,Gre91}. Bo\.{z}ejko and
	Speicher realized the relation on a deformed Fock space and used
	it to construct a generalized Brownian motion \cite{BS91}.
	This led to the theory of $q$-Gaussian processes developed by
	Bo\.{z}ejko, K\"ummerer, and Speicher \cite{BKS97}.  For a  real separable Hilbert space $\mathcal{H}_{\mathbb R}$ and $q\in (-1,1)$, the {\em $q$-Gaussian $C^\ast$-algebra} $A_q(\mathcal{H}_{\mathbb R})$ is generated by a family of the so-called $q$-Gaussian operators $\{\ell_i+\ell_i^\ast\}_{1\leq i\leq \dim\mathcal{H}_{\mathbb R}}$, where $\ell_i$ is the left creation operator acting on a $q$-deformed Fock space  satisfying the $q$-commutation relations: $\ell_i^\ast\ell_j-q\ell_j\ell_i^\ast=\delta_{i,j}$. At $q=0$, this recovers Voiculescu's free Gaussian construction \cite{Voi90}.
	
	The structure of $q$-Gaussian $C^\ast$-algebras and the associated von Neumann algebras has been studied extensively over several decades. The case $q=0$
	is the natural starting point. In this case, $A_0(\mathcal H_{\mathbb R})$ is the reduced free
	product of copies of $C([-2,2])$, each equipped with the
	semicircular state \cite{Voi90}. When $\dim_{\mathbb R}\mathcal H_{\mathbb R}\geq 2$,
	Avitzour's theorem gives simplicity and uniqueness of trace,
	and the theorem of Dykema, Haagerup, and R{\o}rdam gives stable
	rank one \cite{Avi82,DHR97}. 
	This analysis becomes considerably more delicate when $q\neq 0$. Bo\.{z}ejko, K\"ummerer, and Speicher proved
	simplicity when $\mathcal H_{\mathbb R}$ is infinite-dimensional
	\cite{BKS97}. Only recently, Amrutam, Jekel, and Wasilewski settled
	the finite-dimensional case for all $q\in (-1,1)$, 
	where they combined a commutator
	spectral gap with rapid decay to obtain the Dixmier averaging
	property, and hence simplicity and uniqueness of trace
	\cite{AJW26}. 
	Klisse and Pagliero \cite{KP26}  subsequently proved
	that $(A_q(\mathcal H_{\mathbb R}),\tau_q)$ is 
	selfless in the sense of Robert \cite{Rob25}, where $\tau_q$ is the vacuum trace. In particular, the algebra has strict
	comparison and stable rank one \cite{Rob25}. 
	
	In this paper, we study simplicity, Dixmier averaging, and selflessness for the
	nontracial analogues of the $q$-Gaussian $C^*$-algebras.
	In 1997, Shlyakhtenko extended Voiculescu's free Gaussian functor
	to the nontracial setting \cite{Shl97}. His construction starts with
	a real Hilbert space $\mathcal H_{\mathbb R}$ and a strongly
	continuous orthogonal representation $U=(U_t)_{t\in\mathbb R}$
	on it. The representation $U$ determines a new inner product on the
	complexification of $\mathcal H_{\mathbb R}$, giving rise to
	the Hilbert space $\mathcal H_U$ upon completion. This
	representation governs the modular group of automorphisms of the vacuum state. Hiai used the $q$-deformed Fock space over
	$\mathcal H_U$ to construct the $q$-Araki-Woods algebras
	\cite{Hia03}.  We write $A_q(\mathcal H_{\mathbb R},U)$ for the unital
	$C^\ast$-algebra generated on this Fock space by the field operators
	$\ell_q(\xi)+\ell_q(\xi)^\ast$, $\xi\in\mathcal H_{\mathbb R}$, and
	$\varphi_q$ for its vacuum state
	(see Section \ref{sec:preliminaries}).
	When $q=0$, this is Shlyakhtenko's free Araki-Woods $C^*$-algebra.
	When $U$ is trivial, it is the $q$-Gaussian $C^*$-algebra with its
	vacuum trace.
	
	Relatively little is known about the structure of $A_q(\mathcal H_{\mathbb R},U)$  for nontrivial $U$.  In the case $q=0$, Shlyakhtenko used reduced free product
	decompositions and Avitzour's averaging argument to prove simplicity
	in dimension at least three and for certain two-dimensional
	representations \cite{Shl97}. For $q=0$ and nontrivial $U$ in
	dimension at least three, Hayes, Kunnawalkam Elayavalli, and Robert
	proved selflessness and pure infiniteness
	\cite[Corollary~7.7]{HKER25}. Klisse and Pagliero subsequently
	proved complete selflessness for all free Araki-Woods
	$C^*$-algebras in dimension at least two \cite{KP26}.
	For $q\neq 0$ and nontrivial $U$, study has been even more limited and earlier results covered two
	regimes. In finite dimension, Nelson showed that
	$A_q(\mathcal H_{\mathbb R},U)$ is isomorphic to
	$A_0(\mathcal H_{\mathbb R},U)$ for sufficiently small $|q|$
	\cite{Nel15}. In infinite dimension,
	Mandal, Mukherjee, and Patri proved simplicity  when the positive generator $A_U$, determined by
	$U_t=A_U^{it}$, is bounded
	\cite{MMP25}.
	
	We settle here the simplicity question in full generality for all
	$q\in(-1,1)$ and all  representations
	$U$ on a real Hilbert space $\mathcal H_{\mathbb R}$ with
	$\dim_{\mathbb R}\mathcal H_{\mathbb R}\geq 2$. 
	In fact, we prove a
	much stronger result, namely that $A_q(\mathcal H_{\mathbb R},U)$ is
	completely selfless.
	Before stating
	the main theorem, we recall the notion of selflessness. Let $(\mathcal A,\varphi)$ be  a $C^\ast$-probability space, that is, $\mathcal{A}$  is a unital $C^\ast$-algebra and $\varphi$ is a state which induces a   faithful GNS representation. For a cofinal ultrafilter
	$\omega$ on an index set $\Lambda$, write $\ell^\infty(\Lambda, \mathcal{A})$ for the uniformly bounded elements in $\mathcal{A}$ indexed by $\Lambda$, and define the {\em $C^\ast$-algebra ultrapower} $\mathcal{A}^\omega$ with state $\varphi^\omega$ by 
	\[
	\mathcal A^\omega
	=
	\ell^\infty(\Lambda,\mathcal A)\big/
	\{(a_i)_{i\in \Lambda}:\lim_{i\to \omega}\|a_i\|=0\},
	\qquad
	\varphi^\omega((a_i)^\omega)=\lim_{i\to \omega}\varphi(a_i)
	\]
	where $(a_i)^\omega$ denotes the element in the quotient for $(a_i)\in \ell^\infty(\Lambda,\mathcal{A})$.   The diagonal embedding $\Delta_{\mathcal A}\colon \mathcal A\to\mathcal A^\omega$ sends
	$a$ to the class of the constant family $(a)_{i\in \Lambda}$. Following Robert \cite{Rob25},
	the pair $(\mathcal A,\varphi)$ is called \emph{selfless} if $\mathcal A\neq\mathbb C$ and,
	for some cofinal ultrafilter $\omega$, there is a state-preserving unital
	$*$-embedding
	$\Theta\colon (\mathcal A,\varphi)*(\mathcal A,\varphi)\longrightarrow(\mathcal A^\omega,\varphi^\omega)$
	whose restriction to the first factor is $\Delta_{\mathcal A}$. Here $*$ denotes
	the reduced free product and state preservation means
	$\varphi^\omega\circ\Theta=\varphi*\varphi$. In particular, the image of
	the second factor is free from $\Delta_{\mathcal A}(\mathcal A)$ with respect to
	$\varphi^\omega$. 
	Ozawa's notion of \emph{complete selflessness} strengthens selflessness by requiring the witnessing embedding $\Theta$ to induce an embedding $(\mathcal A*\mathcal A)\otimes_{\min}\mathcal D\to (\mathcal A\otimes_{\min}\mathcal D)^\omega$ for every $C^*$-algebra $\mathcal D$ \cite{Oza25}. 
	
	As shown in \cite{Rob25}, selfless $C^\ast$-algebras demonstrate many regularity properties.
	In particular, if $(\mathcal A,\varphi)$ is selfless, then $\mathcal A$ is
	simple and has the {\em Dixmier averaging property}: for every $a\in \mathcal{A}$,
	\[
	\varphi(a)1_{\mathcal A}
	\in
	\overline{\operatorname{conv}}^{\|\cdot\|}
	\{uau^*:u\in\mathcal U(\mathcal A)\}.
	\] 
	If $\varphi$ is tracial, then it is the unique tracial state on
	$\mathcal A$, and $\mathcal A$ has strict comparison of positive elements by
	$\varphi$ and stable rank one \cite{Rob25}. If $\varphi$ is 
	nontracial, then $\mathcal A$ is purely infinite \cite{Rob25, Gou26}, that is, for every non-zero positive element $a\in\mathcal A$, there exists  $x\in\mathcal A$ such that $x^\ast ax=1$. After the breakthrough paper \cite{AGKEP25}, which established that many reduced group $C^{\ast }$-algebras with the rapid decay property are selfless, numerous further examples have been discovered; see, for instance, \cite{Oza25, FKOP26, HKPR27, KP26}. We add a new example to this growing list and present the following main result of this paper.

	\begin{theorem}\label{main-theorem}
		Let $q\in (-1,1)$, and let $U$ be a  strongly continuous orthogonal representation of $\mathbb R$ on a real Hilbert space $\mathcal H_{\mathbb R}$. If $\dim_{\mathbb R}  \mathcal H_{\mathbb R}\geq2$, then $(A_q(\mathcal H_{\mathbb R},U),\ph_q)$ is completely selfless. In particular, if $U$ is nontrivial, then $A_q(\mathcal H_{\mathbb R},U)$ is simple and purely infinite.
	\end{theorem}
	Pure infiniteness is a $C^*$-algebraic analogue of the infiniteness
	of type~$\mathrm{III}$ factors. 
	Following partial factoriality
	results of Hiai \cite{Hia03} and of Bikram, Mukherjee, Ricard, and
	Wang \cite{BMRW23}, Kumar, Skalski, and Wasilewski settled the
	factoriality question for $q$-Araki-Woods von Neumann algebras in
	full generality, and in particular proved that the associated von Neumann
	algebra $A_q(\mathcal H_{\mathbb R},U)'' $ is a type~$\mathrm{III}$ factor whenever $U$ is nontrivial \cite{KSW23}.
	Theorem~\ref{main-theorem} above gives a parallel result for the
	$C^\ast$-algebra $A_q(\mathcal H_{\mathbb R},U)$, which is simple and
	purely infinite in this case.

	We now outline the proof strategy for Theorem~\ref{main-theorem}. The idea is to first prove the theorem for $A_q(\mathbb R^2,U)$ when $U$ is nontrivial on a two-dimensional space $\mathbb R^2$. The first step  is to
	establish Dixmier averaging with respect to the vacuum state.
	We obtain this by developing   a  simplicity and averaging criterion  for general $C^\ast$-algebras generated by modular eigenoperators, which we describe  below.
	
	Let $M$ be a von Neumann algebra with a faithful normal state
	$\varphi$, and let $\mathcal A\subseteq M$ be a unital $C^*$-algebra
	generated by eigenoperators of $\sigma^\varphi$ (see Definition \ref{def:eigenoperator} and Section \ref{sec:derivations on free fisher infomration}).  Write
	$\mathcal A^\varphi=\mathcal A\cap M^\varphi$ for the centralizer.  We first show that the Dixmier averaging property for $\mathcal A$ follows from the corresponding property for its centralizer $\mathcal A^\varphi$.
	
	\begin{theorem}[Theorem \ref{thm:sim_centra_implies_simple}]\label{thm:centralizer simplicity implies simplicity for introduction}
		With $M$, $\varphi$, and $\mathcal A$ as above, suppose that $\mathcal A^\varphi$ has the Dixmier averaging property with respect to $\varphi|_{\mathcal A^\varphi}$.
		Then $\mathcal A$ has the Dixmier averaging property with respect to
		$\varphi$, with the averages taken over unitaries in $\mathcal A^\varphi$.
		Similarly, if $\mathcal A^\varphi$ is simple, then $\mathcal A$ is simple.
	\end{theorem}
	
	To obtain Dixmier averaging in the centralizer $\mathcal A^\varphi$, we
	combine a commutator spectral gap with a rapid decay estimate,
	following the tracial approach of Amrutam, Jekel, and Wasilewski
	\cite{AJW26}.  The spectral gap gives $L^2$-decay for
	suitable unitary averages, and the rapid decay estimate lets us
	pass from $L^2$-decay to convergence in operator norm. 
	
	The spectral gap argument has its origins in work on von Neumann
	algebras. Dabrowski used commutator estimates to show that a tracial von
	Neumann algebra generated by at least two self-adjoint elements
	with finite non-microstates free Fisher information has no
	property $\Gamma$ \cite{Dab10}. Nelson extended this approach to
	modular eigenoperators and proved that the centralizer is a
	$\mathrm{II}_1$ factor without property $\Gamma$
	\cite{Nel17}.
	In
	Section~\ref{sec:derivations on free fisher infomration}, we combine
	Nelson's derivations \cite{Nel17} associated  to eigenoperator 
	with Lipschitz
	regularity of the conjugate variables (see Definition \ref{def:conjugate_variable}), to obtain
	the following estimate.

	\begin{proposition}[Proposition \ref{prop:fullgap}]
		With $M$, $\varphi$, and $\mathcal A$ as above, suppose that $M=W^*(G)$
		and $\mathcal A=C^*(1,G)$, where $G=G^*$ is a finite set of modular
		eigenoperators with finite free Fisher information and Lipschitz
		conjugate variables. Let $c\in G$ have modular eigenvalue
		$0<\lambda<1$, and put $\tau=\varphi_{|_{M^\varphi}}$. Then there
		are self-adjoint polynomials
		$s_1,s_2\in\mathbb C\langle c^*c,cc^*\rangle$ and a constant
		$C>0$ such that
		\[
		\norm{f-\tau(f)1}_2
		\leq C\sum_{j=1}^2\norm{[f,s_j]}_2,
		\qquad f\in L^2(M^\varphi,\tau).
		\]
	\end{proposition}
	Our passage from this spectral gap to norm averaging is motivated
	by  \cite[Theorem A]{AJW26}. They use a spectral gap to show that repeated
	unitary averaging converges exponentially fast in $L^2$. Proposition~\ref{prop:rapid_decay-property} adapts this argument to $\mathcal A^\varphi$, where the rapid decay hypothesis
	concerns the spaces obtained by multiplying a fixed  word with letter in $G$ with modular eigenvalue $1$
	on the left and right by words in $s_1$ and $s_2$. Together, the two hypotheses yield the Dixmier averaging property and uniqueness of trace for $\mathcal A^\varphi$. 

	In Sections~\ref{sec:operator_norm_estimate}-\ref{section:simplicity-q-Araki-Woods}, we apply our general results to $A_q(\mathbb R^2,U)$
	with $U$ nontrivial on a two-dimensional space $\mathbb R^2$.
	We  verify this rapid decay hypothesis for $A_q(\mathbb R^2,U)$ 
	using an operator-norm estimate in Proposition~\ref{prop:rd}.
	The first author's work on Lipschitz conjugate variables
	\cite{Kum23} supplies the regularity needed for the spectral gap.
	Together, these two ingredients give Dixmier averaging for
	$A_q(\mathbb R^2,U)$ and its  centralizer in
	Proposition~\ref{thm:simplicity-rotation-plane}. 
	
	The next step is to elevate the Dixmier averaging argument to prove selflessness for $A_q(\mathbb R^2,U)$. In fact, we prove  a stronger fact that $A_q(\mathbb R^2,U)$ is Toeplitz selfless in the sense of \cite{GJKEPR26}.
	The calculations for this in Section \ref{sec:selfless_q-Araki-Woods} are motivated by the first and third authors' proof of fullness of $q$-Araki-Woods factors \cite{KW24}. 
	Their estimates for powers of Wick operators on a two-dimensional plane
	guide $L^2$-norm calculations in our proof. As in  \cite{KW24}, we make crucial use of the nontriviality of $(U_t)_{t\in\mathbb R}$.
	First, we fix  unit
	eigenvectors $e$ and $f$ in $\mathbb C^2$ of the positive invertible generator $A_U$ of $U_t$ with eigenvalues $\lambda^{-1}$ and $\lambda$ for some $\lambda<1$, respectively. In Lemma \ref{lemma;Positive element B}, we construct the positive
	centralizer element
	$B=\lim_{n\to\infty}(1-q)^nW(e)^{*n}W(e)^n$
	where the limit is in operator norm and $\varphi_q(B)>0$. Here $W(e)$ is the Wick operator (see Section \ref{sec:preliminaries}).
	We establish in  Proposition \ref{prop:estimate for Toeplitz selflessness} and Corollary \ref{cor: rest-A-Toe-selfless} 
	that the completely positive map
	$\Phi_B: x\to \varphi_q(B^{1/2}xB^{1/2})B$ is Toeplitz selfless. 
	Dixmier averaging property then lets us write $\varphi_q(\cdot)1$ by  convex combinations of conjugates of the map $\Phi_B$, which will then yield the required Toeplitz selflessness of $A_q(\mathbb R^2,U)$ by the permanence
	properties in \cite{GJKEPR26}.
	For a general orthogonal representation $U$ on a real Hilbert space $\mathcal{H}_{\mathbb R}$, we extend these calculations  first to the finite-dimensional case and then invoke  another permanence theorem \cite{Rob25, GJKEPR26} for directed
	norm closures of selfless $C^\ast$-algebras to  finish the proof of Theorem~\ref{main-theorem}.
	
	The outline of the paper is as follows. In Section \ref{sec:preliminaries}, we recall the construction of $q$-Araki-Woods $C^\ast$-algebras and a few of their properties. In Section \ref{sec:derivations on free fisher infomration}, we study Nelson's derivations on the centralizer of von Neumann algebras generated by eigenoperators and prove a spectral gap property. Section \ref{sec:simplicity of C-algebras with eigenoperators} obtains a Dixmier averaging property for such centralizer using a rapid decay property and use it to prove simplicity for the $C^\ast$-algebras. In Section \ref{sec:operator_norm_estimate}, we obtain the required rapid decay property for $A_q(\mathcal{H}_\mathbb{R}, U)$ in the two-dimensional setup, and use it to prove its simplicity in Section \ref{section:simplicity-q-Araki-Woods}. Finally, Section \ref{sec:selfless_q-Araki-Woods} studies Toeplitz selflessness in the nontracial case and proves Theorem \ref{main-theorem}.

	\section{Preliminaries}\label{sec:preliminaries}
	
	In this section, we recall some well-known facts about $q$-Araki-Woods $C^\ast$ and von Neumann algebras introduced by Hiai \cite{Hia03}. As a convention, 
	inner products are linear in the second variable.

	Let $\mathcal H_{\mathbb R}$ be a real Hilbert space with $\dim_{\mathbb R} \mathcal H_{\mathbb R}\ge2$, and let $(U_t)_{t\in\mathbb{R}}$ be a  strongly continuous orthogonal representation of $\mathbb R$ on $\mathcal H_{\mathbb R}$. Consider the complexification $ \mathcal H_{\mathbb C}= \mathcal H_{\mathbb R}\otimes_{\mathbb R}\mathbb C$ of $\mathcal H_{\mathbb R}$ and denote its inner product by $\langle\cdot,\cdot\rangle_{\mathcal H_{\mathbb C}}$. Identify $\mathcal H_{\mathbb R}$ inside $ \mathcal H_{\mathbb C}$ by $\mathcal H_{\mathbb R}\otimes 1$, so that $  \mathcal H_{\mathbb C}=\mathcal H_{\mathbb R}+i\mathcal H_{\mathbb R}$.
	Let $I:\mathcal H_{\mathbb C}\to \mathcal H_{\mathbb C}$ denote the antilinear operator given by $I(\xi+i\eta)=\xi-i\eta$, for $\xi,\eta\in \mathcal H_{\mathbb R}$.
	Extend $(U_t)_{t\in \mathbb R}$ to a strongly continuous one-parameter group of unitaries on $\mathcal H_{\mathbb C}$, which is still denoted by $(U_t)_{t\in\mathbb{R}}$.
	Let $A_U:\mathrm{ dom \;} A_U\subset \mathcal H_{\mathbb C}\to H_{\mathbb C}$ be the positive, self-adjoint, and non-singular (unbounded) operator for which $U_t=A_U^{it}$ for all $t\in \mathbb R$. Note that $IA_U = A_U^{-1}I$.
	As in \cite{Shl97}, define the deformed inner product $\ip\cdot\cdot_U$ on $\mathcal H_{\mathbb C}$ by
	\[
	\ip\xi\eta_U=\ip\xi{2A_U(1+A_U)^{-1}\eta}_{\mathcal H_{\mathbb C}}\quad \forall \xi, \eta\in \mathcal H_{\mathbb C}.
	\]
	Denote by $\mathcal H_U$ the Hilbert space completed by
	$H_{\mathbb C}$ with respect to the inner product $\ip{\cdot}{\cdot}_U$.
	We still denote the inner product of $\mathcal H_U$ by $\ip{\cdot}{\cdot}_U$. Then there is an isometric embedding of $\mathcal H_{\mathbb R}$ into $\mathcal H_U$ such that $\mathcal H_{\mathbb R} +i\mathcal H_{\mathbb R}$ is dense in $\mathcal H_U$ and $\mathcal H_{\mathbb R}\cap i \mathcal H_{\mathbb R} = \{0\}$ (see \cite{Shl97}).
	
	For $q\in (-1,1)$ we form   the {\em $q$-Fock space} $\mathcal F_q(\mathcal H_U)$ introduced in \cite{BS91}: for each $n\geq 1$, consider the map $P_q^{(n)}:\mathcal{H}_U^{\otimes{n}}\to \mathcal{H}_U^{\otimes n}$ on the algebraic tensor product given by 
	\[P_q^{(n)}=\sum_{\sigma\in S_n}q^{i(\sigma)}U_\sigma,\qquad U_\sigma(\xi_1\otimes\ldots\otimes\xi_n)=\xi_{\sigma^{-1}(1)}\otimes\ldots\otimes\xi_{\sigma^{-1}(n)},\qquad\xi_j\in \mathcal{H}_U\]
	where $S_n$ denotes the symmetric group on the set $[n]:=\{1,\ldots,n\}$ and $i(\sigma)=|\{(i,j)\in [n]^2: i<j, \sigma(i)>\sigma(j)\}|$ denotes the number of inversions  for the permutation $\sigma\in S_n$. By \cite[Proposition 1]{BS91},  $P_q^{(n)}$ defines a strictly positive operator on $\mathcal{H}_U^{\otimes n}$; hence it defines an inner product on $\mathcal{H}_U^{\otimes n}$ by 
	\[\langle\xi,\eta\rangle_q=\langle\xi, P_q^{(n)}\eta\rangle_U,\qquad\xi,\eta\in \mathcal{H}_U^{\otimes n}.\]
	Denote by $\mathcal{H}_U^{\otimes_q^n}$ the completion of $\mathcal{H}_U^{\otimes n}$ with respect to this inner product and define the {\em $q$-Fock space} by 
	\begin{equation*}
		\mathcal F_q(\mathcal H_U)=\mathbb C\Omega\oplus\bigoplus_{n\geq 1}\mathcal H_U^{\otimes_q^n},
	\end{equation*}
	where $\Omega$ is a distinguished unit vector in $\mathbb C$ (called {\em vacuum vector}). Remark that the above construction applies to any real or complex Hilbert space and $q\in (-1,1)$. 
	
	For each $\xi\in \mathcal H_U$, define the {\em left q-creation} and {\em left q-annihilation} operators on $\mathcal F_q(\mathcal H_U)$ by
	\begin{align*}
		&\ell_q(\xi)\Omega=\xi,\qquad\ell_q(\xi)(\xi_1\otimes\ldots\otimes\xi_n)=\xi\otimes\xi_1\otimes\ldots\otimes\xi_n,\;\;\;\mbox{ and }\\ 
		&\ell_q^*(\xi)\Omega=0,\qquad\ell_q^*(\xi)(\xi_1\otimes\ldots\otimes\xi_n)=\sum_{i=1}^nq^{i-1}\langle\xi,\xi_i\rangle_U\xi_1\otimes\ldots\otimes\xi_{i-1}\otimes\xi_{i+1}\otimes\ldots\otimes\xi_n, 
	\end{align*}
	for all $\xi_1\otimes\ldots\otimes\xi_n\in\mathcal H_U^{\otimes n}$ and $n\ge 1$. 
	Then $\ell_q(\xi)$ and $\ell_q^*(\xi)$ extend to bounded linear operators on $\mathcal F_q(\mathcal H_U)$ and are adjoint to each other \cite{BS91}.
	Moreover, they satisfy the following $q$-commutation relation:
	\begin{equation*}
		\ell_q^*(\xi)\ell_q(\eta)-q\ell_q(\eta)\ell_q^*(\xi)=\langle\xi,\eta\rangle_U,\qquad \xi,\eta\in \mathcal H_U.
	\end{equation*}
	We now consider the $q$-Araki-Woods $C^*$-algebras and von Neumann algebras associated with $U$ arising from Hiai's construction \cite{Hia03}.
	Define the real-linear map
	$s_q:\mathcal H_{\mathbb R}\to B(\mathcal F_q(\mathcal H_U))$ by
	$s_q(\xi)=\ell_q(\xi)+\ell^*_q(\xi)$ for $\xi\in \mathcal H_{\mathbb R}$.
	The {\em $q$-Araki-Woods} $C^\ast$-algebra is defined by
	$$A_q(\mathcal H_{\mathbb R},U)=C^*(1,s_q(\xi):\xi\in \mathcal H_{\mathbb R}).$$
	The algebra $M_q(\mathcal H_{\mathbb R},U):=A_q(\mathcal H_{\mathbb R},U)''$ is called the {\em $q$-deformed Araki-Woods von Neumann algebra}. It is  known that $M_q(\mathcal H_{\mathbb R},U)$ is a factor \cite{KSW23}, so we will call it the {\em  $q$-Araki-Woods factor}.
	
	The vacuum state $\varphi_q = \ip{\Omega}{\cdot\Omega}_q$ is a normal faithful state on $M_q(\mathcal H_{\mathbb R},U)$, called the {\em $q$-quasi free state}, and $\mathcal F_q(\mathcal H_U)$ is the GNS  Hilbert space of $M_q(\mathcal H_{\mathbb R},U)$. Since $\Omega$ is a cyclic and separating vector for $M_q(\mathcal H_{\mathbb R},U)$, we identify
	$L^2\big(M_q(\mathcal H_{\mathbb R},U),\ph_q\big)$ with $\mathcal F_q(\mathcal H_U)$ by $x\Omega_{\ph_q}\mapsto x\Omega$.
	For algebraic tensors over $\mathcal H_{\mathbb C}$, the linear Wick map (see \cite[Section 2]{KW24})
	\[
	W: \bigoplus_{n\geq0}\mathcal H_{\mathbb C}^{\ten n}
	\to A_q(\mathcal H_{\mathbb R}, U)
	\]
	is the unique map such that  $W(\xi)\Omega=\xi$ for $\xi\in \mathcal H_{\mathbb C}^{\ten n}$. In particular
	$W(\xi)=\ell_q(\xi)+\ell_q(I\xi)^*$ for $\xi\in \mathcal H_{\mathbb C}$.
	We adopt the following notations for our convenience.
	\begin{enumerate}
		\item $\xi_1\ldots\xi_n=\xi_1\otimes\ldots\otimes\xi_n$, for $\xi_i\in\mathcal H_U$, $1\leq i\leq n$.
		\item $\xi^n=\xi^{\otimes n}$, $n\geq 1$ and $\xi^0=\Omega$ for all $\xi\in\mathcal H_U$. 
		\item $d_0=1$, $d_n=\prod_{j=1}^n(1-q^j)$ for $n\geq 1$ and $d_\infty=\prod_{j=1}^\infty(1-q^j)$.
		\item $C_q=\prod_{j=1}^\infty\frac{1}{1-|q|^j}$.
		\item $[n]_q=\sum_{j=0}^{n-1}q^j=\frac{1-q^n}{1-q}$.
		\item $[n]_q!=\prod_{j=1}^n[j]_q=\frac{d_n}{(1-q)^n}$.
		\item $\begin{pmatrix}
			n\\k
		\end{pmatrix}_q=\frac{[n]_q!}{[k]_q![n-k]_q!} = \frac{d_n}{d_kd_{n-k}}$ for any $1\leq k\leq n$.
	\end{enumerate}
	
	We collect below several well-known results from \cite{BKS97,BMRW23,KW24} that will be used frequently in the final two sections.
	\begin{lemma}\label{lemma;some facts}
		For any $\xi\in\mathcal H_U$ and $n\ge1$, we have
		\begin{enumerate}
			\item[(i)] for any $\eta\in \mathcal H_U^{\otimes_q^n}$ and $\zeta\in \mathcal H_U^{\otimes_q^m}$ with $m\ge 1$,
			$$\ell_q^*(\xi)(\eta\otimes\zeta)=(\ell_q^*(\xi)\eta)\otimes\zeta+q^n\eta\otimes \ell_q^*(\xi)\zeta.$$
			
			\item[(ii)]  We have $\|\ell_q(\xi)^n\|\le C_q(1-q)^{-\frac{n}{2}}\|\xi\|_U^n$.
			
			\item[(iii)] For $m\ge n$, 
			\begin{align}\label{eq:normal-ordering}
				\ell_q^*(\xi)^n\ell_q(\xi)^m=\frac{[m]_q!}{[m-n]_q!}\|\xi\|_U^{2n}\ell_q(\xi)^{m-n}+\sum_{k=0}^{n-1}\frac{[m]_q!}{[m-k]_q!}q^{m-k}\|\xi\|_U^{2k}\ell_q^*(\xi)^{n-k-1}\ell_q(\xi)^{m-k}\ell_q^*(\xi).
			\end{align}
			
			\item[(iv)] For $0\le i\le m$,
			\begin{equation}\label{eq:q-fractional relations}
				[m-i]_q\binom{m}{i}_q
				=[m]_q\binom{m-1}{i}_q.
			\end{equation}
			
			\item[(v)]  For $0\le i\le m$,
			\begin{equation}\label{eq;Pascal identity}
				\binom{m}{i}_q + q^{i+1}\binom{m}{i+1}_q
				=\binom{m+1}{i+1}_q.
			\end{equation}
			\item [(vi)]  $d_n\leq C_q$, $\frac{1}{d_n}\leq C_q$, $\begin{pmatrix}
				n\\k
			\end{pmatrix}_q\leq C_q$
		\end{enumerate}
	\end{lemma}

	\section{Derivations from finite free fisher information}\label{sec:derivations on free fisher infomration}
	In this section, we follow a general framework in which a $C^\ast$-algebra is generated by finitely many eigenoperators in the sense of \cite{Nel17}. The $q$-Araki-Woods $C^\ast$-algebras $A_q(\mathcal H_{\mathbb R}, U)$ are particular examples when $\dim_{\mathbb R} \mathcal H_{\mathbb R}<\infty$. Motivated by \cite{AJW26}, we first establish general estimates for conjugate variables and then use them to prove the simplicity of $C^\ast$-algebras satisfying a certain rapid decay property. We finally apply them to prove the selflessness of $q$-Araki-Woods $C^\ast$-algebras.
	
	We begin by fixing some notation.
	Let $(M,\varphi)$ be a von Neumann algebra   equipped with a normal faithful state $\varphi$. Let $(L^2(M,\varphi), \pi_\varphi,\Omega_\varphi)$ denote the GNS space, that is, $L^2(M,\varphi)$ is a Hilbert space, $\pi_\varphi:M\to B(L^2(M,\varphi))$ is a normal $\ast$-homomorphism and $\Omega_\varphi$  is the cyclic and separating vector in $L^2(M,\varphi)$ such that $\varphi(x)=\langle\Omega_\varphi,\pi_\varphi(x)\Omega_\varphi\rangle$ for $x\in M$. We  denote the $L^2$-norm and inner product on $L^2(M,\varphi)$ by $\|\cdot\|_{2,\varphi}$ and $\langle\cdot,\cdot\rangle_\varphi$, respectively, and simply write $\|\cdot\|_2$ when $\varphi$ is clear from the context.  
	
	The Tomita operator $S_\varphi$ is the unbounded antilinear operator on $L^2(M,\varphi)$ given by $x\Omega_\varphi\mapsto x^\ast\Omega_\varphi, x\in M$. Let $S_\varphi=J_\varphi\Delta_\varphi^{1/2}$ be the polar decomposition. Here $J_\varphi$ is an antilinear unitary and is called the modular conjugation operator. Moreover, $\Delta_\varphi^{it}\pi_\varphi(M)\Delta_\varphi^{-it}=\pi_\varphi(M)$  for $t\in\mathbb R$, and the {\em modular automorphism group} $(\sigma_t^\varphi)_{t\in \mathbb R}$ on $M$ associated to $\varphi$ is defined by $\sigma_t^\varphi(x):=\pi_\varphi^{-1}(\Delta_\varphi^{it}\pi_\varphi(x)\Delta_\varphi^{-it})$ for $x\in M$. We often suppress the notation $\pi_\varphi$. 
	Denote by $M^\varphi$ the centralizer
	\begin{align*}
		M^\varphi:=\{x\in M: \varphi(xy)=\varphi(yx),\;\;\forall y\in M\}.
	\end{align*}
	Let $M^{\mathrm {op}}$ denote the opposite von Neumann algebra with the associated state $\ph^{\mathrm{op}}$. 
	We naturally identify $(M^\ph)^{\mathrm{op}}$ with $(M^{\mathrm{op}})^{\ph^{\mathrm{op}}}$.
	We now follow the conventions of \cite{Nel17}. For any $x^{\mathrm{op}}\in M^{\mathrm{op}}$, we suppress the superscript $\mathrm{op}$ and simply write $x$. We consider the tensor
	multiplication $\#$, the involution $\dagger$, and the outer bimodule actions given by
	\[
	(x\ten y)\#(z\ten w)=xz\ten wy,\qquad (x\ten y)^\dagger=y^*\ten x^*,\qquad
	a\cdot(x\ten y)\cdot b=ax\ten yb,
	\]
	for $x,y,z,w,a,b\in M$. Write $M\otimes M^{\mathrm{op}}$ for the algebraic tensor product of $M$ and $M^{\mathrm{op}}$.
	
	\subsection{Modular derivations and conjugate variables} In von Neumann algebras, modular derivations have proved to be a powerful tool in Popa's deformation/rigidity theory \cite{Pet09, DI16, KW26} and Free probability \cite{Dab10, Nel17}. Here, we study  $\mu$-modular derivations in Nelson’s sense\cite{Nel17}  and their conjugate variables.    
	\begin{definition}
		Let $B\subseteq M$ be a unital weakly dense $\ast$-subalgebra. A {\em derivation} on $ B$ is a $\mathbb C$-linear map $\delta :  B \to M \otimes M^{\mathrm{op}}$ satisfying
		the Leibniz rule: $\delta(ab)=a\delta(b)+\delta(a)b$ for $a,b\in B$.
		We call $ B$ the domain of $\delta$ and write $\dom (\delta) :=  B$. The {\em conjugate derivation} to $\delta$, denoted by $\hat\delta$, is a
		derivation with $\dom(\hat\delta) = \dom (\delta)$ defined by $\hat \delta(x) = \delta(x^\ast)^\dagger$ for $x\in\dom (\delta)$.
	\end{definition}

	\begin{definition}\label{def:conjugate_variable}
		The {\em conjugate variable} to a derivation $\delta$ is defined to be the vector $\xi=\delta^\ast(1\otimes 1)\in L^2(M,\varphi)$ whenever it exists. In this case, $\delta$ is closable \cite[Lemma 3.9]{Nel17}. Moreover, the conjugate variable is called {\em Lipschitz} if $\xi\in \dom(\bar\delta)$ and $\bar\delta(\xi)\in M\overline\otimes M^{\mathrm{op}}$. 
	\end{definition}

	Let $\delta:\dom (\delta)\to M\otimes M^{\mathrm{op}}$ be a derivation. Throughout,
	we assume that $\dom (\delta)\subseteq M_{\infty}$  and $\delta\big(\dom (\delta)\big) \subseteq M_\infty\otimes M_\infty^{\mathrm{op}}$, where $M_\infty$ is the subalgebras of $M$ consisting of analytic elements of $\sigma^\varphi$.
	Following \cite[Definition 3.2]{Nel17},  we say that $\delta$ is {\em $\lambda$-modular} for $\lambda>0$ if $\sigma_t^\varphi(\dom\delta)\subseteq\dom\delta$ for all $t\in\mathbb R$ and it satisfies
	\[\delta\circ\sigma_t^\ph(x)=\lambda^{it}(\sigma_t^\ph\otimes \sigma_t^\ph)\circ\delta(x),\qquad x\in \dom (\delta), t\in\mathbb R.\]
	We identify $M$ with the dense subspace of $L^2(M,\varphi)$ via the map $x\mapsto x\Omega_\varphi$.
	We treat $\delta$ as a densely defined operator from $L^2(M,\varphi)$ to $L^2(M,\varphi)\otimes L^2(M^{\mathrm{op}}, \varphi^{\mathrm{op}})$. So, 
	if $\delta$ is closable as an operator on $L^2(M,\varphi)$ with closure $\bar\delta$, then $\lambda$-modularity is equivalent to
	\[\bar\delta(\eta)=\lambda^{z}(\Delta_\ph^z\otimes\Delta^z_\ph)\circ\bar\delta(\Delta_\ph^{-z}\eta),\qquad \eta\in\dom \bar\delta,\;\; \Delta_\ph^{-z}\eta\in\dom\bar\delta\] 
	for all $z\in\mathbb C$.
	It is straightforward that if $\delta$ is $\lambda$-modular for some $\lambda>0$, then $\hat\delta$ is $\lambda^{-1}$-modular such that $\hat{\hat\delta}=\delta$. One easily verifies that if $a\in G$ is an eigenoperator of $\sigma^\varphi$ with eigenvalue $\lambda$, then $\delta_a$ is $\lambda$-modular.

	The following proposition is a nontracial analogue of the one-sided adjoint estimate in \cite[Proposition~3.5 and Remark~3.6]{AJW26}.
	
	\begin{proposition}\label{prop:generalduals}
		Let $\delta:\dom(\delta)\to M_\infty\otimes M_\infty^{\mathrm{op}}$ be a $\lambda$-modular derivation for some $\lambda>0$. Suppose that $\dom(\delta)$ is a weakly dense unital $\ast$-subalgebra spanned by  eigenoperators of $\sigma^\varphi$,
		that the conjugate variable $\xi$ of $\delta$ exists and is Lipschitz, and that $(\id\ten \hat\delta )\delta =(\delta\ten\id) \hat\delta$. 
		Then $y:=(\id\ten\ph)\overline{\delta }\xi\in M^\varphi$, and the map $D_\delta:=(\id\ten\ph)\delta $ defined on $\dom(\delta)$
		extends to a bounded operator on $L^2(M,\ph)$ such that
		$\norm{D_\delta}\leq\norm{y}^{1/2}$. 
	\end{proposition}
	\begin{proof}
		
		We first show that $y\in M^\varphi$.
		Applying \cite[Lemma~3.10]{Nel17} yields $\xi \in \dom(\Delta_\ph^z)$ and $\Delta_\ph^z\xi =\lambda^z\xi $ for all $z\in \mathbb C$. Since $\delta $ is $\lambda$-modular, we have $\lambda^{z}(\Delta_\ph^z\otimes\Delta_\ph^z)\bar\delta \xi =\bar\delta (\Delta_\varphi^z\xi )$, which yields 
		\begin{align*}
			\Delta_\ph^z y=(\mathrm{id}\otimes\ph)(\Delta_\ph^{z}\otimes\Delta_\ph^z)\bar\delta \xi =\lambda^{-z}(\mathrm{id}\otimes\ph)\bar\delta (\Delta_\ph^z\xi )=(\mathrm{id}\otimes\ph)\bar\delta \xi =y,
		\end{align*}
		for all $z\in\mathbb C$. It follows that $y\in M^\varphi$.
		We next prove that $D_\delta$ is bounded. Since for $p\in{\mathbb C(G)}$, $D_\delta^\ast(p)=\delta^\ast(p\otimes1)$, it suffices to estimate the norm $\|\delta^\ast(p\otimes 1)\|$. Fix $p\in\dom\delta$.
		By \cite[Lemma~3.9]{Nel17},  we have $p\otimes1\in \dom(\delta ^\ast)$ and
		\begin{equation*}
			\delta ^\ast(p\otimes1)=p\xi -(\mathrm{id}\otimes\ph)\hat\delta(p).
		\end{equation*}
		Define the map $T:\dom(\delta)\to \dom(\delta)$ by  $T(p)=(\id\ten\ph)\bigl( \hat\delta (p^*)\cdot p\bigr)$ for $p\in \dom(\delta)$.
		Applying the Leibniz rule to the preceding identity, we obtain
		\begin{align*}
			\delta ^\ast(p^*p\otimes1)&=p^\ast p\xi -(\mathrm{id}\otimes\ph) \hat\delta (p^\ast p)\\
			&=p^\ast p\xi -(\mathrm{id}\otimes\ph)(p^\ast\hat\delta (p))-T(p)\\
			&=p^\ast (p\xi -(\mathrm{id}\otimes\ph)(\hat\delta (p)))-T(p)\\
			&=p^*\delta ^\ast(p\otimes1)-T(p).
		\end{align*}
		Since $\hat{\hat{\delta}}=\delta$, it follows that $T(p)^*=(\ph\ten\id)(p^*\cdot\delta p)$. Combining this with the identity $(\id\ten \hat\delta )\delta p
		=(\delta\ten\id) \hat\delta  p$ 
		yields
		\[
		\begin{aligned}
			\ip{T(p)}{\xi }&=(\ph\otimes\ph)(\delta (Tp)^\ast)\\
			&=\ph\otimes\ph\big( \hat\delta (T(p)^\ast )\big)\\
			&=\varphi\otimes\varphi\Big(\hat\delta\big((\varphi\otimes\mathrm{id})(p^\ast\delta p)\big)\Big)\\
			&=\ph\otimes\ph\otimes\ph \big((p^*\ten1\ten1)(\id\ten \hat\delta )\delta p\big)\\
			&=\ph\otimes\ph\otimes\ph \big((p^*\ten1\ten1)(\delta\ten\id) \hat\delta  p\big)\\
			&=\varphi\otimes\varphi\big((p^\ast\otimes1)\delta((\mathrm{id}\otimes\varphi)\hat{\delta}p)\big)\\
			&=\ip{\delta ^\ast(p\ten1)}{(\mathrm{id}\otimes\ph) \hat\delta (p)}.
		\end{aligned}
		\]
		Consequently, 
		\begin{equation*}
			\begin{aligned}
				\norm{\delta ^\ast(p\otimes1)}_2^2
				&=\ip{\delta ^\ast(p\otimes 1)}{p\xi}-\ip{\delta ^\ast(p\otimes1)}{(\mathrm{id}\otimes\ph) \hat\delta (p)}\\
				&=\ip{p^*\delta ^\ast(p\otimes1)-T(p)}{\xi }\\
				&  =\ip{\delta ^\ast(p^*p\otimes1)}{\xi }\\
				&=\ip{p^*p\ten1}{\overline\delta \xi }\\
				&  =\ph(p^*py).
			\end{aligned}
		\end{equation*}
		By Lemma \ref{lem:dense elements in centralizer} and Kaplansky's density theorem,  $\ph(xy)\geq0$ for all $x\in M^\ph$ positive. Since $y\in M^{\varphi}$, it follows that $y$ is positive and then
		\begin{align*}
			\norm{\delta ^\ast(p\otimes1)}_2^2  =\ph(y^{1/2}p^\ast py^{1/2})=\|py^{1/2}\Omega_\varphi\|^2=\|p\Omega_\varphi y^{1/2}\|^2
			\leq \|y\|\,\|p\|^2_2.
		\end{align*}
		Therefore, $D_\delta^\ast$ is bounded with $\|D_\delta^\ast\|\leq \|y\|^{1/2}$.
	\end{proof}
	
	\subsection{Eigenoperators and free difference quotient} We now  consider Voiculescu’s and Nelson's
	free difference quotients derivations.
	\begin{definition}\label{def:eigenoperator}
		An element $y\in M$ is called an {\em eigenoperator of $\sigma^\varphi$ with eigenvalue $\lambda>0$} if $\sigma_t^\varphi(y)=\lambda^{it}y$ for all $t\in\mathbb R$.
	\end{definition}
	
	Throughout Sections \ref{sec:derivations on free fisher infomration}--\ref{sec:simplicity of C-algebras with eigenoperators}, 
	we assume that $(M,\varphi)$ is generated by a finite set $G = G^*$ consisting of eigenoperators of $\sigma^\varphi$, and assume that $G$ is algebraically free.
	Denote by $\mathbb C\langle G\rangle$ the unital $\ast$-algebra generated by $G$.
	For any $a\in G$, the {\em free difference quotient with respect to $a$} is the derivation 
	$\delta_a:\mathbb C\langle G\rangle\to\mathbb C\langle G\rangle\otimes \mathbb C\langle G\rangle^{\mathrm{op}}$ given by
	\begin{align*}
		\delta_a(a)=1\otimes 1,\;\;\delta_a(b)=0\;\;\;\mbox{ for } b \in G\setminus \{a\}.
	\end{align*}
	Note that $\delta_{a^\ast}=\hat\delta_a$  and  $(\id\ten \hat\delta_a )\delta_a 
	=(\delta_a\ten\id) \hat\delta_a  $.
	We say that the generating set $G$ has {\em finite free Fisher information} if 
	$\delta_a^\ast(1\otimes1)$ exists 
	for every $a\in G$. 
	
	A {\em word} in $G$ is an element of the form $g_1\ldots g_n$ for some $g_j\in G$. Every word in $G$ is an eigenoperator of $\sigma^{\varphi}$. Indeed, if $\sigma_t^\varphi(g_j)=\lambda_j^{it}g_j$ for some $\lambda_j>0$ and every $1\leq j\leq n$ and $t\in\mathbb R$, then $\sigma_t^\varphi(g_1\ldots g_n)=\lambda^{it}g_1\ldots g_n$, where $\lambda=\lambda_1\ldots\lambda_n$. 
	
	\begin{definition}\label{def:balanced word}
		A word $g=g_1\ldots g_n$ in $G$ is called {\em balanced} if $g\in M^\varphi$, that is, $g$ is an eigenoperator of $\sigma^\varphi$ with eigenvalue 1. A polynomial $p\in \mathbb C\langle G\rangle$ is called {\em balanced} if each monomial appearing in $p$ is balanced.      
	\end{definition}

	\medskip
	We now characterize the centralizer  $M^\varphi$. The following lemma is the special case $\lambda=1$ of \cite[Lemma~2.8]{Nel17}.
	
	\begin{lemma}\label{lem:dense elements in centralizer}
		The centralizer $M^\varphi$ is generated by the balanced words in $G$. Equivalently, $M^\varphi$ is the $\sigma$-weakly closed linear span of
		$\bigl\{p^*p:p\in\mathbb C\langle G\rangle\cap M^\varphi\bigr\}$.
	\end{lemma}


	For $x\in M_\infty$, define $R(x)=J_\ph\sigma_{-i/2}^\ph(x^\ast)J_\ph$, that is, $R(x)y\Omega_\ph=yx\Omega_\ph$ for any $y\in M$.
	
	\begin{corollary}\label{cor:boundedconjugates}
		Suppose that $G$ has
		finite free Fisher information and  Lipschitz conjugate variables. 
		Then for each $c\in G$, the operator $D_c:=(\mathrm{id}\otimes\varphi)\bar\delta_c$ is bounded, and there exists an eigenoperator $x_c\in M$ of $\sigma^\varphi$ such that $\xi_c=x_c\Omega_\varphi$ and
		$D_c^*=R(x_c)-D_{c^\ast}$.
	\end{corollary}
	\begin{proof}
		Fix an eigenoperator $c\in G$, and suppose that $c$ has eigenvalue $\lambda>0$. It follows from Proposition \ref{prop:generalduals} that $D_c$ and $D_{c^\ast}$ are bounded operators on $L^2(M,\varphi)$.
		Since $\delta_{c^\ast}=\hat\delta_c$,
		as noted in the proof of Proposition \ref{prop:generalduals}, we have
		$D_{c^\ast}(p)=p\xi_c-D_c^\ast(p)$ for any $p\in\dom(\delta)$. 
		It follows that $p\mapsto p\xi_c$ extends to a bounded operator on $L^2(M,\varphi)$. Thus $\xi_c$ is a $\varphi$-right bounded vector, that is,  there exists $y_c\in M$ such that $\xi_c=\Omega_\varphi y_c=J_\varphi y_c^\ast \Omega_\varphi$. Hence $D_{c^\ast}=J_\varphi y_c^\ast J_\varphi-D_c^\ast$.
		Since $\delta_c$ is $\lambda$-modular and $\delta_{c^\ast}$ is $\lambda^{-1}$-modular, we have $\Delta_\ph^{it}D_c\Delta_\ph^{-it}=\lambda^{-it}D_c$ and
		$\Delta_\ph^{it}D_{c^\ast}\Delta_\ph^{-it}=\lambda^{it}D_{c^\ast}$ for all $t\in\mathbb R$. This implies that $\Delta_\varphi^{it}J_\varphi y_c^\ast J_\varphi\Delta_\varphi^{-it}=\lambda^{it}J_\varphi y_c^\ast J_\varphi$ and then $\sigma_t^\varphi(y_c^\ast)=\lambda^{-it}y_c^\ast$ for $t\in \mathbb R$. By setting $x_c=\sigma_{-i/2}^\varphi(y_c)$, we obtain the required result. In fact, $x_c$ is an eigenoperator of $\sigma^\varphi$ with the eigenvalue $\lambda$.
	\end{proof}

	\subsection{Spectral gap on \texorpdfstring{$C^*$}{C*}-algebras generated by eigenoperators}\label{sec:cancellation}
	We continue under the standing assumption that $(M,\varphi)$ is generated by a finite self-adjoint set $G$ of eigenoperators of $\sigma^\varphi$, and assume in addition that $G$ has finite free Fisher information and Lipschitz conjugate variables.

	Let $c\in G$ be an eigenoperator of $\sigma^\ph$ with eigenvalue $0<\lambda< 1$.
	Write the polar decompositions of $c$ and $c^*$ as $c=vH_1$ and $c^\ast=v^\ast H_2$, where
	$H_1=(c^*c)^{1/2}$ and $H_2=(cc^*)^{1/2}$. Since $G$ has finite free Fisher information,
	it follows from
	\cite[Theorem 5.9, Corollary~5.10 and Example~5.12]{Nel17} that
	$c^*c$ has no atoms with trivial kernel and $cc^\ast$ has an atom at $0$ of size $1-\lambda$. Hence $v^*v=1$ and $\tau(vv^\ast)=\lambda$. Moreover, 
	$\sigma_t^\ph(v)=\lambda^{it}v$ for $t\in\mathbb R$.
	Let $ E_\varphi:M\to {M^\ph}$ be the $\ph$-preserving normal faithful conditional
	expectation. Following \cite[Definition 6.1]{Nel17}, define derivations $\widetilde\delta_{j}:\mathbb C\langle G\rangle \to L^2(M\otimes M^{\mathrm{op}})$ by
	\[
	\begin{aligned}
		\widetilde \delta_{1}(x)
		&=\delta_cx\#(v\ten1)+\delta_{c^*}x\#(1\ten v^*),\\
		\widetilde \delta_{2}(x)
		&=\delta_{c^*}x\#(v^*\ten1)+\delta_cx\#(1\ten v).
	\end{aligned}
	\]
	Note that $\widetilde{\delta}_1=\delta_{|c|}$ and $\widetilde\delta_{2}=\delta_{|c^\ast|}$ in the notation of \cite{Nel17}.  By \cite[Lemma 6.1]{Nel17}, each $\widetilde\delta_j$ is closable such that $H_k\in \overline{\widetilde{\delta_j}}$ for $1\leq k,j\leq 2$.

	We warn the reader that the notation and results used below appear in version~3 of the arXiv preprint of \cite{Nel17}, but not in the published version.
	Set  $\mathcal D = \mathbb C\langle G,H_1,H_2\rangle\cap M^\varphi$, and define derivations $\delta_j:\mathcal D\to L^2(M^\varphi\otimes (M^\varphi)^{\mathrm{op}})$ by 
	\[ \delta_j =  (E_\varphi\otimes E_\varphi^{\mathrm{op}}) \circ\overline{\widetilde{\delta}_j}_{|_{\mathcal D}}, \qquad j=1,2. \] 
	Define the derivation $\delta_0:\mathcal D\to M^\ph\otimes(M^\ph)^{\mathrm{op}}$ by
	\[
	\delta_0(x)=\sum_{j=1}^2
	\delta_j(x)\#(H_j\ten1-1\ten H_j)
	-(x\ten1-1\ten x),\qquad x\in \mathcal{D}.\]
	Set $\tau=\ph_{|_{M^\ph}}$, $U_1=1\ten1$ and  $U_2=1\ten1-(1-vv^\ast)\ten(1-vv^\ast)$. The following result combines Lemmas~6.3 and~6.4 from version~3 of the arXiv preprint of \cite{Nel17}. We note that Equation~(16) in Lemma~6.4 of that preprint contains a sign error. The correct formula is the one stated below and follows from the proof given there.
	
	\begin{lemma}\label{lem:nelson-cancellation}
		Let $\delta_j$ and $\mathcal D$ be as above for $j=0,1,2$. Then the following hold.
		\begin{enumerate}
			\item[(i)] For $j = 1,2$, $\delta_j$ is a closable real derivation such that
			$1\ten1\in\dom\delta_j^*$,
			$\delta_j(H_k)=\delta_{jk}U_j$ and $
			\delta_j(x)\#U_j=\delta_j(x).$
			Here reality means $\delta_j(x^*)^\dagger=\delta_j(x)$.
			In particular, $\delta_j^*(U_j)=\delta_j^*(1\ten1)$.
			
			\item[(ii)] We have $1\otimes1\in\dom(\delta_0^\ast)$  and hence $\delta_0$ is closable. Its closure $\bar\delta_0$ satisfies $L^2(W^*(H_1,H_2),\tau)\subseteq\ker\overline{\delta_0}$. Moreover, for
			$x\in \mathcal{D}$,
			\begin{equation*}
				\|[x,U_2]\|_2^2
				={}\langle\bar\delta_0(x),[x,1\otimes1]\rangle+\sum_{j=1}^2\Bigl(
				\langle[x,H_j],[x,\delta_j^\ast(1\otimes1)]\rangle_\tau
				+2\operatorname{Re}\langle (\tau\otimes1-1\otimes\tau)\circ\bar\delta_j(x),[x,H_j]\rangle_\tau
				\Bigr).
			\end{equation*}
		\end{enumerate}
	\end{lemma}

	Before recording the next lemma, we recall some notation from \cite{Nel17}. For any self-adjoint element $a\in M$ with spectrum contained in a compact  interval $I\subseteq\mathbb R$, consider $LR_a:C(I\times I)\to M\otimes M^{\mathrm{op}}$ defined by 
	\[LR_a(f\otimes g)=f(a)\otimes g(a), \qquad f,g\in C(I)\]
	where $f\otimes g\in C(I)\otimes C(I)=C(I\times I)$. For $f\in C^1(I)$, divided difference of $f$:
	\[
	\widetilde f(s,t)=
	\begin{cases}
		\dfrac{f(t)-f(s)}{t-s},&s\ne t,\\[6pt]
		f'(t),&s=t.
	\end{cases}
	\]

	\begin{lemma}\label{lem:polar-part-domain}
		Let $c$, $v$ be as above and let $D_c$ be the bounded operator from Corollary \ref{cor:boundedconjugates}.
		Then $v\in \operatorname{dom}\bar\delta_c$.
		Moreover,  $D_c(v\Omega_\varphi), D_{c^\ast}(v^\ast\Omega_\varphi)\in L^2(M^\varphi,\tau)$ such that $J_\tau D_c(v\Omega_\varphi)=D_c(v\Omega_\varphi)$ and $J_\tau D_{c^\ast}(v^\ast\Omega_\varphi)=D_{c^\ast}(v^\ast\Omega_\varphi)$.
	\end{lemma}
	\begin{proof}
		To show $v\in\dom(\bar\delta_c)$, since $v_\varepsilon:=c(c^\ast c+\epsilon)^{-1/2}\to v$ strongly as $\varepsilon\to 0$, it suffices to show that $v_\varepsilon\in\dom(\bar{\delta}_c)$ for every $\varepsilon>0$ and that $\{\bar{\delta}_c(v_\varepsilon)\}_{\varepsilon>0}$ is bounded in $L^2$-norm. Indeed, any weakly convergent subnet of $\{\bar{\delta}_c(v_\varepsilon)\}_{\varepsilon>0}$, together with the weak closedness of the graph of $\bar{\delta}_c$, yields $v\in\dom(\bar{\delta}_c)$.
		
		To prove the two claims, write $a=c^\ast c$ and $I=[0, \|a\|]\subseteq\mathbb R$.  For $f\in C^1(I)$, \cite[Lemma 4.3]{Nel17} gives
		$\bar\delta_c(f(a))
		=LR_a(\widetilde f)\#(c^*\otimes1).$
		Since $f(a)\in M^\varphi$, \cite[Lemma 4.2]{Nel17} gives $cf(a)\in\operatorname{dom}\bar\delta_c$ and
		\begin{align}\label{eq:formula for cf(a)}
			\bar\delta_c(cf(a))
			=1\otimes f(a)+c\cdot\bigl(LR_a(\widetilde f)\#(c^*\otimes1)\bigr).
		\end{align}
		Set $f_\varepsilon(t)=(t+\varepsilon)^{-1/2}, t\in I$ and $\varepsilon>0$. Then $v_\varepsilon=cf_\varepsilon(a)\in \dom(\bar\delta_c)$. This proves the first claim.
		Now if we define 
		$K_\varepsilon:I\times I\to [0,\infty)$  by 
		$K_\varepsilon(s,t)=f_\varepsilon(t)+s\widetilde f_\varepsilon(s,t),$ for $ s,t\in I$, then we get $\bar\delta_c(cf_\varepsilon(a))=K_\varepsilon(cc^\ast\otimes1, 1\otimes c^\ast c)$. To see this, first take $f,g\in C^1(I)$, then
		\[c\cdot(LR_a(f\otimes g)\# (c^\ast\otimes1))=cf(c^\ast c)c^\ast\otimes g(c^\ast c)=cc^\ast f(cc^\ast)\otimes g(c^\ast c)\]
		and then by approximation, we have $c\cdot(LR_a(\widetilde f_\varepsilon)\# (c^\ast\otimes1))=(cc^\ast\otimes 1) \widetilde f_\varepsilon(cc^\ast\otimes, 1\otimes c^\ast c) $, which finally gives the formula for $\bar\delta_c(cf_\varepsilon(c^\ast c))$.
		As $K_\varepsilon$ is the divided difference of the increasing function $tf_\varepsilon(t)$, it is nonnegative. A simple calculation shows  $\widetilde{f}_\varepsilon\leq 0$, so we have $0\leq K_\varepsilon(s,t)\leq f_\varepsilon(t)$ on $I\times I$. The joint spectral calculus gives 
		\begin{align*}
			\|\bar\delta_c(v_\varepsilon)\|_2=\|K_\varepsilon(cc^\ast \otimes1, 1\otimes c^\ast c)\|_2\leq \|1\otimes f_\varepsilon(a)\|_2=\varphi((a+\varepsilon)^{-1})^{1/2}.
		\end{align*}
		We now use the boundedness of $D_c$ to obtain a uniform bound on the right-hand side.
		Set $z_\varepsilon=(a+\varepsilon)^{-1}$. By considering the function $(t+\varepsilon)^{-1}$ in \eqref{eq:formula for cf(a)}, we get $\bar\delta_c(cz_\varepsilon)=1\otimes z_\varepsilon-cz_\varepsilon c^\ast\otimes z_\varepsilon$. Hence $D_c(cz_\varepsilon)=\varphi(z_\varepsilon)(1-cz_\varepsilon c^\ast)$. In particular, multiplying on the left by $1-vv^*$ gives $(1-vv^*)D_c(cz_\varepsilon)=\varphi(z_\varepsilon)(1-vv^*)$, since $(1-vv^*)c=0$.
		Since $\varphi(vv^\ast)=\lambda$, we have
		\begin{align*}
			\varphi(z_\varepsilon)^2(1-\lambda)=\|\varphi(z_\varepsilon)(1-vv^\ast)\|^2_2\leq \|D_c\|^2\|cz_\varepsilon\|^2_2\leq \|D_c\|^2\,\varphi(z_\varepsilon c^\ast cz_\varepsilon)\leq \|D_c\|^2\,\varphi(z_\varepsilon),
		\end{align*}
		where the last inequality follows from $z_\varepsilon c^\ast cz_\varepsilon\leq z_\varepsilon$. This  implies that $\varphi(z_\varepsilon)\leq \frac{1}{1-\lambda}\|D_c\|^2$, and hence $\sup_{\varepsilon>0}\|\bar\delta_c(v_\varepsilon)\|_2<\infty$ proving the second claim. 
		
		Finally, we have $\bar\delta_c(v_\varepsilon)\in M^\varphi\bar\otimes(M^\varphi)^{\mathrm{op}}$ and $\bar\delta_c(v_\varepsilon)\ge0$. Its slice $D_c(v_\varepsilon\Omega_\varphi)$ is therefore represented by a positive element of $M^\varphi$, so it lies in $L^2(M^\varphi,\tau)$ and is fixed by $J_\tau$. Boundedness of $D_c$ gives
		$D_c(v_\varepsilon\Omega_\varphi)\to D_c(v\Omega_\varphi)$ in $L^2(M,\varphi)$ as $\varepsilon\to 0$ and hence $J_\tau D_c(v\Omega_\varphi)=D_c(v\Omega_\varphi)$.  
		Moreover, we use \cite[Lemma 4.1]{Nel17} to conclude that $v^\ast\in \dom(\bar\delta_{c^\ast})$ and $\bar\delta_{c^\ast}(v^\ast)=\overline{\hat\delta_c}(v^\ast)=\bar\delta_c(v)^{\dagger}$. Applying $\mathrm{id}\otimes\varphi$ will give the required result for $D_{c^\ast}(v^\ast\Omega_\varphi)$.
	\end{proof}

	\begin{lemma}\label{lem:spectral gap fo kernel of delta0}
		\label{lem:an estimate on L2-norm in terms of H1,H2}
		With  $\delta_0$  as above,
		there exists  $C>0$ such that for  $f\in  \ker(\bar\delta_0)\subseteq L^2(M^\varphi,\tau)$
		\begin{equation*}
			\norm{f-\tau(f)1}_2
			\leq C\bigl(\norm{[f,{H_1}]}_2+\norm{[f,{H_2}]}_2\bigr).
		\end{equation*}
	\end{lemma}
	\begin{proof}
		Set $q=vv^\ast$. Note that
		$U_2=1\otimes q+q\otimes(1-q)=q\otimes 1+(1-q)\otimes q.$  For $f\in \mathcal{D}$ with $\tau(f)=0$, since $\tau(q)=\lambda$,  we have
		\[
		\|f U_2\|_2^2
		=\|fq\|_2^2+\lambda\|f(1-q)\|_2^2,\quad
		\|U_2 f\|_2^2
		=\|qf\|_2^2+\lambda\|(1-q)f\|_2^2,\quad \langle fU_2,U_2f\rangle= -|\tau(qf)|^2.
		\]
		Therefore,
		\[
		\begin{aligned}
			\|[f,U_2]\|_2^2
			&={}\|fq\|_2^2+\lambda\|f(1-q)\|_2^2+\|qf\|_2^2+\lambda\|(1-q)f\|_2^2
			+2|\tau(qf)|^2\\
			&\geq \lambda\|fq\|_2^2+\lambda\|f(1-q)\|_2^2+\lambda\|qf\|_2^2+\lambda \|(1-q)f\|_2^2\\
			&=2\lambda\|f\|_2^2.
		\end{aligned}
		\]
		Thus it suffices to estimate $\|[f,U_2]\|_2$. Recall the bounded operators $D_c$ and $D_{c^\ast}$ from Corollary \ref{cor:boundedconjugates}. By Lemma \ref{lem:polar-part-domain},
		$\zeta_1:=D_c(v\Omega_\varphi)$ and $\zeta_2:=D_{c^\ast}(v^\ast\Omega_\varphi)$ belong to $L^2(M^\ph,\tau)$ and satisfy $J_\tau\zeta_j=\zeta_j$ for $j = 1,2$. Here $J_\tau$ is the antiunitary involution on $L^2(M^\ph,\tau)$ given by $x\Omega_\varphi\mapsto x^\ast\Omega_\varphi$. Recall that $R(v)=J_\varphi \sigma_{-i/2}^\varphi(v^\ast)J_\varphi=\lambda^{-1/2}J_\varphi v^\ast J_\varphi$.
		We claim that for every $x\in\mathcal{D}$, we have
		\begin{align*}
			&(\mathrm{id}\otimes \tau)\circ\delta_1(x)=R(v)D_c(x\Omega_\varphi)+\lambda^{-1}D_{c^\ast} R(v^\ast)(x\Omega_\varphi)-\lambda^{-1}x\zeta_2,\\
			&(\mathrm{id}\otimes \tau)\circ\delta_2(x)=R(v^\ast)D_{c^\ast}(x\Omega_\varphi)+\lambda D_c R(v)(x\Omega_\varphi)-\lambda x\zeta_1.
		\end{align*}
		To prove the claim, first let $x\in \mathbb C\langle G\rangle\cap M^\ph$ and write $ \delta_{c^\ast} (x)=\sum_j{z_j\otimes w_j}$ as a finite sum. By the KMS condition, $\ph(v^\ast w_j)=\ph(w_j \sigma_{-i}^\ph(v^\ast))=\lambda^{-1}\ph(w_jv^\ast)$. Then we obtain
		\begin{align*}
			(\mathrm{id}\otimes  \tau)\circ\delta_1(x)&
			={E_\varphi}\left(((\mathrm{id}\otimes \ph)\delta_c(x))v+\sum_j \ph(v^\ast w_j)z_j\right)\\
			&={E_\varphi}\left(R(v)D_c(x\Omega_\varphi)+\lambda^{-1}\sum_j \ph( w_jv^\ast )z_j\right)\\
			&={E_\varphi}\left(R(v)D_c(x\Omega_\varphi)+\lambda^{-1}((\mathrm{id}\otimes \ph)( \delta_{c^\ast} (x)v^\ast)\right)\\
			&={E_\varphi}\left(R(v)D_c(x\Omega_\varphi)+\lambda^{-1}(\mathrm{id}\otimes \ph)( \bar\delta_{c^\ast} (xv^\ast)-x \bar\delta_{c^\ast} (v^\ast))\right)\\
			&={E_\varphi}\left(R(v)D_c(x\Omega_\varphi)+\lambda^{-1}D_{c^\ast}(R(v^\ast)x\Omega_\varphi)-\lambda^{-1}x D_{c^\ast}(v^\ast\Omega_\varphi)\right),
		\end{align*}
		where for the last two equalities, we have used  $v^\ast, xv^\ast\in \dom(\bar\delta_{c^\ast})$, which follows from Lemma \ref{lem:polar-part-domain}. 
		The same argument proves the second formula for $(\mathrm{id}\otimes \tau)\circ\delta_2$. Since both formulas have now been established for $x\in\mathbb C\langle G\rangle\cap M^\varphi$, it remains to extend them to $x\in\mathcal D$.
		By \cite[Theorem~4.8]{Nel17}, there exists a sequence $x_n\in\mathbb C\langle G\rangle\cap M^\varphi$ such that $\|x_n\|\leq\|x\|$ for all $n$ and $x_n\to x$ in the graph norm of $\widetilde{\delta}_j$. Passing to the limit yields the desired formulas.
		
		By Lemma \ref{lem:nelson-cancellation}, $1\otimes1  \in \dom(\delta_j^\ast)$. We now give an explicit formula for $\xi_j:=\delta_j^\ast(1\otimes1)$.  Corollary \ref{cor:boundedconjugates} yields an eigenoperator $a$  in $M_\infty$ with eigenvalue $\lambda$ such that $D_{c}^\ast=R(a)-D_{c^\ast}$. Hence
		\[R(a)^\ast\Omega_\varphi=(J_\varphi\sigma_{-i/2}^\varphi(a^\ast)J_\varphi)^\ast\Omega_\varphi=(J_\varphi\sigma_{i/2}^\varphi(a)J_\varphi)\Omega_\varphi=\lambda^{-1} a^\ast\Omega_\varphi.\] 
		A similar calculation yields the identity $R(v)^\ast\Omega_\varphi=\lambda^{-1}v^\ast\Omega_\varphi$.
		Thus for $x\in \mathcal{D}$,  using  $D_{c^\ast}^\ast\Omega_\varphi=R(a)^\ast\Omega_\varphi-D_c\Omega_\varphi=\lambda^{-1}a^\ast\Omega_\varphi $, $v\Omega_\varphi=\lambda^{-1/2}\Omega_\varphi v$ and $J_\tau\zeta_2=\zeta_2$, we have
		\begin{align*}
			\langle x\Omega_\varphi,\delta_1^*(1\otimes1)\rangle_\varphi
			&=\langle R(v)D_c(x\Omega_\varphi)+\lambda^{-1}D_{c^\ast} R(v^\ast)(x\Omega_\varphi)-\lambda^{-1}x\zeta_2,  \Omega_\varphi\rangle_{\varphi}\\
			&=\langle x\Omega_\varphi, D_c^\ast(R(v)^\ast\Omega_\varphi)\rangle+\lambda^{-1}\langle xv^\ast\Omega_\varphi, D_{c^\ast}^\ast\Omega_\varphi\rangle_\varphi-\lambda^{-1}\langle\zeta_2, x^\ast\rangle_{\tau}\\
			&=\lambda^{-1}\Big\langle x\Omega_\varphi,\,
			D_c^*(v^*\Omega_\varphi)+a^*v\Omega_\varphi-\zeta_2
			\Big\rangle_\varphi\\
			&=\lambda^{-1}\Big\langle x\Omega_\varphi,\,
			(v^*a+a^*v)\Omega_\varphi-2\zeta_2
			\Big\rangle_\varphi,
		\end{align*}
		where the last equality uses
		$D_c^*(v^*\Omega_\varphi)=v^*a\Omega_\varphi-\zeta_2$.
		This shows that $\xi_1=\lambda^{-1}(v^\ast a+a^\ast v)\Omega_\varphi-2\lambda^{-1}\zeta_2$. A similar calculation would show that $\xi_2=(va^\ast +av^\ast)\Omega_\varphi-2\lambda\zeta_1$.  This also shows that
		\begin{align*}
			&(\mathrm{id}\otimes \tau)\circ\delta_1(x)=R(v)D_c(x\Omega_\varphi)+\lambda^{-1}D_{c^\ast} R(v^\ast)(x\Omega_\varphi)-\frac{1}{2}\lambda^{-1}x(v^\ast a+a^\ast v)\Omega_\varphi+\frac{1}{2}x\xi_1,\\
			&(\mathrm{id}\otimes \tau)\circ\delta_2(x)=R(v^\ast)D_{c^\ast}(x\Omega_\varphi)+\lambda D_c R(v)(x\Omega_\varphi)-\frac{1}{2}x(va^\ast +av^\ast)\Omega_\varphi+\frac{1}{2}x\xi_2.
		\end{align*}
		Write 
		\begin{align*}
			&Q_1=R(v)D_c+\lambda^{-1}D_{c^\ast} R(v^*)-\frac{1}{2}\lambda^{-1}R(v^*a+a^*v))\\
			&Q_2=R(v^*)D_{c^\ast}+\lambda D_c R(v)-\frac{1}{2}R(va^*+av^*)
		\end{align*}
		so that we get
		\[
		(\mathrm{id}\otimes\tau)\delta_j(x)\Omega_\varphi
		=Q_j(x\Omega_\varphi)+\frac12 x\xi_j,
		\qquad x\in \mathcal{D}.
		\]
		For  $x\in\mathcal{D}$, since $\delta_j(x)^\dagger=\delta_j(x^\ast)$, it follows that $(\tau\otimes \mathrm{id})\circ\delta_j(x^\ast)=((\mathrm{id}\otimes\tau)\delta_j(x))^\ast$. 
		Consequently,
		\[
		(\tau\otimes\mathrm{id}-\mathrm{id}\otimes\tau)\delta_j(x)\Omega_\varphi
		=J_\tau Q_jJ_\tau(x\Omega_\varphi)-Q_j(x\Omega_\varphi)
		-\frac12[x,\xi_j].
		\]
		Therefore, taking the real part in the identity of part (ii) of  Lemma \ref{lem:nelson-cancellation} provides the following equation: for  $x\in \mathcal{D}$ with $\tau(x)=0$,
		\begin{align*}
			\|[x,U_2]\|_2^2&=\mathrm{Re}\langle\bar\delta_0(x), [x,1\otimes1]\rangle+\mathrm{Re}\sum_{j=1}^2 \left\langle [x,H_j], [x, \xi_j]\right\rangle+2\mathrm{Re}\sum_{j=1}^2
			\langle (\tau\otimes\mathrm{id}-\mathrm{id}\otimes\tau)\circ\delta_j(x),[x,H_j]\rangle\\
			&=\mathrm{Re}\langle\bar\delta_0(x), [x,1\otimes1]\rangle+\mathrm{Re}\sum_{j=1}^2 \left\langle [x,H_j], [x, \xi_j]\right\rangle+\mathrm{Re}\sum_{j=1}^2
			\langle 2J_\tau Q_jJ_\tau(x)-2Q_j(x)-[x,\xi_j],[x,H_j]\rangle\\
			&=\mathrm{Re}\langle\bar\delta_0(x), [x,1\otimes1]\rangle+2\mathrm{Re}\sum_{j=1}^2
			\langle J_\tau Q_jJ_\tau(x)-Q_j(x),[x,H_j]\rangle.
		\end{align*}
		This implies that
		\begin{align*}
			\lambda\|x\|_2^2\leq 
			\|[x,U_2]\|_2^2\leq 2\|\bar\delta_0(x)\|_2\,\|x\|_2 +4\|x\|_2\sum_{j=1}^2
			\|Q_j\|\,\|[x,H_j]\|_2
		\end{align*}
		and hence
		\[\|x\|_2\leq \frac{2}{\lambda}\|\bar\delta_0(x)\|_2+ C\sum_{j=1}^2\|[x,H_j]\|_2\,\]
		for some constant $C>0$.
		Now fix $f\in \ker(\bar\delta_0)$ with $\tau(f)=0$. Since $\mathcal{D}$ is a core for $\bar\delta_0$, there exists a sequence $f_n\in \mathcal{D}$ with $\tau(f_n)=0$ such that $f_n\to f$ and $\bar\delta_0(f_n)\to 0$ in $L^2$-norm. The above calculation  therefore yields
		\begin{align*}
			\|f\|_2=\lim_{n\to\infty}\|f_n\|_2\leq C\lim_{n\to\infty}\sum_{j=1}^2\|[f_n,H_j]\|_2=C\sum_{j=1}^2\|[f,H_j]\|_2,
		\end{align*}
		which finishes the proof.
	\end{proof}

	\begin{proposition}\label{prop:fullgap}
		For any $c\in G$ with modular eigenvalue $\lambda\in (0,1)$, 
		there are self-adjoint balanced polynomials $s_1,s_2\in  \mathbb C\langle c^\ast c,cc^\ast\rangle$
		and $C>0$ such that
		\begin{equation*}
			\norm{f-\tau(f)1}_2\leq C\sum_{i=1}^2\norm{[f,s_i]}_2,
			\qquad f\in L^2({M^\ph},\tau).
		\end{equation*}
	\end{proposition}
	\begin{proof}
		Put $N_0=W^*( H_1, H_2)$ and
		$\mathcal K=\ker\overline{\delta_0}\subseteq L^2(M^\varphi,\tau)$, where $H_1, H_2$ and $\delta_0$ are as in Lemma \ref{lem:an estimate on L2-norm in terms of H1,H2}. By the same lemma, we have $L^2(N_0,\tau)\subseteq \mathcal K$. Let $C_0$ denote the constant appearing in Lemma \ref{lem:an estimate on L2-norm in terms of H1,H2}.
		Since $H_j\in W^*(H_j^2)$, the continuous functional calculus allows us to choose real polynomials $p_j$ such that
		\[
		2C_0\sum_{j=1}^2\|H_j-s_j\|<\frac12,
		\]
		where $s_j:=p_j(H_j^2)$.
		For  $f\in \mathcal K$ with $\tau(f)=0$, applying
		Lemma \ref{lem:an estimate on L2-norm in terms of H1,H2}  again gives
		\begin{align*}
			\|f\|_2\leq C_0\sum_{j=1}^2\|[f, H_j]\|_2\leq C_0\sum_{j=1}^2\|[f,s_j]\|_2+2C_0\sum_{j=1}^2
			\|H_j-s_j\|\|f\|_2 \leq C_0\sum_{j=1}^2\|[f,s_j]\|_2+\frac{1}{2}\|f\|_2,
		\end{align*}
		and hence
		\[\|f\|_2\leq 2C_0\sum_{j=1}^2\|[f,s_j]\|_2.\]
		In particular, $N_0$ is a factor. Indeed, if $x\in \mathcal{Z}(N_0)$ is in the center, then $x\in\mathcal K$ and $[x, s_j]=0$ for $j=1,2$.
		Applying the preceding estimate to $f:=x-\tau(x)1$ yields $x=\tau(x)1$ and hence $\mathcal Z(N_0)=\mathbb C1$.
		Since $H_1$ is atomless by \cite[Corollary 5.10]{Nel17}, it follows that $N_0$ is infinite-dimensional. Thus $N_0$ is a $\mathrm{II}_1$ factor. Moreover, the preceding estimate shows that $\{s_1,s_2\}$ is a non-$\Gamma$ set for $N_0$ in the sense of \cite[Definition~2.8]{DI16}. By
		\cite[Definition 2.4, Lemma~2.10]{DI16}, it is a non-amenability set
		for $N_0$, that is, there is a constant $C_1>0$ such that
		\[
		\|\eta\|_2\leq C_1\sum_{j=1}^2\|s_j\eta-\eta s_j\|_2,
		\qquad
		\eta\in L^2(N_0,\tau)\otimes_2 L^2(N_0^{\mathrm{op}},\tau).
		\]
		We next extend this estimate to vectors in $\mathcal K^\perp$ and subsequently to vectors in $L^2(M^\varphi,\tau)$. 
		Set
		\[a_j=s_j\otimes1-1\otimes s_j^{\mathrm{op}}\in N_0\overline{\otimes}N_0^{\mathrm{op}}\]
		which is a self-adjoint element.
		In the coarse representation of $N_0\bar\otimes N_0^{\mathrm{op}}$ on $L^2(N,\tau)\otimes L^2(N^{\mathrm{op}}\tau)$,
		the operator $a_j$ acts by $\eta\mapsto s_j\eta-\eta s_j$.  The previous estimate then gives
		\[
		\sum_{j=1}^2\langle\eta, a_j^2\eta\rangle=\sum_{j=1}^2\|s_j\eta-\eta s_j\|_2^2\geq \frac{1}{2C_1^2}\,\|\eta\|^2,\qquad\eta\in L^2(N_0,\tau)\otimes_2 L^2(N_0^{\mathrm{op}},\tau)
		\]
		that is $a_1^2+a_2^2\geq \frac{1}{2C_1^2}1$ as an operator. The same inequality for $a_1^2+a_2^2$ holds on $ L^2(M^\varphi,\tau)\otimes_2 L^2((M^\varphi)^{\mathrm{op}},\tau)$. Thus we have
		\begin{equation}\label{eq:proof-in-Pro-3.10}
			\|\xi\|_2\leq \sqrt{2}C_1\sum_{j=1}^2\|s_j\xi-\xi s_j\|_2,
			\qquad
			\xi\in L^2(M^\varphi,\tau)\otimes_2 L^2((M^\varphi)^{\mathrm{op}},\tau).
		\end{equation}
		The derivation $\delta_0$  is
		$N_0$-bimodular since $N_0\subseteq \mathcal K:=\ker\bar\delta_0$ and hence $\delta_0(afb)=a\delta_0(f)b$ for 
		$a,b\in\mathbb C\langle H_1, H_2\rangle$ and $f$ in the core. Hence $K$ is an $N_0$-$N_0$ Hilbert bimodule by the same Leibniz rule.
		Let $\bar\delta_0=V(\delta_0^\ast\bar\delta_0)^{1/2}$ be the polar decomposition of the closure $\bar\delta_0$. Bimodularity of $\delta_0$ implies that $V$ is also $N_0$-bimodular, so its restriction to
		$\mathcal K^\perp:=L^2(M^\ph,\tau)\ominus\mathcal K$ is an $N_0$-bimodular isometry into
		$L^2(M^\varphi,\tau)\otimes_2 L^2((M^\varphi)^{\mathrm{op}},\tau)$. Applying \eqref{eq:proof-in-Pro-3.10} with $\xi=Vf$, we obtain
		\[
		\|f\|=\|Vf\|
		\leq \sqrt{2}C_1\sum_{j=1}^2\|V([f,s_j])\|_2
		\leq \sqrt{2}C_1\sum_{j=1}^2\|[f,s_j]\|_2,
		\qquad f\in\mathcal K^\perp.
		\] 
		Combining this with the estimate in 
		Lemma~\ref{lem:an estimate on L2-norm in terms of H1,H2} for $f\in\mathcal K$ yields
		\begin{align*}
			\|f-\tau(f)1\|\leq C_2\sum_{j=1}^2\|[f, s_j]\|_2,\qquad f\in L^2(M^\varphi,\tau)
		\end{align*}
		for some constant $C_2>0$. 
	\end{proof}

	\section{Simplicity of \texorpdfstring{$C^\ast$}{C*}-algebras generated by eigenoperators}\label{sec:simplicity of C-algebras with eigenoperators}
	In this section, we establish additional spectral gap properties for $C^\ast$-algebras generated by eigenoperators. The following result shows that the Dixmier property of the centralizer passes to the whole algebra.
	
	\begin{theorem}\label{thm:sim_centra_implies_simple}
		Let $(M,\varphi)$ be a von Neumann algebra equipped with a normal faithful state $\varphi$ and let $\mathcal A\subset M$ be a unital $C^*$-algebra generated by eigenoperators of $\sigma^\varphi$.
		\begin{enumerate}
			\item[(i)] If $\mathcal A^\varphi$ is simple, then $\mathcal A$ is simple.
			\item[(ii)]  Suppose that for every $b\in\mathcal A^\varphi$,
			$\varphi(b)1\in
			\overline{\operatorname{conv}}^{\|\cdot\|}
			\{ubu^*:u\in\mathcal U(\mathcal A^\varphi)\}.$
			Then for every $a\in\mathcal A$,
			$\varphi(a)1\in
			\overline{\operatorname{conv}}^{\|\cdot\|}
			\{uau^*:u\in\mathcal U(\mathcal A^\varphi)\}.$
			In particular, $\mathcal A$ is simple and has the Dixmier averaging property.
		\end{enumerate}
	\end{theorem}
	\begin{proof}
		By working with $(W^\ast(\mathcal A), \varphi)$, we may assume that $\mathcal A$ is  dense in $M$ in the strong operator topology. It follows that $\mathcal A^\varphi$ is strongly dense in $M^\varphi$.
		
		(i) Assume that $\mathcal A^\varphi$ is simple. We use standard averaging arguments to prove simplicity of $\mathcal A$.
		For each $u\in \mathcal{U}(\mathcal A^\varphi)$, let $\mathrm{Ad}(u)$ denote the inner automorphism of $M$ given by $x\mapsto uxu^\ast$. Consider the family of averaging unital completely positive maps
		\[
		\mathcal T=\operatorname{conv}
		\{\operatorname{Ad}(u):u\in\mathcal U(\mathcal A^\varphi)\},
		\]
		where $\mathrm{conv}$ denotes the convex hull.
		Let $y_1,\ldots,y_n\in \mathcal A$ be modular eigenoperators with respective eigenvalues
		$\mu_1,\ldots,\mu_n$, none of which equals $1$. We show that 
		\begin{align}\label{eq:0 belong to comvex hull}
			0\in \overline{\mathrm{conv}}^{\|\cdot\|_2}\{(Ty_1, \cdots, Ty_n): T\in \mathcal{T}\}\subseteq L^{2}(M,\varphi)^{\oplus n}.
		\end{align}
		Let $\xi\in L^{2}(M,\varphi)^{\oplus n}$ be the unique vector of minimal
		norm in the closed convex hull.   Ultraweak compactness of the unit ball of $M$ yields elements $z_1,\ldots,z_n\in M$ with $\|z_j\|\leq \|y_j\|$ such that $\xi=(z_1\Omega_\varphi, \ldots,z_n\Omega_\varphi)$.  The uniqueness of $\xi$ indicates that $uz_j=z_ju$ for all $u\in \mathcal{U}(\mathcal A^\varphi)$, that is, $z_j\in (M^\varphi)'\cap M$. Moreover, since each $uy_ju^\ast$ is an eigenoperator with eigenvalue $\mu_j$ for $u\in M^\varphi$, it follows that each $z_j$ is also  an eigenoperator with eigenvalue $\mu_j$.

		Write $z_j=v_j|z_j|$, with initial and final projections $p_j=v_j^*v_j$ and $q_j=v_jv_j^*$.  Both projections lie in $M^\varphi$ (see \cite[Proposition 2.6]{Nel17}). Since $z_j\in (M^\varphi)'\cap M$, we have  $v_j=v_jp_j=p_jv_j$, which yields $q_j\le p_j$. Applying the same argument to $v_j^*$ yields  $p_j\le q_j$ and then $p_j=q_j$. Since $\sigma_{-i}^\varphi(v_j)=\mu_j v_j$ by \cite[Proposition 2.6]{Nel17},  the KMS condition gives 
		\[
		\varphi(q_j)=\varphi(v_j^\ast\sigma_{-i}^\varphi(v_j))=\mu_j\varphi(p_j).
		\]
		Since $\mu_j\ne1$, we obtain $\varphi(p_j)=0$. The faithfulness of $\varphi$ then gives $z_j=0$.

		Now let $J$ be a proper closed two-sided ideal of $\mathcal A$, and let
		$\pi\colon \mathcal A\to \mathcal A/J$ be the quotient map. Since $\mathcal A^\varphi$ is simple, the restriction of $\pi$ to $\mathcal A^\varphi$ is injective.
		The state $\pi(b)\mapsto\varphi(b)$ on $\pi(\mathcal A^\varphi)$ therefore
		extends to a state $\psi$ on $\mathcal A/J$.
		For each  $T\in \mathcal{T}$, put $\omega_T=\psi\circ\pi\circ T$. We prove that 
		\begin{align}\label{eq:estimate on phi(x)}
			\varphi(x)\in \overline{\{\omega_T(x): T\in \mathcal{T}\}},\qquad x\in \mathcal{A}.\end{align}
		This will finish the proof because $\omega_T$ vanishes on $J$ for each $T\in\mathcal{T}$, as  $T(J)\subseteq J$. It will then follow that $\varphi$ vanishes on $J$, and the faithfulness of $\varphi$  forces $J$ to be trivial. Since any word of eigenoperators of $\sigma^\varphi$ in $\mathcal A$ is again an eigenoperator, the set of finite sums of eigenoperators is norm dense in $\mathcal A$ by the hypothesis. So it suffices to prove \eqref{eq:estimate on phi(x)} for $x$ of the form $x=b+y_1+\ldots+y_n$, where each $y_j$ is an eigenoperator of $\sigma^\varphi$ with eigenvalue $\mu_j\neq 1$  and $b$ is an eigenoperator with eigenvalue $1$. Fix such $x$ and let $T\in\mathcal{T}$. Then
		$b\in \mathcal A^\varphi$ and  $T(b)\in \mathcal A^\varphi$, so
		$\omega_T(b)=\varphi(T(b))=\varphi(b)$. For $1\leq j\leq n$, it follows from $\varphi(y_j)=\varphi(\sigma_{-i}^\varphi(y_j))=\mu_j \varphi(y_j)$ that $T(y_j)^*T(y_j)\in \mathcal A^\varphi$ and
		$\varphi(y_j)=0$. By the Cauchy-Schwarz inequality,
		\[
		|\omega_T(y_j)|^2
		\leq \psi\bigl(\pi(T(y_j)^*T(y_j))\bigr)
		=\varphi(T(y_j)^*T(y_j))
		=\|T(y_j)\|_{2}^2.
		\]
		By \eqref{eq:0 belong to comvex hull}, there exists a net $T_k\in \mathcal{T}$ such that $\|T_k(y_j)\|_2\to 0$ for every $j=1,\ldots,n$. The preceding estimate then gives $\omega_{T_k}(y_j)\to 0$. Since $\varphi(y_j)=0$  and $\omega_{T_k}(b)=\varphi(b)$, it follows that $\omega_{T_k}(x)\to\varphi(x)$, as required. Applying the same argument simultaneously to finitely many elements of $\mathcal{A}$ shows that $\varphi\in \overline{\{\omega_T: T\in \mathcal{T}\}}^{\mathrm{wk}^\ast}$.
		
		\medskip
		
		(ii) Assume that $\mathcal A^\varphi$ has the Dixmier averaging property stated in the hypothesis. As in (i), it suffices to show that for every eigenoperator $y\in \mathcal A$, $\varphi(y)1$ can be approximated in norm  by  elements of the form $T(y)$ with $T\in \mathcal{T}$.  Fix such $y$ with an eigenvalue $\mu$. If $\mu=1$, then $y\in\mathcal A^\varphi$, and so the hypothesis gives the required norm approximation for  $\varphi(y)1$. Assume that $0<\mu<1$. Then $y^\ast y, yy^\ast\in\mathcal A^\varphi$ and the KMS condition gives
		\[
		\varphi(yy^*)=\mu\varphi(y^*y)\leq \mu\|y\|^2.
		\]
		By the Dixmier averaging property of $A^\varphi$, there is $T_1\in\mathcal{T}$ such that
		\[
		\|T_1(yy^*)-\varphi(yy^*)1\|\leq  \frac{1-\mu}{2}\|y\|^2.
		\]
		Since $T_1$ is completely positive, by the Kadison-Schwarz inequality, $T_1(y)T_1(y)^\ast\leq T_1(yy^\ast )$. This combined with the above inequality yields
		\[
		\|T_1(y)\|^2\leq \|T_1(yy^*)\|
		\leq \frac{1+\mu}{2}\|y\|^2.
		\]
		Since $T_1(y)$ is also an eigenoperator with the same eigenvalue $\mu$, we may iterate this argument to obtain $T_2,\ldots, T_n$ for $n\geq1$ such that $\|T_n\ldots T_1(y)\|\leq \left(\frac{1+\mu}{2}\right)^{n/2}\|y\|$. As $\mu<1$ and $\varphi(y)=0$, it follows that $\|T_n\ldots T_1(y)-\varphi(y)1\|\to 0$ as $n\to\infty$. A similar estimate will yield the result when $\mu>1$. This finishes the proof.
	\end{proof}

	The following proposition gives a rapid decay condition for the Dixmier property of the centralizer of a $C^\ast$-algebra generated by  eigenoperators.

	\begin{proposition}\label{prop:rapid_decay-property}
		Let $(M,\varphi)$ be a von Neumann algebra equipped with a faithful
		normal state $\varphi$ and generated by a finite self-adjoint set $G$ of
		eigenoperators of $\sigma^\varphi$. Put $\mathcal A=C^*(1,G)$ and
		$\tau=\varphi|_{M^\varphi}$. Suppose that
		$x_1,\ldots,x_d\in\mathcal A^\varphi$ are self-adjoint polynomials in
		$G$ such that
		\[
		\|f-\tau(f)1\|_2
		\le C\sum_{j=1}^d\|[x_j,f]\|_2,
		\qquad f\in L^2(M^\varphi,\tau).
		\]
		For each balanced word $g$ in $G$ (see Definition \ref{def:balanced word}), and for each $m\ge0$, put
		\[
		V_m(g)=\mathbb C1+
		\operatorname{span}\{agb:
		a,b\text{ are words in }x_1,\ldots,x_d,\ |a|+|b|\le m\},
		\qquad
		\rho_m(g)=\sup_{0\ne y\in V_m(g)}\frac{\|y\|}{\|y\|_2}.
		\]
		If $\log\rho_m(g)=o(m)$ for every balanced word $g$ in $G$, then
		$\mathcal A^\varphi$ has the Dixmier averaging property. More precisely,
		for every $a\in\mathcal A^\varphi$,
		\[
		\tau(a)1\in
		\overline{\operatorname{conv}}^{\|\cdot\|}
		\{uau^*:u\in\mathcal U(C^*(1,x_1,\ldots,x_d))\}.
		\]
		In particular, $\mathcal A^\varphi$ is simple and
		$\tau|_{\mathcal A^\varphi}$ is its unique tracial state.
	\end{proposition}
	\begin{proof}
		It suffices to establish the stated averaging property, since it implies simplicity and uniqueness of the tracial state.
		By rescaling the $x_j$, we may assume that $\|x_j\|\le 1/2$ for all $1\le j\le d$.
		Put $\mathcal B=C^*(1,x_1,\ldots,x_d)$, and
		$u_j=e^{ix_j}$. To avoid confusion, we shall use $\|\cdot\|_\infty$ for operator norm and $\|\cdot\|_2$ for the $L^2$-norm for elements in $M^\varphi$. Define the averaging operator $T$ on $M^{\varphi}$ by
		\[
		T(y)=\frac1{2d}\sum_{j=1}^d
		\bigl(u_jyu_j^*+u_j^*yu_j\bigr),\qquad y\in M^{\varphi}.
		\]
		Note that $T$ extends to an operator on $L^2(M^\varphi,\tau)$. We shall treat $T$ both as an operator on $M^\varphi$ and on $L^2(M^\varphi,\tau)$, and we shall make this distinction clear.
		Define the operators $A_j$ on  $L^2(M^{\varphi},\tau)$ by   $A_j\eta=[x_j,\eta]=x_j\eta-\eta x_j$, that is, $A_j=x_j-J_\tau x_j J_\tau$, where $J_\tau$ is the conjugation anti-unitary on $L^2(M^\varphi,\tau)$. Since each $x_j$ is self-adjoint, each $A_j$ is also self-adjoint with $\|A_j\|\le1$. 
		Since $x_j$ and $J_\tau x_j J_\tau$ commute, it follows that 
		$e^{itA_j}=e^{itx_j}e^{-itJ_\tau x_jJ_\tau}=e^{itx_j}J_\tau e^{itx_j}J_\tau$ for $t\in\mathbb R$. Equivalently, $e^{it A_j}\eta=e^{itx_j}\eta e^{-itx_j}$ for $\eta\in L^2(M^{\varphi},\tau)$. Then viewed as an operator on $L^2(M^\varphi,\tau)$, $T$ can be written as
		\[T=\frac{1}{d}\sum_{j=1}^d\frac{ e^{iA_j}+e^{-iA_j}}{2}=\frac{1}{d}\sum_{j=1}^d\cos(A_j).\]
		Since $\cos(t)\geq\cos(1)$ for all $t\in [-1,1]$, we have $T\geq \cos(1)\mathrm{id}$  on $L^2(M^\varphi,\tau)$.
		Using $1-\cos t\ge t^2/4$ for $|t|\le1$ and the
		spectral gap hypothesis, we obtain
		\[
		\langle y,(\mathrm{id}-T)y\rangle_\tau\geq \frac{1}{4d}\sum_{j=1}^d\langle y, A_j^2y\rangle_\tau= 
		\frac1{4d}\sum_{j=1}^d\|[x_j,y]\|_2^2
		\ge \frac1{4d^2C^2}\|y\|_2^2
		\]
		for $y\in M^{\varphi}\ominus\mathbb C1:=\{z\in M^\varphi; \tau(z)=0\}$. Thus $T\leq (1-\frac{1}{4d^2C^2})\mathrm{id}$ as an operator on $L^2(M^\varphi,\tau)\ominus\mathbb C$. Therefore, there is a constant $c>0$ such that
		\begin{align}\label{eq:Tn L2-norm}
			\|T^n(\xi)\|_2\le e^{-cn}\|\xi\|_2,
			\qquad \xi\in L^2(M^{\varphi},\tau)\ominus\mathbb C1.  
		\end{align}
		We now approximate $T^n$ by truncating the Taylor series of the full iterate at total degree $m$. Regard $A_j(y)=[x_j,y]$ as a bounded
		operator on $M^{\varphi}$ with its operator norm $\|\cdot\|_{\mathrm{op}}$. Then we still have $\|A_j\|_{\mathrm{op}}\le1$.
		For $z\in\mathbb C$, set
		\[
		F(z)=\frac1{2d}\sum_{j=1}^d
		\bigl(e^{izA_j}+e^{-izA_j}\bigr).
		\]
		Then $F(1)=T$ and $\|F(z)\|_{\mathrm{op}}\leq e^{|\mathrm{Im\,} z|}$. For $n,m\geq 1$, write
		\[
		F(z)^n=\sum_{k=0}^\infty z^k T_{n,k}\quad \text{and}
		\quad
		S_{n,m}=\sum_{k=0}^m T_{n,k}
		\]
		for some bounded operators $\{T_{n,k}\}_{k\geq 0}$ on $M^\varphi$.
		The operator $T_{n,k}$ is a linear combination of maps $y\mapsto ayb$, where $a$ and $b$ are words in the $x_j$ with $|a|+|b|=k$. Indeed, $T_{n,k}$ is a linear combination of operators of  the form $A_{j_1}\ldots A_{j_k}$ with $1\leq j_1,\ldots,j_k\leq d$, and  each occurrence of $A_j$ acts by the difference of left and right multiplication by $x_j$.
		This also implies that $T_{n,k}1=0$ for all $k\geq 1$. Since $T_{n,0}=\mathrm{id}$, we have $S_{n,m}1=1$ for all $n,m\geq 1$. In particular, for any balanced word $g$, $S_{n,m}(g)$ is a linear combination of terms $agb$, where $a$ and $b$ are words with $|a|+|b|\leq m$. Hence, for $y=g-\tau(g)1$, $S_{n,m}(y)\in V_m(g)$.

		On the circle $|z|=2$,  Cauchy’s formula gives
		\[
		T_{n,k}=\frac{1}{2\pi i}\int_{|z|=2}\frac{F(z)^n}{z^{k+1}}\,dz.
		\] 
		Since $\|F(z)^n\|_{\mathrm{op}}\leq e^{n|\mathrm{Im\,} z|}\leq e^{2n}$ on $|z|=2$, it follows that
		$\|T_{n,k}\|_{\mathrm{op}}\le e^{2n}2^{-k}$, and consequently
		\begin{align}\label{eq:estimate on Tn-Snm}
			\|T^n-S_{n,m}\|_{\mathrm{op}}\leq \sum_{k=m+1}^\infty\|T_{n,k}\|_{\mathrm{op}}\leq \sum_{k=m+1}^\infty e^{2n}2^{-k}
			=e^{2n}2^{-m}.
		\end{align}
		One could use Cauchy's estimate for any $r>1$ instead of $2$ to get a decay factor $r^{-k}$.
		Choose $m_n=\lceil bn\rceil$, where $b>0$ is a constant with $b\log2>2+c$ and $c$ is the constant in \eqref{eq:Tn L2-norm}. Put
		$\varepsilon_n=e^{2n}2^{-m_n}$. Then 
		$\varepsilon_n\le e^{-cn}$.  Since $S_{n,m_n}(y)\in V_{m_n}(g)$,  the hypothesis and \eqref{eq:estimate on Tn-Snm} yield
		\begin{align*}
			\|T^n(y)\|_\infty\leq \|S_{n,m_n}(y)\|_\infty+\|(T^n-S_{n,m_n})(y)\|_\infty\le \rho_{m_n}(g)\|S_{n,m_n}(y)\|_2+\varepsilon_n\|y\|_\infty.
		\end{align*}
		However, since $y\in M^\varphi\ominus\mathbb C1$, \eqref{eq:Tn L2-norm} and \eqref{eq:estimate on Tn-Snm} give
		\begin{align*}
			\|S_{n,m_{n}}(y)\|_2\leq \|T^n(y)\|_2+\|(T^n-S_{n,m_n})(y)\|_2\leq \|T^n(y)\|_2+\|(T^n-S_{n,m_n})(y)\|_\infty\leq e^{-cn}\|y\|_2+\varepsilon_n\|y\|_\infty.
		\end{align*}
		Using $\varepsilon_n\leq e^{-cn}$,
		we finally get
		\[\|T^n(g)-\tau(g)1\|_\infty=\|T^n(y)\|_\infty\leq e^{-cn}\rho_{m_n}(g)\left(\|y\|_2+\|y\|_\infty\right)+e^{-cn}\|y\|_\infty.\]
		By hypothesis, $\log\rho_{m_n}(g)=o({m_n})=o(n)$, hence $\rho_{m_n}(g)e^{-cn}\to 0$ as $n\to\infty$. Therefore, $T^n(g)\to\tau(g)1$ in norm for every balanced word $g$. Note that $T^n(g)\in \mathrm{conv}\{ugu^\ast: u\in \mathcal{U}(\mathcal{B})\}$, thus we have proved that $\tau(g)1\in \overline{\mathrm{conv}}^{\|\cdot\|}\{ugu^\ast: u\in \mathcal{U}(\mathcal{B})\}$.
		Since the balanced words span a norm-dense subspace of $\mathcal A^\varphi$ by Lemma \ref{lem:dense elements in centralizer}, 
		the result follows.
	\end{proof}
	
	As a corollary of Proposition \ref{prop:rapid_decay-property}, we recover the following result of \cite{AJW26}.
	
	\begin{corollary}[{\cite[Theorem~A]{AJW26}}]\label{cor:ajw}
		Let $\mathcal A=C^*(1,x_1,\ldots,x_d)$ be a unital $C^*$-algebra
		with a faithful tracial state $\tau$. Suppose that there is $C>0$
		such that
		\[
		\|f-\tau(f)1\|_2
		\leq C\sum_{j=1}^d
		\bigl(\|[f,x_j]\|_2+\|[f^*,x_j]\|_2\bigr),
		\qquad f\in\mathcal A.
		\]
		For $m\geq 0$, put
		\[
		\rho(m)=\sup\left\{
		\frac{\|p(x)\|}{\|p(x)\|_2}:
		p\text{ is a $*$-polynomial of degree at most }m,\quad p(x)\ne0
		\right\}.
		\]
		If $\log\rho(n^2)=o(n)$, then $\mathcal A$ has the Dixmier averaging
		property. In particular, $\mathcal A$ is simple and $\tau$ is its
		unique tracial state.
	\end{corollary}
	
	\begin{proof}
		By considering the real and imaginary parts of each $x_j$ for the spectral gap hypothesis, we may as well assume that each $x_j$ is self-adjoint.
		Since $\tau$ is a trace, every word in the set 
		$G=\{x_1,\ldots,x_d\}$ is balanced.
		Fix a word $g$ in $G$ of length $k$. Consider the set $V_m(g)$ and $\rho_m(g)$ as in  Proposition \ref{prop:rapid_decay-property}. Then every element of
		$V_m(g)$ is a polynomial of degree at most $m+k$ in $G$ and
		$\rho_m(g)\leq \rho(m+k)$.
		By monotonicity of $\rho$, the hypothesis $\log\rho(n^2)=o(n)$
		implies $\log\rho(m)=o(\sqrt m)$. Thus $\log\rho_m(g)=o(m)$ for
		every word $g$ in $G$. Proposition~\ref{prop:rapid_decay-property}
		now applies to $\mathcal A^\tau=\mathcal A$.
	\end{proof}

	\section{Operator-norm estimates for balanced polynomials on \texorpdfstring{$q$}{q}-Araki-Woods algebras}\label{sec:operator_norm_estimate}
	In the remainder of this paper, we specialize to $q$-Araki-Woods algebras ${A}_q(\mathcal H_{\mathbb R}, U)$ with $\dim_{\mathbb R} \mathcal H_{\mathbb R}\geq 2$. If $\mathcal H_{\mathbb R}$ is of dimension $n<\infty$, then for a fixed basis $\{e_k\}_{1\leq k\leq n}$ of $\mathcal H_{\mathbb C}$ consisting of eigenvectors of $A_U$ with eigenvalues $\lambda_k>0$,   the generators $\{W(e_k)\}_{1\leq k\leq  n}$  are eigenoperators of $\sigma^{\varphi_q}$. Indeed, \cite[Lemma 1.4]{Hia03} yields $\sigma_t^{\varphi_q}(W(e_k))=W(U_{-t}e_k)=\lambda_k^{-it}W(e_k)$ for $t\in\mathbb R$. Moreover, the set $\{W(e_k)\}_{1\leq k\leq  n}$ has  finite free Fisher information with Lipschitz conjugate variables.  Miyagawa and Speicher \cite[Theorem 1.1]{MS23} proved this in the tracial case  and the first author \cite[Corollary 3.13]{Kum23} proved it in the nontracial case. Hence, ${A}_q(\mathcal H_{\mathbb R}, U)$ satisfies the hypothesis of Proposition \ref{prop:fullgap} and hence has a commutator spectral gap.  To show  ${A}_q(\mathcal H_{\mathbb R}, U)$ is simple, we establish the operator-norm estimates needed to verify the rapid decay condition in Proposition~\ref{prop:rapid_decay-property}.
	
	\subsection{A tensor product estimate}
	We first prove that the constant in the lower bound of \cite[lemma 1]{Nou04} is independent of $n$ and $k$. This answers the question raised in the remark following that lemma.
	Let $\mathcal H$ be a complex Hilbert space, $q\in (-1,1)$ and let $\mathcal{F}_q(\mathcal{H})=\oplus_{n\geq0}\mathcal{H}^{\otimes_q^n}$ be the $q$-Fock space. We recall some notation from \cite{BS91, Boz98, Nou04, ABW18}.
	For every $k,l\geq 0$, consider the operator $R_{k,l}^\ast:\mathcal{H}^{\otimes_q^{k+l}}\to \mathcal{H}^{\otimes_q^k}\otimes\mathcal{H}^{\otimes_q^l}$, which acts on the simple tensor  
	$\xi=\xi_1\otimes\cdots\otimes\xi_n$, for $\xi_j\in \mathcal{H}$ as
	\[
	R_{k,l}^*\xi
	=\sum_{\substack{J\subseteq[k+l]\\ |J|=k}}
	q^{i(J)}\,\xi_J\otimes_2 \xi_{J^c},
	\qquad
	i(J)=\sum_{r=1}^{k}(j_r-r).
	\]
	Here $J=\{j_1<\cdots<j_k\}\subseteq[k+l]$, $\xi_J=\xi_{j_1}\otimes\cdots\otimes\xi_{j_k}$, and $\xi_{J^c}$ is formed from $J^c=[k+l]\setminus J$ in increasing order. It is known that $R_{k,l}^\ast$ extends to a bounded operator with $\|R_{k,l}^\ast \|\leq C_q$ (\cite[Equation 3]{Nou04}). Also it is easy to verify that the adjoint $R_{k,l}$ is just the tensor concatenation:
	$R_{k,l}(\zeta\ten\eta)=\zeta\ten\eta$ on the algebraic tensors $\mathcal{H}^{\otimes k}\otimes\mathcal{H}^{\otimes l}$.
	
	The following result shows that $R_{k,l}^{-1}$ is also a well-defined operator whose operator norm is bounded by a constant independent of $k,l$. 
	Recall the strictly positive operator $P_q^{(n)}$ on the algebraic tensor product $\mathcal{H}^{\otimes n}$ for each $n\geq 0$ which defines the $q$-inner product.

	\begin{proposition}\label{prop:upper and lower bounds on tensor products}
		There is  a constant $ B_q >0$, depending only on $q$, such that
		\[
		B_q (P_q^{(k)}\ten P_q^{(l)})\leq P_{q}^{k+l}\leq {C_q}(P_q^{(k)}\ten P_q^{(l)})
		\]
		for every $k,l\geq0$. Consequently
		$\norm{R_{k,l}}\leq {C_q}^{1/2}$ and
		$\norm{R_{k,l}^{-1}}\leq B_q ^{-1/2}$.
	\end{proposition}
	\begin{proof}
		For convenience, we write $P_n=P_q^{(n)}$.
		The upper comparison is \cite[Section~2, equation~(3)]{Nou04}. So we will prove the lower comparison.
		For $k=0$ or $l=0$,  the inequality is obvious. Assume that $k,l\geq1$.
		By \cite[Theorem 6]{Boz98} and the proof of \cite[Theorem 3.1]{BS94}, there are constants $b_q$ and $d_q$ depending only on $q$ such that 
		\[
		b_q(I\ten P_{n-1})\leq P_n\leq d_q(I\ten P_{n-1}),\;\;n\geq 1.
		\] 
		The following  inclusion--exclusion calculation is the mechanism of \cite[Lemma~2.6 and proof of Theorem~2.4]{BS94}.
		Set $P_0=Q_0=1$, and for $n\geq1$, put $Q_n=P_n^{-1}$. Define the operator
		$K_n:\mathcal H^{\ten n}\to \mathcal H^{\ten n}$ by
		\[
		K_n=\sum_C(-1)^{|C|}\bigotimes_{B\in\pi_C}Q_{|B|},
		\] 
		where $C\subseteq{[n-1]}$ varies over all subsets. Here $\pi_C$ is the partition of $\{1,\ldots, n\}$ into consecutive blocks,  obtained by placing a cut immediately after each position in  $C$, and $|B|$ denotes the cardinality of each block $B\in\pi_C$. For example, take $n=7$ and $C=\{2,5\}$. The cuts give
		$
		1\,2\mid3\,4\,5\mid6\,7,
		$
		so
		$\pi_C=\bigl\{\{1,2\},\{3,4,5\},\{6,7\}\bigr\}.$
		M\"obius inversion on the lattice of interval partitions (see for example, \cite[Proposition 10.6]{NS06}) gives
		\[Q_n=\sum_\pi\bigotimes_{B\in\pi}K_{|B|},\]
		where $\pi$ runs over interval partitions of $[n]$: every block must consist of consecutive integers. One can also prove this by induction.
		Write $n=k+l$ and $A=P_k\ten P_l$ so that $A^{-1}=Q_k\otimes Q_l$. Any interval partition $\pi$ of $[n]$  consists  either of an arbitrary interval partition of $\{1,\ldots,k\}$ and an arbitrary interval partition of $\{k+1,\ldots,n\}$, or with block containing both $k$ and $k+1$. In  the latter case, the 
		block has the form $B=\{u,u+1,\ldots,v\}$ with $1\le u\le k<v\le n$. The former case yields $Q_k\otimes Q_l$, and the resulting identity is
		\begin{equation*}
			Q_n=A^{-1}+\sum_{1\leq u\leq k<v\leq n}
			Q_{u-1}\ten K_{v-u+1}\ten Q_{n-v}.
		\end{equation*}
		For fixed  $u\leq k<v$, set $L=v-u+1$. Iterating the one-letter upper bound gives
		$A\leq d_q^L(P_{u-1}\ten 1_{H^{\ten L}}\ten P_{n-v})$. Indeed, 
		\[
		P_k
		\leq{d_q}(P_{k-1}\otimes 1_H)\leq d_q^2(P_{k-2}\otimes 1_H^{\otimes2})\leq\cdots\leq d_q^i
		\bigl(P_{u-1}\otimes 1_{H^{\otimes i}}\bigr),\;\;\;i=k-u+1.
		\]
		Similarly, 
		$P_{n-k}\leq d_q^j
		\bigl(1_{H^{\otimes j}}\otimes P_{n-v}\bigr)$ for $j=v-k$. Therefore,
		\[
		A
		\leq d_q^{i+j}
		\bigl(P_{u-1}\otimes 1_{H^{\otimes i}}
		\otimes 1_{H^{\otimes j}}\otimes P_{n-v}\bigr)=d_q^L
		\bigl(P_{u-1}\otimes 1_{H^{\otimes L}}\otimes P_{n-v}\bigr).
		\]
		It follows that
		\[
		\norm{A^{1/2}(Q_{u-1}\ten K_L\ten Q_{n-v})A^{1/2}}
		\leq d_q^L\norm{K_L}.
		\]
		Summing these inequalities gives
		$A^{1/2}Q_nA^{1/2}\leq B_q ^{-1}\mathrm{id}$, where
		\[
		A^{1/2}Q_nA^{1/2}\leq \left(1+\sum_{1\leq u\leq k<v\leq n}d_q^{v-u+1}\|K_{v-u+1}\|\right)\mathrm{id}.
		\]
		We now estimate the norms of $P_nK_n$.
		For $J\subseteq{[n-1]}$, let $W_J\subseteq S_n$ be the
		subgroup generated by $(i\ i+1)$, $i\in J$, and put
		\[
		{D_J}=\{\sigma\in S_n:\sigma(i)<\sigma(i+1)\ \,\,\forall i\in J\}.
		\]
		Every $\sigma\in S_n$ has a unique expression $\sigma=uv$,
		$u\in {D_J}$, $v\in W_J$, with
		$i(\sigma)=i(u)+i(v)$. Note that ${D_J}$ is precisely the set that appears before \cite[Lemma 2.6]{BS94}.
		Consequently, we have
		\[
		P_n=\sum_{\sigma\in S_n}q^{i(\sigma)}U_\sigma=(\sum_{u\in {D_J}}q^{i(u)}U_u)
		(\sum_{v\in W_J}q^{i(v)}U_v).
		\]
		The second factor is the tensor product $\otimes_{B\in\pi_{C}}P_{|B|}$, where $C={[n-1]}\setminus J$. Indeed, it follows from
		$W_J=\prod_{B\in\pi_C}S_B,$
		where $S_B$ denotes the set of permutations of $B$, fixing everything outside $B$. Therefore, 
		\begin{align*}
			P_nK_n=\sum_{J}(-1)^{n-1-|J|}\sum_{u\in {D_J}}q^{i(u)}U_u,
		\end{align*}
		where the sum runs over all subsets $J$ of ${[n-1]}$. If we set
		$A(\sigma)=\{i\in{[n-1]}:\sigma(i)<\sigma(i+1)\}$, then 
		$\sigma\in D_J\iff J\subseteq A(\sigma)$, so we have
		\begin{align*}
			P_nK_n
			=(-1)^{n-1}\sum_{\sigma\in S_n}
			q^{i(\sigma)}U_\sigma
			{\sum_{J\subseteq A(\sigma)}(-1)^{|J|}}.
		\end{align*}
		However, $\sum_{J\subseteq A(\sigma)}(-1)^{|J|}=0$ if $A(\sigma)\neq\emptyset$, while it is equal to $1$ precisely for the permutation $\omega_n$ given by $\omega_n(i)=n+1-i$, $i\in [n]$. Thus we get
		$P_nK_n=(-1)^{n-1}q^{n(n-1)/2}U_{\omega_n}.$
		It follows that $\|P_nK_n\|\leq |q|^{n(n-1)/2}$ and hence
		$\norm{K_n}\leq |q|^{n(n-1)/2}b_q^{-(n-1)}$.
		Thus we finally get
		\begin{align*}
			Q_n \leq   \left(1+\sum_{1\leq u\leq k<v\leq n}d_q^{v-u+1}b_q^{-(v-u)}|q|^{(v-u+1)(v-u)/2}\right)A^{-1}.
		\end{align*}
		Setting
		\[
		B_q^{-1}
		=1+\sum_{L\ge2}(L-1)d_q^L b_q^{-(L-1)}|q|^{L(L-1)/2},
		\]
		we obtain $B_q(P_k\otimes P_l)\leq P_{k+l}$.
	\end{proof}

	\subsection{A Haagerup-type Inequality} For the remainder of this section,
	let $(U_t)_{t\in\mathbb R}$ be a nontrivial  orthogonal representation  on a two-dimensional space $\mathbb R^2$. Our aim is now to estimate the operator norms of certain balanced Wick operators in terms of their $L^2$-norms. 
	
	Since $U$ is nontrivial, there is $0<\lambda<1$ and unit vectors $e,f\in \mathcal{H}_{\mathbb C}$  such that $A_Ue=\lambda^{-1} e$ and $A_Uf=\lambda f$. We have $Ie=\lambda^{1/2}f$, where recall that $I: \mathbb C^2\to  \mathbb C^2$ is the antilinear operator. Note that $W(e)^\ast=\lambda^{1/2}W(f)$.
	
	We assign each of $e$ and $f$ a charge in $\{-1,1\}$.
	Define 
	$\chi:\{e,f\}\to\{-1,1\}$ by  $\chi(e)=1, \chi(f)=-1.$ For a simple tensor $\xi=\xi_1\otimes\ldots\otimes\xi_n$ with
	$\xi_j\in \{e,f\}$,  extend the charge by setting
	\begin{align}\label{eq:height of word}
		\chi(\Omega)=0,\qquad\chi(\xi_1\otimes \ldots \otimes\xi_k)=\sum_{j=1}^k\chi(\xi_j),
		\quad1\leq k\leq n,\qquad\mbox{and }\;\;
		\operatorname{h}(\xi)=\max_{0\leq k\leq n}\chi(\xi_1\otimes\ldots\otimes \xi_k).
	\end{align}
	Call $\operatorname{h}(\xi)$ the {\em height} of $\xi$ which is the
	largest prefix charge. 
	We assign the same charge and height to scalar multiples of these simple tensors.
	Since $\sigma_t^{\varphi_q}(W(e))=W(U_{-t}e)=\lambda^{it} W(e)$ and $\sigma_t^{\varphi_q}(W(f))=\lambda^{-it}W(f)$,  it follows that
	\[\sigma_t^{\varphi_q}(W(\xi))=\lambda^{it\chi(\xi)}W(\xi),\qquad \sigma_t^{\varphi_q}(W(\xi_1)\ldots W(\xi_n))=\lambda^{it\chi(\xi)}W(\xi_1)\ldots W(\xi_n),\;\;\;t\in\mathbb R\]  for the tensor $\xi$. In particular, $\chi(\xi)=0$ if and only if $W(\xi)$ is in the centralizer $M_q^{\varphi_q}$.  We emphasize that in the notation of \cite[Section 5.2]{Nel17}, a word $w=W(\xi_1)\ldots W(\xi_n)$ in $\{W(e),W(e)^\ast\}$ belongs to $E_{\leq1}^\infty$ if and only if each suffix word $
	W(\xi_k)\ldots  W(\xi_n)$ has a modular eigenvalue $\leq 1$, which is equivalent to  $\chi(\xi_1\otimes\ldots\otimes \xi_k)\leq \chi(\xi_1\otimes\ldots \otimes\xi_n)$ for all $1\leq k\leq n$.
	
	A tensor $\xi=\xi_1\otimes\ldots\otimes\xi_n$ is {\em balanced} if  $\chi(\xi)=0$, that is, $e$ and $f$ occur equally often
	in the tensor. This terminology agrees with that used for operators in Section \ref{sec:derivations on free fisher infomration}. In particular, $\xi$ is balanced if and only if $W(\xi)$ and $W(\xi_1)\ldots W(\xi_n)$ are balanced.  
	For $R\geq0$, define $Q_R$ to be the orthogonal projection on
	$\mathcal F_q(\mathcal H_U)$ with range
	\[
	\mathrm{ran\;} Q_R= \overline{\operatorname{span}}^{\,\|\cdot\|_q}
	\bigl\{\xi=\xi_1\otimes\ldots\otimes\xi_k: k\geq0, \xi_j\in \{e,f\},\ \xi\text{ is balanced and }\operatorname{h}(\xi)\leq R\bigr\}.
	\]
	The vacuum vector $\Omega$ belongs to $\ran Q_R$.

	\begin{lemma}\label{lem:heightwick}
		There is a constant $C=C(q,\lambda)$ such that for all $n\in\mathbb N$ and $R>0$,
		\[
		\norm{W(\eta)}\leq C(n+1)\lambda^{-R/2}\norm\eta_q,
		\qquad \eta\in \mathcal H_U^{\otimes_q^ n}\cap\mathrm{ran\;} Q_R.
		\]
	\end{lemma}
	\begin{proof}
		For each $k\geq 1$,
		let $\widehat I_k:\mathcal{H}_{\mathbb C}^{\otimes k}\to\overline{\mathcal{H}}_{\mathbb C}^{\otimes k}$ be the linear map on the algebraic tensor product defined by 
		$\widehat I_k(\zeta)=\overline{I^{\otimes k}\zeta}$. Here $\overline{\mathcal{K}}$ denotes the conjugate Hilbert space for any Hilbert space $\mathcal{K}$, and we identify $(\overline{\mathcal{H}}_{\mathbb C})^{\otimes k}$ with $\overline{\mathcal{H}_{\mathbb C}^{\otimes k}}$.
		The Wick formula \cite[Proposition~2.16]{ABW18} gives
		\[
		W(\eta)=\sum_{k=0}^n
		U_{n-k,k}\bigl((\id\ten\widehat I_k)R_{n-k,k}^*\eta\bigr),\qquad\eta\in \mathcal{H}_{\mathbb C}^{\otimes n} 
		\]
		where $R_{n-k,k}$ is the concatenation map defined above, and $U_{k,l}:\mathcal H_{U}^{\otimes_q^{k}}\otimes \overline {\mathcal H}_U^{\otimes_q^ l}\to B(\mathcal{F}_q(\mathcal H_U))$ is given by
		\[U_{k,l}\bigl((\xi_1\ten\cdots\ten\xi_k)
		\ten{(\bar \zeta_1\ten\cdots\ten\bar\zeta_l)}\bigr)
		={}\ell_q(\xi_1)\cdots\ell_q(\xi_k)
		\ell_q(\zeta_1)^*\cdots\ell_q(\zeta_l)^*.\]   
		It is known that $U_{k,l}$ is bounded with $\|U_{k,l}\|\leq C_q$ (\cite[Lemma 4 and Corollary 1]{Nou04}).
		Fix $0\leq k\leq n$ and set $m=n-k$. In the polar decomposition $I=\mathcal{J}|I|=\mathcal{J}A_U^{-1/2}$, put
		$D_{k}=\id\ten|I|^{\otimes k}$ on $\mathcal H_{\mathbb C}^{\otimes m}\otimes \mathcal H_{\mathbb C}^{\otimes k}$.   Since  $|I|(e)=\lambda^{1/2}e$ and $|I|(f)=\lambda^{-1/2}f$, it follows that   
		\[D_{k}(\xi\otimes\zeta)=\lambda^{s/2}\xi\otimes\zeta,\]
		whenever $\xi\in \mathcal H_{\mathbb C}^{\otimes{m}},\zeta\in \mathcal H_{\mathbb C}^{\otimes k}$ are words in $\{e,f\}$ with the charge $\chi(\zeta)=s$. Moreover, if such a word $\xi\otimes\zeta$ belongs to $\mathrm{ran\; } Q_R$, then it  is balanced by definition with $\chi(\xi)\leq R$ and  $\chi(\xi)+\chi(\zeta)=0$. This implies  $\chi(\zeta)\geq -R$, that is, $\lambda^{\chi(\zeta)/2}\leq \lambda^{-R/2}$. Since vectors of different charge are orthogonal, it follows that whenever $\eta\in \mathcal H_{\mathbb C}^{\otimes n}\cap\mathrm{ran\;} Q_R$, we have
		\[\|D_{k}R_{m,k}^{-1}\eta\|_q^2\leq \lambda^{-R}\|R_{m,k}^{-1}\eta\|_{\mathcal H_U^{\otimes_q^{m}}\otimes \mathcal H_{U}^{\otimes_q ^k}}^2\leq B_q^{-1}\lambda^{-R}\|\eta\|_q^2,\]
		where the last inequality follows from Proposition \ref{prop:upper and lower bounds on tensor products}.
		The standard estimate now gives
		\begin{align*}
			\|W(\eta)\|&\leq  \sum_{k=0}^n \|U_{m,k}(\id\ten \widehat I_k)R_{m,k}^*\eta\|       \leq C_q(n+1)\max_{0\leq k\leq n}\|D_{k}R_{m,k}^*\eta\|,
		\end{align*}
		where we have used $\|U_{m,k}\|\leq C_q$. For any $0\leq s\leq \min\{m, k\}$, write $\operatorname{Flip}_s$
		for interchange of the middle two factors in
		$\mathcal H_U^{\ten_q^{(m-s)}}\ten \mathcal H_U^{\ten_q^ s}\ten
		\mathcal H_U^{\ten_q^ s}\ten \mathcal H_U^{\ten_q^{(k-s)}}
		$ which is a unitary operator. We now show that
		\[
		D_kR_{m,k}^*R_{m,k}D_k^{-1}
		=\sum_{s=0}^{\min\{m,k\}}q^{s^2}
		(R_{m-s, s}\otimes R_{s,k-s})\circ\operatorname{Flip}_s\circ
		\bigl(\id\ten|I|^{\otimes s}\ten(|I|^{-1})^{\otimes s}\ten\id\bigr)\circ(R_{m-s, s}^\ast\otimes R_{s,k-s}^\ast).
		\]
		Let $\xi=\xi_1\cdots\xi_m$ and $\zeta=\zeta_1\cdots\zeta_k$ be words
		in $\{e,f\}$. For $J=\{j_1<\cdots<j_r\}\subseteq[m]$, put
		$\xi_J=\xi_{j_1}\cdots\xi_{j_r}$.
		Similarly, define $\xi_{J^c}$ for $J^c=[m]\setminus J$.
		The empty
		subword is interpreted as the vacuum vector. Therefore, we have
		\begin{align*}
			R_{m,k}^*R_{m,k}(\xi\ten\zeta)=R_{m,k}^\ast(\xi\otimes\zeta)=\sum_{s=0}^{\min\{m,k\}}
			\sum_{\substack{J\subseteq[m],\ |J|=m-s\\
					L\subseteq[k],\ |L|=s}}
			q^{s^2+i(J)+i(L)}
			(\xi_J\zeta_L)\ten(\xi_{J^c}\zeta_{L^c}),
		\end{align*}
		which implies that 
		\begin{align*}
			D_kR_{m,k}^*R_{m,k}(\xi\ten\zeta)
			&=\sum_{s=0}^{\min\{m,k\}}
			\sum_{\substack{J\subseteq[m],\ |J|=m-s\\
					L\subseteq[k], |L|=s}}q^{s^2+i(J)+i(L)}
			(\xi_J\zeta_L)\ten |I|^{\otimes k}(\xi_{J^c}\zeta_{L^c})\\
			&=\sum_{s=0}^{\min\{m,k\}}
			\sum_{\substack{J\subseteq[m],\ |J|=m-s\\
					L\subseteq[k],\ |L|=s}}
			q^{s^2+i(J)+i(L)}
			\lambda^{\frac{\chi(\xi_{J^c})+\chi(\zeta_{L^c})}{2}}
			(\xi_J\zeta_L)\ten(\xi_{J^c}\zeta_{L^c}).
		\end{align*}
		Since $D_k^{-1}(\xi\ten\zeta)=\lambda^{-\chi(\zeta)/2}\xi\ten\zeta$
		and $\chi(\zeta)=\chi(\zeta_L)+\chi(\zeta_{L^c})$, it follows that
		\begin{align*}
			D_kR_{m,k}^*R_{m,k}D_k^{-1}(\xi\ten\zeta)&=
			\sum_{s=0}^{\min\{m,k\}}
			\sum_{\substack{J\subseteq[m],\ |J|=m-s\\
					L\subseteq[k],\ |L|=s}}
			q^{s^2+i(J)+i(L)}
			\lambda^{\frac{\chi(\xi_{J^c})-\chi(\zeta_L)}{2}}
			(\xi_J\zeta_L)\ten(\xi_{J^c}\zeta_{L^c})\\
			&=\sum_{s=0}^{\min\{m,k\}}
			q^{s^2}\sum_{\substack{J\subseteq[m],\ |J|=m-s\\
					L\subseteq[k],\ |L|=s}}
			q^{i(J)+i(L)}
			\lambda^{\frac{\chi(\xi_{J^c})-\chi(\zeta_L)}{2}}
			(\xi_J\zeta_L)\ten(\xi_{J^c}\zeta_{L^c}).
		\end{align*}
		A similar calculation would show that the last term above is equal to 
		\[\sum_{s=0}^{\min\{m,k\}}q^{s^2}
		(R_{m-s, s}\otimes R_{s, k-s})\circ\operatorname{Flip}_s
		\bigl(\id\ten|I|^{\otimes s}\ten(|I|^{-1})^{\otimes s}\ten\id\bigr)\circ(R_{m-s, s}^\ast\otimes R_{s,k-s}^\ast)(\xi\otimes\zeta).\]
		Since $\operatorname{Flip}_s$ is unitary and $\norm{R_{m-s,s}}\leq C_q$, and $\|I\|, \|I^{-1}\|\leq \lambda ^{-1/2}$,
		it follows that
		\[
		\norm{D_kR_{m,k}^*R_{m,k}D_k^{-1}}
		\leq C_q^2\sum_{s\geq0}|q|^{s^2}\lambda^{-s}:=K_{q,\lambda}<\infty.
		\]
		Consequently, for any  $\eta\in H_q^{\ten n}\cap\ran Q_R$,
		\[
		\norm{D_k R_{m,k}^*\eta}_q
		\leq \norm{D_k R_{m,k}^*R_{m,k}D_k^{-1}}\,
		\norm{D_k R_{m,k}^{-1}\eta}_q
		\leq K_{q,\lambda}B_q^{-1/2}\lambda^{-R/2}\norm{\eta}_q
		\]
		which finally yields
		\begin{align*}
			\|W(\eta)\|\leq C_q(n+1)K_{q,\lambda}B_q^{-1/2}\lambda^{-R/2}\norm{\eta}_q.
		\end{align*}
		This finishes the proof.
	\end{proof}
	
	\subsection{Norm estimates of balanced polynomials}
	We now restrict to the two-dimensional plane spanned by $e$ and $f$, and use the notation following \cite{KSW23,MS23}. Let $\mathcal P_{1,2}(n)$ be the set of all partitions of  $\{1,\ldots,n\}$ into blocks of size one or two. For a given partition $\pi\in \mathcal{P}_{1,2}(n)$, write
	$p(\pi)$ for the set of pairs and $s(\pi)$ for the set of singletons. Let $\mathrm{cross}(\pi)$ denote the number of crossings in $\pi$ in the following sense:
	\begin{align*}
		\mathrm{cross}(\pi)  = & \#\{((i,j),(k,l))\in p(\pi)\times p(\pi):\ i<k<j<l\}\\
		& + \#\{((i,j),k)\in p(\pi)\times s(\pi):i<k<j\}.
	\end{align*}
	Set $\xi_j=w_j\Omega$, and for $s(\pi)=\{j_1<\ldots<j_s\}$ write
	\[\delta_{p(\pi)}^w=\prod_{\underset{i<j}{(i,j)\in p(\pi)}}\varphi_q(w_iw_j)\quad \text{and}\quad \xi_{s(\pi)}^w=\xi_{j_1}\otimes\ldots\otimes \xi_{j_s}.\]

	\begin{lemma}(\cite[Lemma 3.10]{BKM23})\label{lem:formula for a word on omega}
		For $n\geq1$ and $w_1,\ldots,w_n\in \{W(e), W(e)^\ast\}$,
		\[
		w_1\cdots w_n\Omega
		=\sum_{\pi\in\mathcal P_{1,2}(n)}
		q^{\mathrm{cross}(\pi)}\,\delta_{p(\pi)}^w\,\xi_{s(\pi)}^w.
		\]
	\end{lemma}
	
	\begin{lemma}\label{lem:estimate on the sum of q powers cross}
		For any $q\in (-1,1)$ and $n\geq 1$ we have
		\begin{equation*}
			\sum_{\pi\in\mathcal P_{1,2}(n)}q^{\operatorname{cross}(\pi)}
			\le \left(1+\frac1{1-|q|}\right)^n.
		\end{equation*}
	\end{lemma}
	\begin{proof}
		Without loss of generality, assume that $q>0$.
		We sketch the argument.
		Let $z_n$ denote the sum on the left.
		Every partition of $\{1,\ldots,n\}$ is obtained from a partition $\sigma$ of $\{2,\ldots,n\}$ in one of two ways. One may add $1$ as a singleton, in which case the crossing number is unchanged, or pair $1$ with the $j$th singleton of $\sigma$, in which case the crossing number increases by $j-1$.
		This yields
		\begin{align*}
			z_n=\sum_{\sigma\in \mathcal{P}_{1,2}(\{2,\ldots, n\})}
			q^{\operatorname{cross}(\sigma)}
			\left(1+\sum_{j=1}^{|s(\sigma)|}q^{j-1}\right)
			\le \left(1+\frac1{1-q}\right)z_{n-1}.
		\end{align*}
		The lemma then follows by induction on $n$.
	\end{proof}

	As before, let $U$ be a nontrivial representation, with $e,f$  as above. The following result establishes the rapid-decay estimates for balanced polynomials required in Proposition \ref{prop:rapid_decay-property}.
	
	\begin{proposition}\label{prop:rd}
		Let $x_1,\ldots,x_d$ be any fixed balanced polynomials in ${W(e)},{W(e)}^*$.
		There are constants $C,L>0$  such that every polynomial $g$ of
		degree at most $m$ in $x_1,\ldots,x_d,x_1^*,\ldots,x_d^*$ satisfies
		\[
		\norm g\leq C(m+1)^{3/2}L^{\sqrt{m}}\norm{g\Omega}_q.
		\]
	\end{proposition}
	\begin{proof}
		\textit{Step 1}:
		We estimate the coefficients of $g$ as a polynomial  in ${W(e)}$ and ${W(e)}^*$. More precisely, write $g=\sum_wc_ww$, where $c_w\in \mathbb C$ and words $w$ consist of letters in $\{W(e), W(e)^\ast\}$. Let $N_0$ denote the maximum of the degrees of the generators $x_j$. Then $g$ has degree at most $N_0m$. Setting $N:=N_0m$, we show that there is a constant $C>0$ independent of $g$ and $m$ such that \begin{align}\label{eq:estimate of coeff of g}
			\sum_w|c_w|\leq C^{N}\|g\Omega\|_q.
		\end{align}
		
		To this end, we first express $g$ as a linear combination of the Appell polynomials \cite{Ans04}, and then estimate the modulus of the coefficients in this expansion.
		Write $c=W(e)$. Let $\mathcal{B}$ be the unital $\ast$-subalgebra generated by $c$ and $c^\ast$.
		Let $NC(n)$ be the set of non-crossing partitions of $[n]$. Recall that a partition is called {\em non-crossing} if
		there do not exist distinct blocks $B,C$ and indices
		$i<j<k<l$ with $i,k\in B$ and $j,l\in C$.
		For each $n\geq1$, let 		$\kappa_n:\mathcal{B}^n\to\mathbb C$ be the free cumulant with respect to the vacuum state $\varphi_q$, which is defined   recursively by 
		\[
		\varphi_q(w_1\cdots w_n)
		=\sum_{\pi\in NC(n)}
		\prod_{\{i_1<\cdots<i_k\}\in\pi}
		\kappa_k(w_{i_1},\ldots,w_{i_k}),
		\]
		for any $w_1,\ldots,w_n\in \mathcal{B}$ (\cite[Definition 11.3]{NS06}).
		The noncrossing-partition inversion formula gives		\begin{align}\label{eq:cumulatn bound}
			|\kappa_n(w_1,\ldots,w_n)|\leq 16^n \prod_{j=1}^n\|w_j\|
		\end{align}
		(see, for example \cite[Proposition 13.15]{NS06} for one variable, with the same proof for several variables).	
		For words $w=w_0w_1\ldots w_n$ in $\{{c},{c}^*\}$, define a polynomial $P_w$, called a free Appell polynomial, in $\{c, c^\ast\}$ recursively by $P_\emptyset=1$ and
		\begin{align}\label{eq:recursion formula for Pw}
			P_{w_0\ldots w_n}=w_0P_{w_1\ldots w_n}-\sum_{j=1}^n
			\kappa_{j+1}(w_0,w_1,\ldots,w_{j})P_{w_{j+1}\ldots w_n}.
		\end{align}
		The recursive formula above is precisely \cite[Proposition 3.12(b)]{Ans04}, although the Appell recursion  begins with $j=0$ in \cite{Ans04}. To see this, we claim that $\kappa_1(w_0)=\varphi_q(w_0)=0$.  For any word $w\ne\emptyset$, we have $\ph_q(P_w)=0$. Indeed, $\varphi_q(w_0)=0$, and by recursion,  we have
		\begin{align*}
			\varphi_q(P_{w_0w_1\ldots w_n})=\varphi_q(w_0P_{w_1\ldots w_n})-\kappa_{n+1}(w_0,\ldots, w_n).
		\end{align*}
		Then the claim will follow if we  prove  that 
		\[
		\varphi_q(aP_{w_1\cdots w_n})=\kappa_{n+1}(a,w_1,\ldots,w_n),\qquad a\in\mathcal{B}.
		\]
		We proceed by induction. For $n=0$, $\varphi_q(a)=\kappa_1(a).$  For $n\ge1$, apply the recursion \eqref{eq:recursion formula for Pw} to $P_{w_1\cdots w_n}$, multiply by $a$ on the left, and apply the induction hypothesis to obtain
		\[
		\begin{aligned}
			\varphi_q(aP_{w_1\cdots w_n})
			=\kappa_n(aw_1,w_2,\ldots,w_n)-\sum_{j=1}^{n}
			\kappa_j(w_1,\ldots,w_j)\,
			\kappa_{n-j+1}(a,w_{j+1},\ldots,w_n).
		\end{aligned}
		\]
		By using the cumulant formula with the product $aw_1$ in its first slot, the right hand side above is equal to $\kappa_{n+1}(a,w_1,\ldots,w_n)$ by \cite[Theorem 11.12]{NS06}. This proves the claim that $\varphi_q(P_w)=0$ for all $w\neq\emptyset.$
		
		Next, for any word $w=w_1\ldots w_n$ in $\{c, c^\ast\}$ and $w_0\in\{{c},{c}^*\}$, 
		the moment recursion and polynomial differentiation from \cite[Theorem 10]{Ans04} give
		\begin{align}\label{eq:formula for delta on appell}
			\delta_{w_0}P_w=\sum_{w=u w_0 v}P_u\ten P_v.
		\end{align}
		One can also verify this  directly by induction using the recursive formula \eqref{eq:recursion formula for Pw}.
		Consider the bounded operators $D_c, D_{c^\ast}$ on $L^2(M,\varphi_q)$ as defined in Corollary \ref{cor:boundedconjugates}.
		For any nonempty word $w=w_1\ldots w_n$ in $\{c,c^\ast\}$, the  formulas for $\delta_c P_w$  from \eqref{eq:formula for delta on appell}
		and the fact that
		$\ph(P_w)=0$ give
		\[
		D_{c}(P_{w_1\cdots w_n})=\sum_{\underset{w_j=c}{j=1}}^nP_{w_1\ldots w_{j-1}}\varphi_q(P_{w_{j+1}\ldots w_n})
		=1_{\{c=w_n\}}P_{w_1\cdots w_{n-1}}.
		\]   
		The family of Appell polynomials $\{P_w: w \mbox{ is a word in }\{c,c^\ast\}\}$ forms a linear basis for $\mathcal{B}$ (\cite[Proposition 12(a)]{Ans04}), so we can write $g=\sum_w\alpha_wP_w$ as a finite sum with scalars $\alpha_w$. The above computation for
		$w=w_1\cdots w_n$ yields $\alpha_w=\ip\Omega{D_{w_1}\cdots D_{w_n}g\Omega},$
		so that 
		\[|\alpha_w|\leq \left(\max\{1, \|D_c\|, \|D_{c^\ast }\|\}\right)^n\,\|g\Omega\|_q.\]
		For any polynomial $P\in\mathcal{B}$, let $|P|_{\mathrm{coeff}}$ denote the sum of absolute values of coefficients in $P$, that is, if $P=\sum_{u}\alpha_u u$ with $u$ words in $\{c, c^\ast\}$, then $|P|_{\mathrm{coeff}}=\sum_{w}|\alpha_w|$. In particular,  the coefficients of the Appell polynomials  for any word $w$ satisfy 	
		\[|P_w|_{\mathrm{coeff}}\leq C_0^{|w|},\qquad C_0=32\max\{1,\|c\|\}.\]  This follows from the recursion in \eqref{eq:recursion formula for Pw}  and the bound of cumulants in \eqref{eq:cumulatn bound}.
		Therefore, if $g=\sum_w\alpha_wP_w$ is a linear sum of $P_w$'s, then
		\begin{align*}
			|g|_{\mathrm{coeff}}\leq \sum_{w}|\alpha_w|\,|P_w|_{\mathrm{coeff}}\leq \sum_wC^{|w|}\|g\Omega\|_q\leq\sum_{j=0}^{N_0m} (2C)^{j}\,\|g\Omega\|\leq (2C)^{N_0m+1}\|g\Omega\|,
		\end{align*}
		where  $C=C_0\max\{1, \|D_c\|, \|D_{c^\ast}\|\}$.
		This proves \eqref{eq:estimate of coeff of g} required in  the first step.
		\medskip
		
		To clearly state the next step, we fix some notation.
		Using Lemma \ref{lem:formula for a word on omega},  we first express $g=\sum_wc_w w$ in terms of Wick operators as follows: 
		\begin{align*}
			g=\sum_{w}c_w  \sum_{\pi\in\mathcal P_{1,2}(|w|)}
			q^{\mathrm{cross}(\pi)}\delta_{p(\pi)}^w\,W(\xi_{s(\pi)}^w).
		\end{align*}
		For  $R>0$, let $g_{>R}$ be the sum of the Wick terms in the
		preceding expansion having $\mathrm{h}(\xi^w_{s(\pi)})>R$, that is
		\begin{align}\label{eq:formula for g greater R}
			g_{>R}
			=
			\sum_w c_w
			\sum_{\substack{\pi\in\mathcal P_{1,2}(|w|)\\
					\mathrm{h}(\xi^w_{s(\pi)})>R}}
			q^{\operatorname{cross}(\pi)}
			\delta_{p(\pi)}^w\,W(\xi^w_{s(\pi)}),
		\end{align}
		where  $h(\xi^w_{s(\pi)})$ is the height of the word $\xi^w_{s(\pi)}$ as defined in \eqref{eq:height of word}.
		
		\medskip
		\textit{Step 2:} We next apply Lemma \ref{lem:heightwick} to estimate separately the operator norms $\|g_{>R}\|$ and $\|g-g_{>R}\|$ in terms of the $L^2$-norm $\|g\Omega\|$ for large $R>0$. 
		
		
		We first estimate $\|g_{>R}\|$ for large $R$, to be defined later. For any word $w=w_1\ldots w_n$ with $w_j\in \{c,c^\ast\}$, set $\xi_j=w_j\Omega$ and $\xi^w=\xi_1\otimes\ldots\otimes\xi_n$. Note that $\xi_j\in \{e, \sqrt{\lambda}f\}$.
		Let $R_0$ be the maximum of $\mathrm{h}(\xi^w)$ for all words $w$ of $\{c,c^\ast\}$ appearing in the polynomials $x_1,\ldots,x_d, x_1^\ast,\ldots,x_d^\ast$. Since $w$ is a product of monomials of $x_j, x_j^\ast$'s and since $w$ is balanced by hypothesis, it follows that $\mathrm{h}(\xi^w)\leq R_0$.
		
		Fix a word $w=w_1\ldots w_n$ of $\{c,c^\ast\}$ appearing as a monomial in $g$. Also fix  a partition $\pi\in \mathcal{P}_{1,2}(n)$  with $\delta_{p(\pi)}^w\ne0$.  Note that $\|\xi_j\|\leq \max\{\|e\|, \sqrt{\lambda}\|f\|\}\leq1.$ Since $\delta_{p(\pi)}^w\neq 0$, we must have   $\varphi_q(w_iw_j)\neq 0$ for each  $(i,j)\in p(\pi)$, and because
		$\varphi_q(c^2)=\varphi_q((c^*)^2)=0$, this is equivalent to $w_i= w_j^\ast$. Hence $\chi(\xi_i)+\chi(\xi_j)=0$. This further implies that $\chi(\xi^w_{s(\pi)})=\chi(\xi^w)=0$, that  is, $\xi^w_{s(\pi)}$ is balanced and that $\varphi_q(w_iw_j)=\varphi_q(w_j^\ast w_j)=\|\xi_j\|^2\leq 1$. Hence $|\delta_{p(\pi)}^w|\leq 1$.

		For $0\le m\le n$, set
		\[
		t_m(\xi^w)=\sum_{i=1}^{m}\chi(\xi_i),
		\qquad
		t_m^\pi(\xi^w)=\sum_{\substack{i\in s(\pi)\\ i\le m}}\chi(\xi_i)
		\]
		that is, $t_m(\xi^w)$ sums the charge of each letter in the prefix $\xi_1\otimes\ldots \otimes\xi_m$, while $t_m^\pi(\xi^w)$ counts only the singletons in that prefix. 
		Then it follows from $\chi(\xi_i)+\chi(\xi_j)=0$ for each pair $(i,j)\in p(\pi)$ with $i<j\leq m$ that
		\[
		t_m(\xi^w)-t_m^\pi(\xi^w)
		=\sum_{\substack{(i,j)\in p(\pi)\\ i\le m<j}}\chi(\xi_i).
		\]
		Now suppose that $\mathrm{h}(\xi^w_{s(\pi)})>R_0$.  Put $k:=k_{w}=h(\xi_{s(\pi)}^w)-R_0>0$. We need to estimate the operator norm of each term in \eqref{eq:formula for g greater R}. We first claim that $\mathrm{cross}(\pi)\geq k(k+1)/2$.
		For each $1\le \ell\le k$, let
		\[
		m_\ell=\min\{m\in [n]:t_m^\pi(\xi^w)=R_0+\ell\},
		\]
		which exists because $\mathrm{h}(\xi^w_{s(\pi)})\geq R_0+\ell$.
		The positions $m_1,\ldots,m_k$ are distinct singletons, since
		$t_m^\pi(w)$ changes only at singleton positions.
		Since $t_{m_\ell}(\xi^w)\le h(\xi^w)\leq  R_0$, we have
		\[
		\sum_{\substack{(i,j)\in p(\pi)\\ i\le m_\ell<j}}\chi(\xi_i)=t_{m_\ell}(\xi^w)-t_{m_{\ell}}^\pi(\xi^w)
		=t_{m_\ell}(\xi^w)-(R_0+\ell)\le -\ell.
		\]
		Each summand in the left side is $\pm1$, so there are at least $\ell$ pairs in this sum, that is, $\#\{(i,j)\in p(\pi):i<m_\ell<j\}
		\geq \ell$.
		Moreover, $m_\ell$ is a singleton, so every such pair satisfies
		$i<m_\ell<j$. Therefore,
		\begin{equation}\label{eq:rd-cross-lower}
			\operatorname{cross}(\pi)
			\ge \sum_{\ell=1}^k
			\#\{(i,j)\in p(\pi):i<m_\ell<j\}
			\ge \sum_{\ell=1}^k\ell
			=\frac{k(k+1)}2.
		\end{equation}
		Next, if $s(\pi)=\{j_1<\ldots<j_{|s(\pi)|}\}$, then  Lemma \ref{lemma;some facts} yields 
		\[ \|\xi^w_{s(\pi)}\|_q=\|\ell_q(\xi_{j_1})\ldots\ell_q(\xi_{j_{|s(\pi)|}})\Omega\|\le  (1-|q|)^{-(|s(\pi)|)/2}\leq (1-|q|)^{-n/2}.\]
		As noted before, $\xi_{s(\pi)}^w$ is balanced, so $\xi^w_{s(\pi)}\in \mathcal{H}_U^{\otimes |s(\pi)|}\cap\mathrm{ran}Q_{h(\xi_{s(\pi)}^w)}$. Lemma \ref{lem:heightwick} then implies that
		\[
		\|W(\xi^w_{s(\pi)})\|
		\le C'(n+1)\lambda^{-\mathrm{h}(\xi^w_{s(\pi)})/2} (1-|q|)^{-n/2}= C'(n+1)\lambda^{-(k+R_0)/2} (1-|q|)^{-n/2}
		\]
		for some constant $C'>0$. 
		Since we also have $|\delta_{p(\pi)}^w|\leq 1$, and since $k> R-R_0$, it follows using \eqref{eq:rd-cross-lower} that
		\begin{align*}
			\|q^{\operatorname{cross}(\pi)}\delta^w_{p(\pi)}W(\xi_{s(\pi)})\bigr\|&\leq \sqrt{|q|}^{\mathrm{cross}}  C'(n+1)(1-|q|)^{-n/2}|q|^{k(k+1)/4}\lambda^{-(k+R_0)/2}\\
			&\leq 
			\sqrt{|q|}^{\mathrm{cross}}\, C'(n+1)(1-|q|)^{-n/2} \lambda^{-R_0/2} K_{q,\lambda}\,|q|^{(R-R_0)^2/8}.
		\end{align*}
		where $K_{q,\lambda}:=\sup_{j\geq0} |q|^{j^2/8}
		\left(\sqrt{|q|}\lambda^{-1}\right)^{j/2}<\infty$, and the last equality follows because
		\[
		|q|^{k(k+1)/4}\lambda^{-k/2}
		=|q|^{k^2/8}
		\left[|q|^{k^2/8}\left(\sqrt{|q|}\lambda^{-1}\right)^{k/2}\right]
		\leq K_{q,\lambda} r^{(R-R_0)^2/8}.
		\]
		Then  \eqref{eq:estimate of coeff of g}, \eqref{eq:formula for g greater R}, \eqref{eq:rd-cross-lower} and Lemma \ref{lem:estimate on the sum of q powers cross} with $N=N_0m$ yield
		\begin{align*}
			\|g_{>R}\|&\le \sum_w|c_w|  C'(|w|+1)(1-|q|)^{-|w|/2}\lambda^{-R_0/2}K_{q,\lambda}|q|^{(R-R_0)^2/8}\sum_{\pi\in \mathcal{P}_{1,2}(|w|)}\sqrt{|q|}^{\mathrm{cross}(\pi)}
			\\
			&\leq C''(N+1)(1-|q|)^{-N/2}|q|^{(R-R_0)^2/8}\left(1+\frac{1}{1-\sqrt{|q|}}\right)^{N}\sum_w|c_w| \\
			&\leq C''(N+1)(1-|q|)^{-N/2}\lambda^{-R_0/2}|q|^{(R-R_0)^2/8}\left(1+\frac{1}{1-\sqrt{|q|}}\right)^{N} C^N\|g\Omega\|_q\\
			&\leq C'' (N+1)C_1^{N}|q|^{(R-R_0)^2/8}\|g\Omega\|_q.
		\end{align*}
		where $C_1=(1-|q|)^{-1/2}\left(1+\frac{1}{1-\sqrt{|q|}}\right)C$ and $C''=C'\lambda^{-R_0/2}K_{q,\lambda}$ are constants independent of $N$.
		Now choose $R$ sufficiently large so that $R=R_0+O(\sqrt{N})$ and satisfies $C_1^N|q|^{(R-R_0)^2/8}\leq C_1^{\sqrt{N}}$, e.g. $R=R_0+A \sqrt{N}$ for some constant $A>0$ satisfying $A^2 N\log|q|= 8(\sqrt{N}-N)\log C_1$. This choice of $R$ finally yields
		\begin{align*}
			\|g_{>R}\|\leq C'' (N+1)C_1^{\sqrt{N}}\|g\Omega\|_q.
		\end{align*}
		We next estimate the operator norm $\|g-g_{>R}\|$ for this particular $R$. 
		For any word $w$ and a partition $\pi$,  we have already observed that $\xi_{s(\pi)}^w$ is balanced. Thus, if $h(\xi_{s(\pi)}^w)\leq R$, then $\xi_{s(\pi)}^w\in Q_R$ by the definition of $Q_R$. In particular, $(g-g_{>R})\Omega\in \ran Q_R$. Since each monomial in $g$ is of length at most $N=N_0m$, it follows that $(g-g_{>R})\Omega\in \oplus_{j=0}^N\mathcal{H}_U^{\otimes j}\cap \ran Q_R$. Write $(g-g_{>R})\Omega=\sum_{j=0}^N\eta_j$ for some $\eta_j\in \mathcal{H}_U^{\otimes j}\cap\mathrm{ran}Q_R$. Lemma~\ref{lem:heightwick} yields $\|W(\eta_j)\|\leq C_2(j+1)\lambda^{-R/2}\|\eta_j\|_2$ for some constant $C_2>0$. This together with the Cauchy-Schwarz inequality implies that
		\begin{align*}
			\|g-g_{>R}\|\leq \sum_{j=0}^N\|W(\eta_j)\|\leq C_2\lambda^{-R/2}(N+1)^{3/2}\left(\sum_{j=0}^N\|\eta_j\|_2^2\right)^{1/2}=C_2\lambda^{-R/2}(N+1)^{3/2}\|(g-g_{>R})\Omega\|_q.
		\end{align*}
		Combining $\|(g-g_{>R})\Omega\|_q\leq \|g\Omega\|_q+\|g_{>R}\|$ with the estimate for $\|g_{>R}\|$ and the fact that $R=R_0+O(\sqrt{N})$,
		we obtain
		\begin{align*}
			\|g\| \leq C'''(N+1)^{3/2}C_2^{\sqrt{N}}\|g\Omega\|_q
		\end{align*}
		for some $C_2, C'''>0$.  Writing the above constant in terms of $N_0m$ will yield the required constant for the inequality.
		This finishes the proof.
	\end{proof}

	\section{Simplicity of \texorpdfstring{$q$}{q}-Araki-Woods algebras with two generators}\label{section:simplicity-q-Araki-Woods}
	Let $\dim_{\mathbb R} \mathcal H_{\mathbb R}=2$
	and let $U$ be a nontrivial orthogonal representation on $\mathcal{H}_{\mathbb R}$.
	In this section, we apply the estimates and general theory developed in the preceding sections to  the $q$-Araki–Woods $C^*$-algebra $A_q(\mathbb R^2, U)$.
	As before, there are unit vectors $e,f$ of $\mathcal H_U$ such that
	$Ie=\lambda^{1/2} f$, $If=\lambda^{-1/2}e$, where
	$0<\lambda<1$. 
	
	\begin{proposition}\label{thm:simplicity-rotation-plane}
		Let $q\in (-1,1)$, and let $U$ be a nontrivial orthogonal representation
		on $\mathbb R^2$.
		Put $\mathcal A_q=A_q({\mathbb R}^2,U)$ and let
		$\varphi_q$ be the vacuum state. Then $\mathcal A_q$ and $\mathcal A_q^{\varphi_q}$ have the Dixmier averaging property and
		$\tau={\varphi_q}_{|_{\mathcal A_q^{\varphi_q}}}$ is the unique tracial state on $\mathcal A_q^{\varphi_q}$.
	\end{proposition}
	
	\begin{proof}
		Let $e$ be the unit vector in $\mathbb C^2$ with respect to $\|\cdot\|_U$ such that $A_Ue =\lambda^{-1} e$ for $0<\lambda<1$.
		As mentioned above, the set $G=\{W(e),W(e)^*\}$ consists of modular
		eigenoperators of $\sigma^{\varphi_q}$ and generates $\mathcal A_q$ as a $C^\ast$-algebra. It is also known that the conjugate variables for elements in $G$ exist and are Lipschitz \cite[Corollaries 3.10 and 3.13]{Kum23}. Proposition~\ref{prop:fullgap} is then applicable, so
		there are self-adjoint balanced polynomials $s_1,s_2$ in
		$\{W(e)^\ast W(e),W(e)W(e)^*\}$ such that
		\[
		\|f-\tau(f)1\|_2
		\le C\sum_{j=1}^2\|[s_j,f]\|_2,
		\qquad f\in L^2(M^{\varphi_q},\tau).
		\]
		Fix a balanced word $g$ in $G$ and $m\geq 1$ and recall the set $V_m(g)$ from Proposition~\ref{prop:rapid_decay-property}. Every element of $V_m(g)$ is a
		polynomial of degree at most $m+1$ in the fixed balanced family
		$\{s_1,s_2,g,g^*\}$. Proposition~\ref{prop:rd} therefore gives
		$\rho_m(g)\le C_g(m+2)^{3/2}L_g^{\sqrt{m+1}}$
		for some constants $C_g,L_g\ge1$.
		In particular, $\log\rho_m(g)\leq\log C_g+\frac{3}{2}\log(m+2)+(\log L_g)\sqrt{m+1}=o(m)$. 
		By Proposition~\ref{prop:rapid_decay-property},
		$\mathcal A_q^{\varphi_q}$ has the Dixmier averaging property. Hence $\mathcal A_q$ has the Dixmier averaging property by Theorem \ref{thm:sim_centra_implies_simple}.
	\end{proof}

	\subsection*{An explicit element in the centralizer} We now identify an explicit element $B\in \mathcal{A}_q^{\varphi_q}$ whose unitary conjugates asymptotically average to a scalar. To this end, we establish several technical computations for certain Wick operators. The following extends the uniform boundedness result in \cite[Lemma 3.9]{BMRW23} and \cite[Lemma 3.6]{KW24}. In the sequel, 
	we adopt the following convention:  
	\[\xi_1\cdots\xi_n =\xi_1\otimes\cdots\otimes\xi_n, \qquad\xi^n=\xi^{\otimes n},\qquad
	\xi^0=\Omega,\qquad \mbox{ for } \xi,\xi_1,\ldots,\xi_n\in\mathcal H_U.
	\]

	\begin{lemma}\label{lemma: bounded sequence}
		Let $m,n\ge 0$. The Wick formula for $W(f^ne^m)$ is
		\begin{equation}\label{eq:blockedwick}
			\begin{split}
				W(f^n e^m)={}&\sum_{i=0}^m\sum_{j=0}^n
				\binom mi_q\binom nj_q q^{(m-i)j}\lambda^{(i-j)/2}\ell_q(f)^{n-j}\ell_q(e)^{(m-i)}
				\ell_q(e)^{*j}\ell_q(f)^{*i}.
			\end{split}
		\end{equation} 
		The operators $\left((1-q)^{(m+n)/2}W(f^ne^{m})\right)_{m\geq n\geq 1}$ are uniformly bounded.
	\end{lemma}
	\begin{proof}
		We show this by induction on $n$,
		simultaneously for all $m\geq0$. For $n=0$, using the definition of $W(e)$ and $Ie=\bar e=\sqrt{\lambda}f$, \cite[Lemma 3.2]{BMRW23}) gives
		\begin{equation}\label{eq:powerexpansion}
			W(e^m) = W(e)^m=\sum_{i=0}^m\binom mi_q\lambda^{i/2}
			\ell_q(e)^{m-i}\ell_q(f)^{*i}.
		\end{equation} 
		Similarly, when $m = 0$, we have
		\begin{equation*}
			W(f^n) = W(f)^n=\sum_{j=0}^n\binom nj_q\lambda^{-j/2}
			\ell_q(f)^{n-j}\ell_q(e)^{*j}.
		\end{equation*}
		Suppose that \eqref{eq:blockedwick} holds for a given $n\geq 0$ and every $m\geq0$, and fix $m\geq1$.
		We first establish the following equation
		\begin{equation}\label{eq:7.4}
			W(f^{n+1}e^m) = W(f)W(f^ne^m) - \lambda^{-1/2}q^n[m]_qW(f^ne^{m-1}).
		\end{equation}
		Indeed, since $W(f)=\ell_q(f)+\lambda^{-1/2}\ell_q^\ast(e)$,  it follows from Lemma \ref{lemma;some facts}(i) 
		that 
		$$W(f)(f^ne^m) 
		= f^{n+1}e^m + \lambda^{-1/2} q^n[m]_qf^ne^{m-1},$$
		which yields \eqref{eq:7.4} by the separating property of the vacuum vector $\Omega$. Next, we calculate $W(f)W(f^ne^m)$.  Lemma \ref{lemma;some facts} and   the $q$-commutation relations imply, by induction,  that for $k\ge 0$ and $r\ge 1$ 
		\begin{equation}\label{eq: 7.5}
			\ell_q(e)^*\ell_q(f)^k = q^k\ell_q(f)^k\ell_q(e)^* \quad\text{and}\quad \ell_q(e)^*\ell_q(e)^r = [r]_q\ell_q(e)^{r-1} + q^r \ell_q(e)^r\ell_q(e)^*.
		\end{equation}
		Set 
		\[
		\alpha_{n,m}(i,j)
		=\binom{m}{i}_q\binom{n}{j}_q
		q^{(m-i)j}\lambda^{(i-j)/2}
		\qquad(0\leq i\leq m,\ 0\leq j\leq n),
		\]
		and set $\alpha_{n,m}(i,j)=0$ outside this range.
		By the induction hypothesis, we  apply $W(f)=\ell_q(f)+\lambda^{-1/2}\ell_q^\ast(e)$ to the expression in  \eqref{eq:blockedwick} and then use \eqref{eq: 7.5} to obtain
		\begin{align*}
			W(f)W(f^ne^m) = &  \sum_{i=0}^m\sum_{j=0}^n \alpha_{n,m}(i,j)\ell_q(f)^{n+1-j}\ell_q(e)^{(m-i)}
			\ell_q(e)^{*j}\ell_q(f)^{*i}\\
			& + \sum_{i=0}^{m-1}\sum_{j=0}^n  \lambda^{-1/2}q^{n-j}[m-i]_q \alpha_{n,m}(i,j)\ell_q(f)^{n-j}\ell_q(e)^{m-i-1}
			\ell_q(e)^{*j}\ell_q(f)^{*i}\\
			& + \sum_{i=0}^m\sum_{j=0}^n \lambda^{-1/2}q^{n-j+m-i}\alpha_{n,m}(i,j)\ell_q(f)^{n-j}\ell_q(e)^{(m-i)}
			\ell_q(e)^{*(j+1)}\ell_q(f)^{*i}\\
			& =: \mathrm{I} +\mathrm{II}+\mathrm{III}. 
		\end{align*}
		Using  \eqref{eq:q-fractional relations},
		we obtain $q^{n-j}[m-i]_q\alpha_{n,m}(i,j)
		=q^n[m]_q\alpha_{n,m-1}(i,j)$. Therefore, $$\mathrm{II} = \lambda^{-1/2}q^n[m]_qW(f^ne^{m-1}).$$
		Moreover, the coefficient of $\ell_q(f)^{n+1-j}\ell_q(e)^{(m-i)}
		\ell_q(e)^{*(j)}\ell_q(f)^{*i}$ in $\mathrm{I}+ \mathrm{III}$ 
		is 
		\begin{align*}
			\alpha_{n,m}(i,j) + \lambda^{-1/2}q^{n+1-j+m-i}\alpha_{n,m}(i,j-1)& =  \binom{m}{i}_q q^{(m-i)j}\lambda^{(i-j)/2}
			\left(
			\binom{n}{j}_q
			+q^{n+1-j}\binom{n}{j-1}_q
			\right)\\
			& = \binom{m}{i}_q\binom{n+1}{j}_q
			q^{(m-i)j}\lambda^{(i-j)/2}\\
			& = \alpha_{n+1,m}(i,j),
		\end{align*}
		where the second equation uses \eqref{eq;Pascal identity}.
		Combining this calculation with \eqref{eq:7.4} completes the induction.
		
		It remains to prove the boundedness assertion. Take $m\ge n\ge1$. By Lemma \ref{lemma;some facts}(ii) 
		the powers $(1-q)^{n/2}\ell_q(e)^n$ and $(1-q)^{n/2}\ell_q(f)^n$ are uniformly bounded, which implies that
		\[
		(1-q)^{(m+n)/2}\|W(f^ne^m)\|
		\leq C_q^6
		\sum_{j=0}^n\sum_{i=0}^m
		|q|^{(m-i)j}\lambda^{(i-j)/2}\leq C_q^6
		\sum_{i,j=0}^m
		|q|^{(m-i)j}\lambda^{(i-j)/2}.
		\]
		The right hand side is uniformly bounded in $m$ as proved in \cite[Lemma 3.9]{BMRW23} and \cite[Lemma 3.6]{KW24}. This proves the assertion.
	\end{proof}
	
	
	\begin{lemma}\label{lemma;Positive element B}
		For every $n\ge 0$, we have
		\begin{equation}\label{eq:Bn}
			(1-q)^nW(e)^{*n}W(e)^n
			=d_n^2\sum_{j=0}^n
			\frac{\lambda^{j/2}(1-q)^jW(f^j e^j)}{d_j^2d_{n-j}}.
		\end{equation}
		Moreover, the sequence 
		$ \bigl((1-q)^n W(e)^{*n}W(e)^n\bigr)_{n\geq 1} $
		converges in norm to an element 
		$B\in \mathcal A_q$ satisfying 
		$ \varphi_q(B)=d_\infty $ and $ (1-q)W(e)^* B W(e)=B. $ 
		In particular, $B\in {\mathcal{A}_q}^{\ph_q}_+$.
	\end{lemma}
	\begin{proof}
		From \eqref{eq:powerexpansion}, we have $W(e)^n\Omega=e^n$. So by the same formula and the equality $\ell_q(e)^{*(n-j)}e^n
		=\frac{[n]_q!}{[j]_q!}\,e^{j}$, we obtain
		\begin{align*}
			W(e)^{\ast n}W(e)^n\Omega&= \sum_{j=0}^n\binom nj_q\lambda^{j/2}\ell_q(f)^{j}
			\ell_q(e)^{*(n-j)}e^n\\
			&=\sum_{j=0}^n\binom nj_q\frac{[n]_q!}{[j]_q!}\lambda^{j/2}f^{j}e^{j}\\
			&=\frac{d_n^2}{(1-q)^n}\sum_{j=0}^n\frac{(1-q)^{j}}{d_j^2d_{n-j}}\lambda^{j/2}W(f^{j}e^{j})\Omega,
		\end{align*}
		which proves the required equation. 
		By 
		Lemma \ref{lemma: bounded sequence} and $0<\lambda<1$, the norms of the operator terms on the right-hand side are bounded by a summable sequence. Thus
		\begin{equation*}
			B=\lim_{n\to\infty}(1-q)^nW(e)^{*n}W(e)^n=d_\infty\sum_{m\geq0}
			\frac{\lambda^{m/2}(1-q)^mW(f^m e^m)}{d_m^2}
			\ \in {\mathcal A_q}_+^{\ph_q}.
		\end{equation*}
		The other assertions are obvious now.
	\end{proof}
	
	The following result follows from Proposition \ref{thm:simplicity-rotation-plane}.
	\begin{corollary}\label{cor:full}
		Let $B$ be the positive element in Lemma \ref{lemma;Positive element B}. 
		Then $d_\infty1\in \overline{\mathrm{co}}^{\|\cdot\|}\left\{uBu^*:u\in \mathcal{U}(\mathcal A_q^{\varphi_q})\right\}.$
	\end{corollary}

	\section{Selflessness of \texorpdfstring{$q$}{q}-Araki-Woods \texorpdfstring{$C^*$}{C*}-algebras}\label{sec:selfless_q-Araki-Woods}
	We are now ready to prove our main result Theorem \ref{main-theorem}. We actually prove a stronger result: whenever $U$ is nontrivial, the $q$-Araki-Woods $C^\ast$-algebra $A_q(\mathcal{H}_{\mathbb R}, U)$ is  Toeplitz selfless in the sense of \cite[Definition 5.1, Proposition 5.12]{GJKEPR26}.
	
	We first recall these notions. Let $\mathcal{A}$ be a unital $C^\ast$-algebra and let $\Phi:\mathcal{A}\to\mathcal{A}$ be a completely positive map. 
	The pair $(\mathcal A,\Phi)$ is called {\em Toeplitz
		selfless} if there exists an ultrafilter $\omega$ and $v\in \mathcal{A}^\omega$ such that $v^\ast av=\Phi(a)$ and $\Phi^\omega(avv^\ast a^\ast)=0$ for all $a\in \mathcal{A}$. Here $\Phi^\omega$ is the completely positive map on $\mathcal{A}^\omega$ given by $(a_i)^\omega\mapsto(\Phi(a_i))^\omega$. If   $\varphi$ is a state on $\mathcal{A}$ whose GNS representation is faithful, then   Toeplitz selflessness of $(\mathcal{A}, \varphi(\cdot)1_{\mathcal{A}})$ implies complete selflessness of $(\mathcal{A},\varphi)$, which in
	turn implies selflessness
	\cite[Lemmas~3.10 and~7.2, Proposition~5.12, and Theorem~7.3]{GJKEPR26}.
	

	Let
	$U:\mathbb R\to\mathcal O(\mathcal H_{\mathbb R})$ be a nontrivial orthogonal representation with $\dim_{\mathbb R}\mathcal H_{\mathbb R}\geq 2$.  Put
	$\mathcal A_q=A_q(\mathcal H_{\mathbb R},U)$  and $M_q=\mathcal A_q''$. Put $\mathcal C_q=C^\ast(1,\ell_q(\xi); \xi\in \mathcal{H}_U)$. Then $\mathcal A_q\subset\mathcal C_q$.  Denote the ambient vacuum state by ${\varphi_q}$.
	
	\subsection{Finite-Dimensional Setup}
	We first assume that $\mathcal H_{\mathbb R}$ is finite-dimensional. The first step is to show that the map $x\mapsto \varphi_q(B^{1/2}xB^{1/2})B$ is Toeplitz selfless, where $B$ is as in Lemma \ref{lemma;Positive element B}. This, combined with the Dixmier averaging property in Corollary \ref{cor:full}, then shows that $(\mathcal A_q, \varphi_q)$ is Toeplitz selfless. 
	
	Since $U$ is nontrivial, we find mutually orthogonal unit vectors
	$e,f\in \mathcal H_U$ such that 
	$A_Ue=\lambda^{-1}e$, $ A_Uf=\lambda f
	$ and $Ie=\lambda^{1/2}f$ for some
	$0<\lambda<1$.
	Extend $\{e,f\}$ to an orthonormal basis, say $\Lambda$,  of $\mathcal H_U$.  
	For a word $\xi=\xi_1\cdots \xi_k$ in the basis $\Lambda$, write 
	\[
	|\xi|=k,\qquad \ell_q(\emptyset)=1,\qquad \ell_q(\xi)=\ell_q(\xi_1)\cdots\ell_q(\xi_k).
	\]
	The $q$-commutation relations imply that the linear span of the operators $\ell_q(\xi)\ell_q(\eta)^*$ is norm dense in $\mathcal C_q$.
	For $m\geq1$, define states $\tau_m:\mathcal C_q\to\mathbb C$ by
	\[
	\tau_m(x)=\frac{1}{[m]_q!}
	\ip{f^m}{xf^m}_q=\frac{1}{[m]_q!}
	\ip{\ell_q(f)^m\Omega}{x\ell_q(f)^m\Omega}_q.
	\]
	Let $\tau_f$ be a limit point of the sequence $\{\tau_m\}_{m\geq 1}$ in the weak-$*$ topology, which exists because $\{\tau_m\}_{m\geq 1}$ is uniformly bounded.  
	For $m\geq k$,   a simple calculation yields\ $\tau_m\bigl(\ell_q(f)^k\ell_q(f)^{*k}\bigr)
	=\frac{[m]_q!}{[m-k]_q!}$ which converges to $ (1-q)^{-k}$ as $m\to\infty$.
	All other basis monomials $\ell_q(\xi)\ell_q(\eta)^*$ have value zero. Thus we get
	\begin{align}\label{eq:formula for tau_f}
		\tau_f\bigl(\ell_q(\xi)\ell_q(\eta)^*\bigr)=
		\begin{cases}
			(1-q)^{-k},&\xi=\eta=f^k\text{ for some }k\geq0,\\
			0,&\text{otherwise}.
		\end{cases}
	\end{align}
	Since the linear span of words formed from the basis vectors $\Lambda$ is dense in $\mathcal{F}_q(\mathcal{H}_U)$, we obtain the following orthogonal decomposition of the $q$-Fock space 
	\begin{align}\label{eq:decomposition of Fock space}
		\mathcal F_q(\mathcal H_U)=\bigoplus_{j\geq0}
		\ell_q(f)^j\ker\ell_q(f)^*.
	\end{align}
	Indeed, if $i>j\geq0$ and $\zeta_1,\zeta_2\in \ker\ell_q(f)^\ast$, then the $q$-commutation relation gives $\ell_q(f)^{\ast i}\ell_q(f)^j\zeta_2=\alpha \ell_q(f)^{\ast(i-j)}\zeta_2$ for some scalar $\alpha$. Hence $\langle \ell_q(f)^i\zeta_1, \ell_q(f)^j\zeta_2\rangle_q =\langle\zeta_1, \alpha \ell_q(f)^{\ast(i-j)}\zeta_2\rangle_q =0,$ 
	which proves the asserted orthogonality.
	
	\medskip
	
	The following lemma generalizes \cite[Lemma 3.4]{BMRW23}.
	\begin{lemma}\label{lem:comv of Tnkl}
		For any $k,l\geq0$, the sequence  \[S_n^{k,l}=(1-q)^{n+(k+l)/2}
		\ell_q(f)^{*n}\ell_q(f)^k\ell_q(f)^{*l}\ell_q(f)^n\] converges in norm to $ L_{k-l}$ as $n\to\infty$, where $L_0 := \oplus_{j\geq0}\frac{d_\infty}{d_j}\, \operatorname{id}_{\ell_q(f)^j\ker\ell_q(f)^*}$ is the diagonal operator relative to the orthogonal decomposition of $\mathcal F_q(\mathcal H_U)$ in \eqref{eq:decomposition of Fock space}, 
		and for $s\in\mathbb Z$,
		\[
		L_s :=
		\begin{cases}
			(1-q)^{s/2}L_0\ell_q(f)^s,&s\geq0,\\
			(1-q)^{-s/2}\ell_q(f)^{*(-s)}L_0,&s<0.
		\end{cases}
		\]
	\end{lemma}
	\begin{proof}
		For the case $k=l=0$, the sequence $S_n^{0,0}:=(1-q)^n\ell_q(f)^{\ast n}\ell_q(f)^n$ converges  in norm to $L_0$ by  \cite[Lemma 3.4]{BMRW23}. 
		Fix $k,l\geq0$. It follows from     \eqref{eq:normal-ordering}
		that $S_n^{k,l}$ converges in norm to some operator as $n\to\infty$. Therefore, it suffices to show that $S_n^{k,l}$ converges to $L_{k-l}$ pointwise on each summand $\ell_q(f)^j\ker\ell_q(f)^\ast$, $j\geq 0$.  
		Let $\zeta\in \ker\ell_q(f)^\ast$
		and write $\eta_j=\ell_q(f)^j\zeta$. For large $n$, the $q$-commutation relations give 
		\[\ell_q(f)^{\ast l}\ell_q(f)^n\eta_j=\frac{[n+j]_q!}{[n+j-l]_q!}\ell_q(f)^{n-l}\eta_j.\] 
		Consequently,
		\begin{align*}
			S_{n}^{k,l}\eta_j=(1-q)^{n+(k+l)/2}\frac{[n+j]_q!}{[n+j-l]_q!}\ell_q(f)^{\ast n}\ell_q(f)^{n+k-l}\eta_j.
		\end{align*}
		If $k-l\geq0$, then
		\[S_{n}^{k,l}\eta_j= (1-q)^{(k-l)/2}\frac{d_{n+j}}{d_{n+j-l}}S_n^{0,0}\ell_q(f)^{k-l}\eta_j\to (1-q)^{(k-l)/2}L_0\ell_q(f)^{k-l}\eta_j=L_{k-l}\eta_j;\] 
		and if  $k-l<0$, then
		\begin{align*}
			S_{n}^{k,l}\eta_j&=(1-q)^{(l-k)/2}\frac{d_{n+j}}{d_{n+j-l}}\ell_q(f)^{\ast(l-k)}S_{n+k-l}^{0,0}\eta_j\to (1-q)^{(l-k)/2}\ell_q(f)^{\ast(l-k)}L_0\eta_j=L_{k-l}\eta_j
		\end{align*}
		as $n\to\infty$. This completes the proof.
	\end{proof}
	
	The following is the key step in proving Toeplitz selflessness
	of $(\mathcal A_q,\varphi_q)$. It establishes this property first for a
	certain completely positive map on the larger $C^\ast$-algebra
	$\mathcal C_q$.
	
	\begin{proposition}\label{prop:estimate for Toeplitz selflessness}
		For $n\geq 1$, define $T_n:=(1-q)^{3n/2}\ell_q(f)^nW(e)^{2n}\in\mathcal C_q$. Let $B$ be the positive element in Lemma \ref{lemma;Positive element B} and  $x\in\mathcal C_q$. Then the following hold.
		\begin{enumerate}
			\item[(i)]  $\norm{T_n^*xT_n-d_\infty\tau_f(x)B}\longrightarrow0$.
			\item[(ii)]  $\tau_f(xT_nT_n^*x^*)\longrightarrow0.$
		\end{enumerate}
		In particular,   $(\mathcal{C}_q, \Phi_B)$ is Toeplitz selfless, where 
		$\Phi_B:\mathcal C_q\to\mathcal C_q$ is  the completely positive map  $\Phi_B(x)=d_\infty\tau_f(x)B$ for $x\in\mathcal C_q$. 
	\end{proposition}
	\begin{proof}
		It follows from the uniform bounds in Lemma \ref{lemma: bounded sequence} and Lemma \ref{lemma;some facts}(ii) that $(T_n)_{n\geq1}$ is uniformly bounded. The two assertions then imply that $(\mathcal{C}_q,\Phi_B)$ is Toeplitz selfless. Indeed, if  $\omega$ is a free ultrafilter on $\mathbb N$, then  $T=(T_n)^{\omega}\in \mathcal C_q^{\omega}$ and we have $T^\ast xT=\Phi_B(x)$ and $\Phi_B^{\omega}(xTT^\ast x^\ast)=0$ for all $x\in\mathcal C_q$. So we now prove the two statements.
		
		(i) It suffices to prove the assertion for $x=\ell_q(\xi)\ell_q(\eta)^\ast$, where $\xi=\xi_1\ldots\xi_k$ and $\eta=\eta_1\ldots\eta_l$ are words in the basis $\Lambda$ of $\mathcal{H}_U$. We have
		\begin{align*}
			T_n^\ast xT_n=(1-q)^{3n}W(e)^{\ast 2n}\ell_q(f)^{\ast n}\ell_q(\xi)\ell_q(\eta)^\ast\ell_q(f)^{n}W(e)^{2n}.
		\end{align*}
		Suppose first that $\eta$ contains a basis vector different
		from $f$. In this case,  $\tau_f(x)=0$. So we need to show that $\|T_n^\ast xT_n\|\to 0$.  Let us calculate $\ell_q(\eta)^\ast\ell_q(f)^n$. Write $\eta=f^r g\eta'$, where $g\neq f$, $0\leq r<l$, and $\eta'$ is a word in the basis vectors. For  $n\geq r$,  \eqref{eq:normal-ordering}
		gives
		\begin{align*}
			\ell_q(f)^{*r}\ell_q(f)^n
			=
			\sum_{k=0}^r
			q^{(r-k)(n-k)}
			\binom{r}{k}_q\binom{n}{k}_q[k]_q!\,
			\ell_q(f)^{n-k}\ell_q(f)^{*(r-k)}.
		\end{align*}
		Since $g$ and $f$ are orthogonal, we also have
		$\ell_q(g)^*\ell_q(f)^{n-k}
		=q^{n-k}\ell_q(f)^{n-k}\ell_q(g)^*.$
		Consequently,
		\begin{align*}
			\ell_q(\eta)^*\ell_q(f)^n
			=\sum_{k=0}^r
			q^{(r-k+1)(n-k)}
			\binom{r}{k}_q\binom{n}{k}_q[k]_q!\,
			\ell_q(\eta')^*\ell_q(f)^{n-k}
			\ell_q(g)^*\ell_q(f)^{*(r-k)}.
		\end{align*}
		For $0\leq k\leq r$, we have $(r-k+1)(n-k)\geq n-r$. Thus the uniform bound in 
		Lemma \ref{lemma;some facts}(ii) gives
		\begin{equation}\label{eq:Proof-in-Lemma7.2}
			(1-q)^{n/2}
			\|\ell_q(\eta)^*\ell_q(f)^n\|
			\leq C_q^{6+l}\sum_{k=0}^r |q|^{n-r}\|\eta'\|_q\,\|g\|_q = C_\eta|q|^{n-r},
		\end{equation}
		where $C_\eta>0$  is a constant independent of $n$.
		By (i) and the boundedness of $\bigl((1-q)^nW(e)^{2n}\bigr)_n$, we obtain
		\[
		\|T_n^*xT_n\|
		\leq
		\|T_n\|\,\|\ell_q(\xi)\|\,
		\|\ell_q(\eta)^*T_n\|
		\leq \|T_n\| \,\|\ell_q(\xi)\|\,C_{\eta}|q|^{n-r}\to 0
		\]
		as $n\to\infty$.
		The same argument applied to $x^*$ proves the assertion
		when $\xi$ contains a basis vector different from $f$.
		Therefore, we now assume that $x=\ell_q(f)^k\ell_q(f)^{\ast l}$ for some $k,l\geq0$. Then
		\begin{align*}
			T_n^\ast xT_n=(1-q)^{2n-(k+l)/2}W(e)^{\ast 2n}S_n^{k,l}W(e)^{2n},
		\end{align*}
		where $S_n^{k,l}$ is as in Lemma \ref{lem:comv of Tnkl}. By the same lemma, we have
		\begin{align}\label{eq:convergence of T_nxT}
			\|T_n^\ast xT_n-(1-q)^{2n-(k+l)/2}W(e)^{\ast 2n}L_{k-l}W(e)^{2n}\|\to 0,\qquad n\to\infty.
		\end{align} 
		For $s\in\mathbb Z$, write $D_s=L_s-\delta_{s,0}d_\infty1 $. We shall show that $(1-q)^{2n-(k+l)/2}W(e)^{\ast 2n}D_sW(e)^{2n}\to 0$ as $n\to\infty$, which will then prove that $\|T_n^\ast xT_n-d_\infty \tau_f(x)B\|\to 0$ using the formula in \eqref{eq:formula for tau_f}, convergence in \eqref{eq:convergence of T_nxT} and Lemma \ref{lemma;Positive element B}.
		For any $m\geq1$ and $s\geq 0$, 
		\begin{align*}
			W(e)^{\ast m}D_s W(e)^m=\sum_{i,j=0}^m\binom mj_q\binom mi_q \lambda^{m-(i+j)/2}\ell_q(f)^{m-j}
			\ell_q(e)^{\ast j}D_s\ell_q(e)^{i}\ell_q(f)^{\ast (m-i)}.
		\end{align*}
		We now estimate the terms $\ell_q(e)^{\ast j}D_s\ell_q(e)^{i}$.
		Let $P$ be the orthogonal projection onto
		$\ker\ell_q(f)^*$.  Since the range of $\ell_q(f)$ is closed and equals  the range of $1-P$, $\ell_q(f)$ has a bounded inverse, denoted by $X$, from its range to $\mathcal{F}_q(H_U)$. Using
		$\ell_q(e)^{\ast j}\ell_q(f)=q^{j}\ell_q(f)\ell_q(e)^{\ast j}$, we obtain
		\begin{align*}
			\|\ell_q(e)^{\ast j}(1-P)\|&=\|\ell_q(e)^{\ast j}\ell_q(f)X(1-P)\|=|q|^{j}\|\ell_q(f)\ell_q(e)^{\ast j}X(1-P)\|
			\leq C'|q|^j\,\|\ell_q(e)^{\ast j}\|
		\end{align*}
		for some constant $C'>0$ depending on $q$.  Since $PL_0=L_0P=d_\infty P$, it follows $PD_sP=0$ for every $s\in\mathbb Z$.  This implies that
		\begin{align*}
			\|\ell_q(e)^{\ast j}D_s\ell_q(e)^i\|&\leq \|\ell_q(e)^{\ast j}(1-P)D_s\ell_q(e)^i\|+\|\ell_q(e)^{\ast j}PD_s(1-P)\ell_q(e)^i\|\\
			&\leq C'\left(|q|^i+|q|^j\right)\|D_s\|\,\|\ell_q(e)^{\ast j}\|\;\|\ell_q(e)^i\|.
		\end{align*}
		Thus, the uniform bound in Lemma \ref{lemma;some facts}(ii) implies
		\begin{align*}
			(1-q)^m  \|  W(e)^{\ast m}D_s W(e)^{m}\|&\leq C'C_q^6\|D_s\|\, \sum_{i,j=0}^{m}\lambda^{m-(i+j)/2}(|q|^i+|q|^j)\| \\
			&\leq 2C'C_q^6\|D_s\|\, (m+1)^2(\lambda^{m/2}+|q|^{m/2}).
		\end{align*}
		We obtain the last inequality by  dividing the sum according to whether $i$ and $j$ are at most or at least $m/2$.  
		This shows that $(1-q)^{m+s}\|W(e)^{\ast m}D_sW(e)^{m}\|\to 0$ as $m\to\infty$. Taking $m=2n$ and letting $n\to\infty$ yields the claim.
		
			(ii) As above, it is sufficient to show the limit for  $x=\ell_q(\xi)\ell_q(\eta)^*$ for some words $\xi,\eta$ with letters in the basis. Let $(\pi_f,\mathcal H_f,\Omega_f)$ be the GNS representation of $\tau_f$ so that $\tau_f(x T_nT_n^\ast x^\ast)=\|\pi_f(T_n^\ast x^\ast)\Omega_f\|^2$. 
			If $\xi$ contains a letter different from $f$, then
			$\pi_f(x^*)\Omega_f=0$ and then $\tau_f(x T_nT_n^* x^\ast)=0$. If $\eta$ contains a letter different from $f$ then part (i) above shows that $\|\pi_f(T_n^\ast x^\ast)\Omega_f\|^2\leq\|T_n^\ast x^\ast \|^2=\|T_n^\ast x^\ast xT_n\|\to 0$ as $n\to\infty$.
			We may therefore assume that $\xi=f^k$ and $\eta=f^l$ for $k,l\geq 0$.  
			For $j\in\mathbb Z$,
			put
			\[
			\zeta_j=
			\begin{cases}
				(1-q)^{j/2}\pi_f(\ell_q(f)^j)\Omega_f,&j\geq0,\\
				(1-q)^{-j/2}\pi_f(\ell_q(f)^{*(-j)})\Omega_f,&j<0.
			\end{cases}
			\]
			A simple calculation  gives an orthonormal family $(\zeta_j)_{j\in\mathbb Z}$. Indeed, for $j,l\geq 0$,
			\begin{align*}
				\langle\zeta_j,\zeta_l\rangle=(1-q)^{(j+l)/2}\lim_{m\to\infty}\frac{1}{[m]_q!}\langle \ell_q(f)^j f^m,\ell_q(f)^\ell f^m\rangle_q=(1-q)^{(j+l)/2}\delta_{j,l}\lim_{m\to\infty}\frac{1}{[m]_q!}\|f^{l+m}\|_q^2=\delta_{j,l},
			\end{align*}
			and for $j\geq 0, l<0$,
			\begin{align*}
				\langle\zeta_j,\zeta_l\rangle =(1-q)^{(j-l)/2}\lim_{m\to\infty}\frac{1}{[m]_q!}\langle \ell_q(f)^j f^m,\ell_q(f)^{\ast(-\ell)} f^m\rangle_q 
				=0.
			\end{align*}   
			Next we claim that
			\[
			\sqrt{1-q}\,\pi_f(\ell_q(f))\zeta_j=\zeta_{j+1},\qquad j\in\mathbb Z.
			\]
			This is obvious for $j\geq 0$. When $j<0$, for any words $\mu=\mu_1\ldots\mu_n$ and $\zeta=\zeta_1\ldots\zeta_m$ in basis vectors, we have
			\begin{align*}
				\langle \pi_f(\ell_q(\mu)\ell_q(\zeta)^\ast)\Omega_f, \pi_f(\ell_q(f)\ell_q(f)^{\ast(-j)})\Omega_f\rangle&=\lim_{m\to\infty}\frac{1}{[m]_q!}\langle \ell_q(\mu)\ell_q(\zeta)f^m, \ell_q(f)\ell_q(f)^{\ast(-j)}f^m\rangle_q\\
				&=\lim_{m\to\infty}\frac{1}{[m+j+1]_q!}\langle \ell_q(\mu)\ell_q(\zeta)f^m, \ell_q(f)^{\ast(-j-1)}f^m\rangle_q\\
				&=\langle \pi_f(\ell_q(\mu)\ell_q(\zeta)^\ast)\Omega_f, \pi_f\ell_q(f)^{\ast(-j-1)})\Omega_f\rangle.
			\end{align*}
			A similar calculation gives
			\begin{align*}
				\sqrt{1-q}\,\pi_f(\ell_q(f)^*)\zeta_j=\zeta_{j-1},\qquad j\in\mathbb Z.
			\end{align*}
			Moreover, we have $\pi_f(\ell_q(g)^*)\zeta_j=0$ and $\pi_f(W(e)^{*2n})\zeta_j=\lambda^n\pi_f(\ell_q(f)^{2n})\zeta_j$, for every basis vector $g\ne f$, $j\in \mathbb Z$ and $n\ge0$.  For large $n$, the above calculation yields 
			\begin{align*}
				\pi_f( T_n^\ast x^\ast)\Omega_f&=(1-q)^{3n/2}\pi_f(W(e)^{\ast 2n}\ell_q(f)^{\ast n}\ell_q(f)^l\ell_q(f)^{\ast k})\Omega_f\\
				&=(1-q)^{n-(l+k)/2}\pi_f(W(e)^{\ast 2n})\zeta_{l-k-n}\\
				&= \lambda^{n}(1-q)^{n-(l+k)/2}\pi_f(\ell_q(f)^{2n})\zeta_{l-k-n}\\
				&= \lambda^{n}(1-q)^{-(l+k)/2}\zeta_{l-k+n}.
			\end{align*}
			Since $\zeta_j$ is uniformly bounded and $\lambda^n\to 0$, it follows that $\tau_f(xT_nT_n^\ast x^\ast)=\|\pi_f(T_n^\ast x^\ast)\Omega_f\|^2\to 0$ as $n\to\infty$.
		\end{proof}


		\begin{corollary}\label{cor: rest-A-Toe-selfless}
			With the  map $\Phi_B$ as in Proposition \ref{prop:estimate for Toeplitz selflessness}, we have   $\Phi_B(x)=\varphi_q(xB)B$ for every $x\in \mathcal A_q$. Hence  $\Phi_B(\mathcal{A}_q)\subseteq\mathcal{A}_q$. Moreover,
			$(\mathcal{A}_q,{\Phi_B}_{|_{\mathcal A_q}})$ is Toeplitz selfless. 
		\end{corollary}  
		\begin{proof}

			For $x\in\mathcal A_q$ and $m\geq0$, it follows from the KMS condition and formulas
			$W(e)^{*m}\Omega=\lambda^{m/2}f^m$ and  $\sigma_i^{\varphi_q}(W(e)^m)=\lambda^{-m}W(e)^m$  that
			\begin{align*}
				\varphi_q\bigl(xW(e)^{*m}W(e)^m\bigr)=\varphi_q(\sigma_i^{\varphi_q}(W(e)^m)x(W(e)^{\ast m}))
				=\ip{f^m}{xf^m}_q.
			\end{align*}
			By Lemma \ref{lemma;Positive element B}, we obtain
			$$\varphi_q(xB) = \lim_{m\to\infty} \varphi_q\bigl(x(1-q)^mW(e)^{*m}W(e)^m\bigr) = \lim_{m\to\infty} d_m\tau_m(x) = d_\infty \tau_f(x),$$
			and then $\Phi_B(x)=d_\infty\tau_f(x)B = \varphi_q(xB)B=\varphi_q(B^{1/2}xB^{1/2})B$.
			
			Let $\omega$ be a free ultrafilter on $\mathbb{N}$. Consider the bounded sequence $(T_n)_{n\geq1}$ from Proposition \ref{prop:estimate for Toeplitz selflessness}.  We claim that  $T=(T_n)^{\omega}\in \mathcal A_q^{\omega}$. This will show that $\Phi_B|_{\mathcal A_q}$ is Toeplitz selfless by Proposition  \ref{prop:estimate for Toeplitz selflessness}.
			
			To this end, set $S_n=(1-q)^{3n/2}W(f^ne^{2n})\in \mathcal A_q.$   We shall show that $\|T_n-S_n\|\to0$ as $n\to\infty$, which will then imply $(T_n)^\omega=(S_n)^\omega\in \mathcal{A}_q^\omega$. 
			Taking $m=2n$ in \eqref{eq:blockedwick} and applying \eqref{eq:powerexpansion}, we obtain
			\begin{align*}
				T_n-S_n= (1-q)^{3n/2} \sum_{i=0}^{2n}\sum_{j=1}^n \binom{2n}{i}_q \binom{n}{j}_q q^{(2n-i)j}\lambda^{(i-j)/2}  \ell_q(f)^{n-j}\ell_q(e)^{2n-i}\ell_q(e)^{*j}\ell_q(f)^{*i}
			\end{align*}
			which using  Lemma \ref{lemma;some facts}(ii), (vi) gives
			\begin{align*}
				\|T_n-S_n\|
				\leq C_q^6 \sum_{i=0}^{2n}\sum_{j=1}^n |q|^{(2n-i)j}\lambda^{(i-j)/2}. 
			\end{align*}
			If $q=0$, then $\|T_n-S_n\|\leq \sum_{j=1}^n\lambda^{(2n-j)/2}\leq  \frac{\sqrt{\lambda}^n}{1-\sqrt{\lambda}}\to 0$ as $n\to\infty$. Now assume that $q\neq0$.
			Fix $j$, and set $a_i=|q|^{(2n-i)j}\lambda^{(i-j)/2}$ for $1\le i\le 2n$. Then $\frac{a_{i+1}}{a_i}=\frac{\sqrt{\lambda}}{|q|^j}$ for each $0\leq i\leq 2n$. If $\frac{\sqrt{\lambda}}{|q|^j}<1$  then $(a_{i})_{i=0}^{2n}$ decreases, and hence $a_i\leq a_0=|q|^{2nj}\lambda^{-j/2}\leq |q|^{2n}\lambda^{-1/2}$. On the other hand, if  $\frac{\sqrt{\lambda}}{|q|^j}\geq1$, then $a_i\leq a_{2n}=\lambda^{(2n-j)/2}\leq \lambda^{n/2}$. Thus, in either case, we obtain $a_i\leq \max\{\lambda^{-1/2}|q|^{2n},\lambda^{n/2}\}$. It follows that
			\begin{align*}
				\|T_n-S_n\|\leq C_q^6 n(2n+1)\max\{\lambda^{-1/2}|q|^{2n},\lambda^{n/2}\},
			\end{align*}
			which tends to $0$ as $n\to\infty$. This finishes the proof.
		\end{proof}

		\subsection{Proof of the main result} We are now ready to establish complete selflessness of $q$-Araki-Woods $C^\ast$-algebras in full generality.
		
		\begin{proof}[Proof of Theorem \ref{main-theorem}]
			Write $\mathcal A_q=A_q(\mathcal H_{\mathbb R},U)$ and $M_q=\mathcal A_q''$.   If $U_t=1$ for all $t$, then $\mathcal A_q$ is completely selfless by  \cite[Theorem 6.3]{KP26}. Therefore, assume that $U$ is nontrivial.
			
			First assume that $\dim_{\mathbb R} \mathcal H_{\mathbb R}<\infty$.        
			Let $B\in \mathcal A_q$ be the operator in  Lemma \ref{lemma;Positive element B}.
			By Corollary \ref{cor:full}, there are unitaries
			$u_1,\ldots,u_k\in \mathcal A_q^{\varphi_q}$ and $t_1,\ldots,t_k\geq0$ with
			$\sum_{i=1}^k t_i=1$, such that
			$h:=\sum_{i=1}^k t_iu_iBu_i^*$
			satisfies $\norm{h-d_\infty 1}<d_\infty /2$.
			In particular, $h$ is positive and invertible in $\mathcal A_q^{\varphi_q}$.
			Put $b_i=\sqrt{t_i}\,h^{-1/2}u_i\in \mathcal A_q^{\varphi_q}$.
			Therefore, $\sum_{i=1}^k b_iBb_i^*=1$ and using the fact that $b_j$ belongs to the centralizer
			one obtains
			\[
			\sum_{i,j}^k b_i\,\Phi_B(b_j^*ab_j)\,b_i^*=\sum_{i,j=1}^k\varphi_q(b_j^\ast ab_jB)b_iBb_i^\ast
			=\left(\sum_{j=1}^k\varphi_q(a b_j B b_j^*)\right)
			\left(\sum_{i=1}^k b_i B b_i^*\right)
			=\varphi_q(a)1,
			\]    
			where $\Phi_B$ is the map in Corollary \ref{cor: rest-A-Toe-selfless}. Since $\Phi_B$ is Toeplitz selfless by the same corollary, it follows from \cite[Theorem 5.11(i)]{GJKEPR26} that $(\mathcal A_q,\varphi_q)$ is Toeplitz selfless, and hence completely selfless by \cite[Theorem~13]{Oza25}.
			
			Now we allow arbitrary $\mathcal H_{\mathbb R}$, possibly nonseparable.
			For every finite-dimensional real subspace
			$\mathcal {K}_{\mathbb R}\subseteq \mathcal H_{\mathbb R}$ with dimension at least $2$,
			write $\mathcal{K}_{\mathbb C}=\mathcal{K}_{\mathbb R}+i\mathcal{K}_{\mathbb R}$ and write $P_{\mathcal K}$ for
			the complex orthogonal projection
			$P_{\mathcal{K}}:\mathcal H_{\mathbb C}\to \mathcal {K}_{\mathbb C}$.
			Set
			\[
			G=\frac{2A_U}{1+A_U} \quad\text{and}\quad G_{\mathcal K} = P_{\mathcal K}G_{|_{\mathcal K_{\mathbb C}}}.
			\]
			Since $G$ and $2-G$ are nonsingular positive self-adjoint operators, it follows that $G_{\mathcal K}$ and $2-G_{\mathcal K}$ are injective operators on $\mathcal K_{\mathbb C}$. Finite dimensionality of $\mathcal K_{\mathbb C}$ implies that $G_{\mathcal K}$ and $2-G_{\mathcal K}$ are positive and invertible.

			The relation $IA_UI=A_U^{-1}$ gives $IGI=2-G$. Since $\mathcal{K}_{\mathbb C}$ is invariant under $I$, and  $I^2$ is the identity, we have $P_{\mathcal K}I=IP_{\mathcal K}$. This shows that $IG_{\mathcal{K}}I=2-G_{\mathcal {K}}$ on $\mathcal{K}_{\mathbb C}$.
			Define $A_{\mathcal{K}}:\mathcal{K}_{\mathbb C}\to \mathcal{K}_{\mathbb C}$ by $A_{\mathcal{K}}=G_{\mathcal{K}}(2-G_{\mathcal{K}})^{-1}$. Then $A_{\mathcal K}$ is positive and invertible and  $IA_{\mathcal K}I=A_{\mathcal K}^{-1}$ holds. Hence $IA_{\mathcal K}^{it}I=A_{\mathcal K}^{it}$ so that $A_{\mathcal K}^{it}(\mathcal K_{\mathbb R})\subset \mathcal K_{\mathbb R}$ for all $t\in\mathbb R$. This induces an orthogonal representation $U_{\mathcal{K}}:\mathbb R\to\mathcal O(\mathcal{K}_{\mathbb R})$ by $U_{\mathcal{K}}(t)=A_{\mathcal{K}}^{it}|_{\mathcal{K}_{\mathbb R}}$. We remark that the exact formula for $A_{\mathcal K}$ in terms of $A_U$ is
			\[
			A_{\mathcal K}=\left(P_{\mathcal K}(1+A_U)^{-1}\big|_{{\mathcal K}_{\mathbb C}}\right)^{-1}
			-\mathbf 1_{{\mathcal K}_{\mathbb C}},
			\]
			which is precisely equal to ${A_U}{|_{{\mathcal K}_{\mathbb C}}}$ whenever ${\mathcal K}_{\mathbb C}$ is invariant under $A_U$. Consider the $q$-Araki-Woods $C^\ast$-algebra $A_q({\mathcal K}_{\mathbb R}, U_{\mathcal K})$ associated to the representation $U_{\mathcal K}$. We denote by $s_q^{\mathcal K}(\xi)$ the field operators on the Fock space $\mathcal{F}_q({\mathcal K}_{U_{\mathcal K}})$ which generate $A_q(\mathcal K_{\mathbb R}, U_{\mathcal K})$. Note that the deformed scalar-product on $\mathcal K_{\mathbb C}$ is exactly given by
			\[\langle\xi, \eta\rangle_{U_{\mathcal K}}=\langle\xi, G_{\mathcal K}\eta\rangle_{{\mathcal K}_{\mathbb C}}=\langle\xi, G\eta\rangle_{\mathcal H_{\mathbb C}}=\langle\xi,\eta\rangle_U,\qquad \forall\xi,\eta\in \mathcal K_{\mathbb C}.\]
			The $q$-Wick pairing formula therefore gives equality of all joint vacuum moments of the  operators $s_q^{\mathcal{K}}(\xi), \xi\in \mathcal K_{\mathbb R}$ in $A_q(\mathcal K_{\mathbb R}, U_{\mathcal K})$ and the ambient operators $s_q(\xi)$ in $A_q(\mathcal H_{\mathbb R}, U)$. This induces an isometric unital vacuum state preserving $\ast$-homomorphism
			$\Theta_{\mathcal{K}}:A_q(\mathcal{K}_{\mathbb R},U_{\mathcal{K}})\to
			A_q(\mathcal H_{\mathbb R},U)$ by $ \Theta_{\mathcal{K}}(s_q^{\mathcal{K}}(\xi))=s_q(\xi)$ for $\xi\in \mathcal{K}_{\mathbb R}$
			(See \cite[Theorem 4.11]{NS06}). By the above observation, the $C^\ast$-algebra $A_q(\mathcal K_{\mathbb R}, U_{\mathcal K})$ is Toeplitz selfless if $U_{\mathcal{K}}$ is nontrivial, in which case its range $C^\ast(1, s_q(\xi); \xi\in \mathcal K_{\mathbb R})$ under $\Theta_{\mathcal K}$ is also Toeplitz selfless.
			
			Now since $U$ is nontrivial, choose a finite-dimensional real subspace $\mathcal{K}_0$ such that $G_{\mathcal{K}_0}\neq 1$. Then consider the subalgebras $\{\Theta_{\mathcal{K}}({A}_q({\mathcal{K}, U_{\mathcal{K}}})); \mathcal{K}\supseteq \mathcal{K}_0\}$  forming an upward directed family associated to  finite-dimensional real subspaces $\mathcal{K}$ of $\mathcal{H}_{\mathbb R}$ such that their union is norm dense in $A_q(\mathcal H_{\mathbb R},U)$. Since Toeplitz  selflessness passes to this directed norm closure by \cite[Proposition 5.12]{GJKEPR26}, this proves the theorem for all $U$ and $\mathcal H_{\mathbb R}$, including nonseparable spaces. 
			
			Finally, selflessness implies simplicity by \cite[Theorem~3.1]{Rob25}.
			If $U$ is nontrivial, then $\varphi_q$ is nontracial, hence  $A_q(\mathcal H_{\mathbb R},U)$ is purely infinite by
			\cite[Theorem~3.1]{Rob25}.
		\end{proof}
		
		\section*{Acknowledgment}
		M.~K. was partially supported by the Research Foundation - Flanders (FWO) through a postdoctoral fellowship (1221025N). X.-Q.~L. was partially supported by the China Scholarship Council (Grant No. 202306200023). S.~W. was partially supported by the NSF of China (No. 12671160, No. W2611005, No. W2441002, No.12031004).
		Part of this work was completed during a one-month visit by M.~K. and X.-Q.~L. to Harbin Institute of Technology in July--August 2026, and during their participation in the operator algebras trimester programme at IHP, where they collaborated primarily in September--October 2026. They gratefully acknowledge Harbin Institute of Technology and the organisers of the IHP programme for their warm hospitality and financial support.
		
		GPT-6 Astra was used  to proofread and edit
		the manuscript, and identify potential gaps in
		the proofs. The main calculations in Sections~\ref{sec:operator_norm_estimate}--\ref{sec:selfless_q-Araki-Woods} were developed through iterative
		exchanges with the model. The authors conceived and worked out the mathematical ideas behind
		Sections~\ref{sec:derivations on free fisher infomration}--\ref{sec:simplicity of C-algebras with eigenoperators}.
		The authors take full
		responsibility for the mathematical arguments and conclusions in the
		final manuscript.

	\end{document}